 \documentclass[10pt,reqno]{article}

\usepackage{amscd,amssymb}
\usepackage[mathscr]{eucal}
\usepackage[all]{xy}
\usepackage{mathrsfs}
\usepackage{xypic}
\usepackage{amsfonts}
\usepackage{amsmath}
\usepackage{amsthm}
\usepackage{latexsym}
\usepackage{eufrak}
\usepackage{verbatim}

\theoremstyle{plain}
\newtheorem{theorem}[equation]{Theorem}
\newtheorem{corollary}[equation]{Corollary}
\newtheorem{prop}[equation]{Proposition}
\newtheorem{lemma}[equation]{Lemma}

\theoremstyle{definition}
\newtheorem{defn}[equation]{Definition}

\newtheorem{conj}[equation]{Conjecture}

\newtheorem{assumption}[equation]{Basic Assumptions}
\newtheorem{claim}[equation]{Claim}

\theoremstyle{remark}
\newtheorem{notation}[equation]{Notation}
\newtheorem{rem}[equation]{Remark}

\makeatletter
\renewcommand{\subsection}{\@startsection{subsection}{2}{0pt}{-3ex
plus -1ex minus -0.2ex}{-2mm plus -0pt minus
-2pt}{\normalfont\bfseries}} \makeatother

\numberwithin{equation}{subsection}

\newcommand{\f}[1]{\mathfrak{#1}}
\newcommand{\scr}[1]{\mathscr{#1}}
\renewcommand{\cal}[1]{\mathcal{#1}}

\newcommand{\erem}{\hfill$\lozenge$\end{rem}}

\newcommand{\idot}{_{\:\raisebox{1pt}{\text{\circle*{1.5}}}}}
\newcommand{\hdot}{^{\:\raisebox{3pt}{\text{\circle*{1.5}}}}}

\DeclareMathOperator{\charp}{\mathsf{char}}

\DeclareMathOperator{\Aut}{\mathrm{Aut}}
\DeclareMathOperator{\Hom}{\mathrm{Hom}}

\DeclareMathOperator{\End}{\mathrm{End}}
\DeclareMathOperator{\eend}{{\mathcal{E}^{\!}\text{\it nd}}}
\DeclareMathOperator{\rk}{\mathrm{rk}}
\DeclareMathOperator{\Spec}{\mathrm{Spec}}
\DeclareMathOperator{\sym}{{\mathrm{Sym}}}
\DeclareMathOperator{\Ker}{\mathrm{Ker}}
\DeclareMathOperator{\im}{\mathrm{Im}}
\DeclareMathOperator{\ad}{\mathrm{ad}}
\DeclareMathOperator{\Ad}{\mathrm{Ad}}
\DeclareMathOperator{\Lie}{\mathrm{Lie}}
\DeclareMathOperator{\Tr}{\mathrm{tr}}
\DeclareMathOperator{\ind}{\mathrm{Ind}}
\DeclareMathOperator{\gr}{\mathrm{gr}}
\DeclareMathOperator{\dlog}{\mathrm{dlog}}
\DeclareMathOperator{\diag}{\mathrm{diag}}
\DeclareMathOperator{\rees}{{\boldsymbol{\mathcal{R}ees}}}

\newcommand{\dis}{{\displaystyle}}
\newcommand{\op}{\operatorname}

\newcommand{\vi}{${\sf {(i)}}\;$}
\newcommand{\vii}{${\sf {(ii)}}\;$}
\newcommand{\viii}{${\sf {(iii)}}\;$}
\newcommand{\iv}{${\sf {(iv)}}\;$}

\newcommand{\sset}{\subset}
\newcommand{\sseq}{\subseteq}
\newcommand{\sminus}{\smallsetminus}
\newcommand{\into}{\,\hookrightarrow\,}
\newcommand{\map}{\longrightarrow}
\newcommand{\too}{\,\,\longrightarrow\,\,}
\newcommand{\tooo}{{\;{-\!\!\!-\!\!\!-\!\!\!-\!\!\!\longrightarrow}\;}}
\newcommand{\isto}{{\underset{{}^{\sim}}{\;{-\!\!\!-\!\!\!-\!\!\!-\!\!\!\longrightarrow}\;}}}

\newcommand{\iso}{\stackrel{_\sim}\too}
\newcommand{\mto}{\longmapsto}
\newcommand{\onto}{\,\twoheadrightarrow\,}
\newcommand{\cd}{\!\cdot\!}
\newcommand{\dash}{\mbox{-}}
\newcommand{\pb}{\noindent$\bullet\quad$\parbox[t]{140mm}}
\def\hp{\hphantom{x}}

\newcommand{\inn}{^{\pm 1}}
\newcommand{\muu}{{\mu_{_\U}}}
\newcommand{\mua}{{\mu_{_{\op{alg}}}}}
\newcommand{\bzu}{{\mathbf{z}}_{_\U}}
\newcommand{\bzd}{{\mathbf{z}}_{_\D}}
\newcommand{\bzdr}{{\mathbf{z}}_{_{\rr\D}}}
\newcommand{\bzur}{{\mathbf{z}}_{_{\rr\U}}}

\newcommand{\BX}{{\mathbb{X}}^*}
\newcommand{\ke}{\k[\EE]}

\newcommand{\DD}{{\mathsf{D}}}
\newcommand{\Gm}{{\mathbb{G}}_{\mathbf{m}}}
\newcommand{\Th}{\Theta}
\newcommand{\TT}{{\mathcal T}}
\newcommand{\vol}{{\mathsf{vol}}}

\newcommand{\TTT}{{\mathbb{T}}}
\newcommand{\Z}{{\mathbb{Z}}}
\newcommand{\bz}{{\mathbf{z}}}

\newcommand{\oo}{{\mathcal O}}
\newcommand{\beq}{\begin{equation}\label}
\newcommand{\eeq}{\end{equation}}
\renewcommand{\mid}{\enspace\big|\enspace}
\newcommand{\ct}{{\stackrel{_\circ}{T^*}}X}
\newcommand{\G}{\Gamma}
\newcommand{\ffs}{{\mathbf{s}}}
\newcommand{\si}{\sigma}

\newcommand{\wh}{\widehat}
\newcommand{\whi}{\widehat{I} }

\newcommand{\RF}{\op{RF}}
\newcommand{\Hilb}{{\operatorname{Hilb}^n{\mathbb{A}}^2}}
\newcommand{\Hilgm}{{\operatorname{Hilb}^n(\BC^*\times\BC^*)}}
\newcommand{\RG}{{\operatorname{R\Gamma}}}
\newcommand{\z}{{\mathfrak{Z}}}
\newcommand{\QG}{\BQ^{\text{good}}}
\newcommand{\la}{\lambda}
\newcommand{\A}{{\mathsf{A}}}
\newcommand{\loc}{{{\scr L}oc}}
\newcommand{\az}{{\scr A}}
\newcommand{\B}{{\mathsf B}}
\newcommand{\ah}{{\scr H}}
\newcommand{\af}{{\mathbb{A}}}

\newcommand{\TS}{\widetilde{S}}
\newcommand{\rr}{{\boldsymbol{\mathcal{R}}}}
\newcommand{\ti}{\tilde{i}}

\newcommand{\irr}{{{\op{Irr}}(S_n)}}
\newcommand{\CY}{{\mathcal{Y}}}
\newcommand{\gre}{{$p\text{\em -graded}\enspace$}}
\newcommand{\EE}{{\mathfrak{G}}}

\newcommand{\eer}{{\EE^\circ}}
\newcommand{\Ups}{\Upsilon}

\newcommand{\tors}{_{\op{torsion}}}
\newcommand{\M}{{\mathsf{M}}}
\newcommand{\bv}{{\mathbf{v}}}
\newcommand{\Parab}{P}

\newcommand{\gc}{{\g^\circ}}
\newcommand{\VC}{V^\circ}
\newcommand{\triv}{{\mathsf{triv}}}
\newcommand{\sign}{{\mathsf{sign}}}
\newcommand{\mm}{{\scr M}}

\newcommand{\Si}{{\boldsymbol{\Pi}}}

\newcommand{\mss}{{\mathfrak{M}}^{\operatorname{s}}_0}

\newcommand{\rf}{\rr\CF}

\newcommand{\I}{{\mathcal I}}
\newcommand{\J}{{\scr J}}
\newcommand{\UU}{{\widetilde{U}}}

\newcommand{\ee}{{\scr E}}
\newcommand{\inv}{^{-1}}

\newcommand{\D}{{\scr D}}
\newcommand{\La}{\Lambda}
\newcommand{\e}{{\mathsf{e}}}

\newcommand{\g}{{\f g}}

\newcommand{\bvp}{{\boldsymbol{\varpi}}}

\newcommand{\cm}{{\mathsf{L}}}

\renewcommand{\aa}{{\f a}}
\newcommand{\p}{\f p}
\newcommand{\gl}{{\mathfrak{gl}}}

\newcommand{\h}{{\mathfrak{h}}}
\newcommand{\hreg}{{\mathfrak{h}}^{\operatorname{reg}}}

\newcommand{\tw}{^{(1)}}
\newcommand{\tp}{^{[p]}}
\newcommand{\tpp}{^{[1]}}

\newcommand{\atw}{{\mathfrak{a}}^{1,*}}
\newcommand{\ttw}{{T^{1,*}_{\as(\chi)}}}
\newcommand{\ttwp}{{T^{1,*}_{\as(\psi)}}}
\newcommand{\ttt}{^{1,*}}
\newcommand{\grs}{{\f g}^{\operatorname{rs}}}
\newcommand{\BP}{{\mathbb{P}}}
\newcommand{\hh}{{\mathsf{H}}}
\newcommand{\tf}{\tilde{f}}
\newcommand{\ohh}{{\overline{\mathsf{H}}}}
\newcommand{\HH}{{H\!H}}
\newcommand{\ehe}{{\mathsf{e}\mathsf{H}_c\mathsf{e}}}
\newcommand{\heh}{{\mathsf{H}\mathsf{e}\mathsf{H}}}
\newcommand{\id}{{\operatorname{Id}}}

\newcommand{\GL}{{\operatorname{GL}}}

\newcommand{\pr}{{\operatorname{pr}}}

\newcommand{\upi}{{\uparrow J}}
\newcommand{\U}{{\mathcal{U}}}
\newcommand{\Ug}{{\mathcal{U}\f g}}
\newcommand{\Ua}{{\mathcal{U}\f a}}
\newcommand{\uchi}{{{\mathfrak{u}}_\chi(\f a)}}
\newcommand{\upsi}{{{\mathfrak{u}}_\psi(\f a)}}
\renewcommand{\Up}{{\mathcal{U}\f p}}
\newcommand{\ZZ}{{\mathsf{Z}}}

\renewcommand{\k}{\Bbbk}
\newcommand{\vc}{\check{v}}

\def\ccirc{{{}_{^{\,^\circ}}}}

\newcommand{\nc}{\newcommand}
\nc{\on}{\operatorname}

\nc{\BA}{{\mathbb{A}}}
\nc{\BC}{{\mathbb{C}}}
\nc{\BM}{{\mathbb{M}}}
\nc{\BN}{{\mathbb{N}}}
\nc{\BR}{{\mathbf{R}}}
\nc{\BZ}{{\mathbb{Z}}}
\nc{\BS}{{\mathbb{S}}}
\nc{\BQ}{{\mathbb{Q}}}
\nc{\BF}{{\mathbb{F}}}

\nc{\CA}{{\mathcal{A}}}
\nc{\CB}{{\mathcal{B}}}
\nc{\CC}{{\mathcal{C}}}
\nc{\CE}{{\mathcal{E}}}
\nc{\CF}{{\mathcal{F}}}
\nc{\CG}{{\mathcal{G}}}
\nc{\CH}{{\mathcal{H}}}
\nc{\CL}{{\mathcal{L}}}
\nc{\CM}{{\mathcal{M}}}
\nc{\CN}{{\mathcal{N}}}
\nc{\CO}{{\mathcal{O}}}
\nc{\CP}{{\mathcal{P}}}
\nc{\CQ}{{\mathcal{Q}}}
\nc{\CR}{{\mathcal{R}}}
\nc{\CT}{{\mathcal{T}}}
\nc{\CU}{{\mathcal{U}}}
\nc{\CV}{{\mathcal{V}}}
\nc{\CW}{{\mathcal{W}}}
\nc{\CZ}{{\mathcal{Z}}}

\nc{\cM}{{\check{\mathcal M}}{}}
\nc{\csM}{{\check{\mathcal A}}{}}
\nc{\oM}{{\overset{\circ}{\mathcal M}}{}}
\nc{\obM}{{\overset{\circ}{\mathbf M}}{}}
\nc{\oCA}{{\overset{\circ}{\mathcal A}}{}}
\nc{\obA}{{\overset{\circ}{\mathbf A}}{}}
\nc{\ooM}{{\overset{\circ}{M}}{}}
\nc{\osM}{{\overset{\circ}{\mathsf M}}{}}
\nc{\vM}{{\overset{\bullet}{\mathcal M}}{}}
\nc{\nM}{{\underset{\bullet}{\mathcal M}}{}}
\nc{\oD}{{\overset{\circ}{\mathcal D}}{}}
\nc{\obD}{{\overset{\circ}{\mathbf D}}{}}
\nc{\oA}{{\overset{\circ}{\mathbb A}}{}}
\nc{\cp}{{\overset{\circ}{\mathbf p}}{}}
\nc{\oU}{{\overset{\bullet}{\mathcal U}}{}}
\nc{\oZ}{{\overset{\circ}{\mathcal Z}}{}}
\nc{\ofZ}{{\overset{\circ}{\mathfrak Z}}{}}
\nc{\as}{\varkappa}
\nc{\vkb}{{\boldsymbol{\varkappa}}}
\nc{\fa}{{\mathfrak{a}}}
\nc{\fb}{{\mathfrak{b}}}
\nc{\fg}{{\mathfrak{g}}}
\nc{\fgl}{{\mathfrak{gl}}}
\nc{\fh}{{\mathfrak{h}}}
\nc{\fri}{{\mathfrak{i}}}
\nc{\fj}{{\mathfrak{j}}}
\nc{\fm}{{\mathfrak{m}}}
\nc{\fn}{{\mathfrak{n}}}
\nc{\ft}{{\mathfrak{t}}}
\nc{\fu}{{\mathfrak{u}}}
\nc{\fz}{{\mathfrak{z}}}
\nc{\fp}{{\mathfrak{p}}}
\nc{\frr}{{\mathfrak{r}}}
\nc{\fs}{{\mathfrak{s}}}
\nc{\fsl}{{\mathfrak{sl}}}
\nc{\hsl}{{\widehat{\mathfrak{sl}}}}
\nc{\hgl}{{\widehat{\mathfrak{gl}}}}
\nc{\hg}{{\widehat{\mathfrak{g}}}}
\nc{\chg}{{\widehat{\mathfrak{g}}}{}^\vee}
\nc{\hn}{{\widehat{\mathfrak{n}}}}
\nc{\chn}{{\widehat{\mathfrak{n}}}{}^\vee}

\nc{\fA}{{\mathfrak{A}}}
\nc{\fB}{{\mathfrak{B}}}
\nc{\fD}{{\mathfrak{D}}}
\nc{\fE}{{\mathfrak{E}}}
\nc{\fF}{{\mathfrak{F}}}
\nc{\fG}{{\mathfrak{G}}}
\nc{\fK}{{\mathfrak{K}}}
\nc{\fL}{{\mathfrak{L}}}
\nc{\fM}{{\mathfrak{M}}}
\nc{\fN}{{\mathfrak{N}}}
\nc{\frP}{{\mathfrak{P}}}
\nc{\fU}{{\mathfrak{U}}}
\nc{\fZ}{{\mathfrak{Z}}}

\nc{\bc}{{\mathbf{c}}}
\nc{\be}{{\mathbf{e}}}
\nc{\bj}{{\mathbf{j}}}
\nc{\bn}{{\mathbf{n}}}
\nc{\bp}{{\mathbf{p}}}
\nc{\bq}{{\mathbf{q}}}
\nc{\bfu}{{\mathbf{u}}}

\nc{\bx}{{\mathbf{x}}}
\nc{\by}{{\mathbf{y}}}
\nc{\bw}{{\mathbf{w}}}
\nc{\bA}{{\mathbf{A}}}
\nc{\bB}{{\mathbf{B}}}
\nc{\bC}{{\mathbf{C}}}
\nc{\bD}{{\mathbf{D}}}
\nc{\bF}{{\mathbf{F}}}
\nc{\bH}{{\mathbf{H}}}
\nc{\bK}{{\mathbf{K}}}
\nc{\bM}{{\mathbf{M}}}
\nc{\bN}{{\mathbf{N}}}
\nc{\bO}{{\mathbf{O}}}
\nc{\bT}{{\mathbf{T}}}
\nc{\bV}{{\mathbf{V}}}
\nc{\bW}{{\mathbf{W}}}
\nc{\bX}{{\mathbf{X}}}
\nc{\bP}{{\mathbf{P}}}
\nc{\bZ}{{\mathbf{Z}}}

\nc{\sA}{{\mathsf{A}}}
\nc{\sB}{{\mathsf{B}}}
\nc{\sC}{{\mathsf{C}}}
\nc{\sD}{{\mathsf{D}}}
\nc{\sF}{{\mathsf{F}}}
\nc{\sK}{{\mathsf{K}}}
\nc{\sM}{{\mathsf{M}}}
\nc{\sO}{{\mathsf{O}}}
\nc{\sQ}{{\mathsf{Q}}}
\nc{\sP}{{\mathsf{P}}}
\nc{\sT}{{\mathsf{T}}}
\nc{\sZ}{{\mathsf{Z}}}
\nc{\sfp}{{\boldsymbol{\phi}}}
\nc{\sr}{{\mathsf{r}}}
\nc{\sfb}{{\mathsf{b}}}
\nc{\sfc}{{\mathsf{c}}}
\nc{\sd}{{\mathsf{d}}}

\nc{\BK}{{\bar{K}}}

\nc{\tA}{{\widetilde{\mathbf{A}}}}
\nc{\tB}{{\widetilde{\mathcal{B}}}}
\nc{\tg}{{\widetilde{\mathfrak{g}}}}
\nc{\tG}{{\widetilde{G}}}
\nc{\TM}{{\widetilde{\mathbb{M}}}{}}
\nc{\tO}{{\widetilde{\mathsf{O}}}{}}
\nc{\tU}{{\widetilde{\mathfrak{U}}}{}}
\nc{\TZ}{{\tilde{Z}}}
\nc{\tx}{{\tilde{x}}}
\nc{\tbv}{{\tilde{\bv}}}
\nc{\tfP}{{\widetilde{\mathfrak{P}}}{}}
\nc{\tz}{{\tilde{\zeta}}}
\nc{\tmu}{{\tilde{\mu}}}

\nc{\urho}{\underline{\rho}}
\nc{\uB}{\underline{B}}
\nc{\uC}{{\underline{\mathbb{C}}}}
\nc{\ui}{\underline{i}}
\nc{\uj}{\underline{j}}
\nc{\ofP}{{\overline{\mathfrak{P}}}}
\nc{\oB}{{\overline{\mathcal{B}}}}
\nc{\og}{{\overline{\mathfrak{g}}}}
\nc{\oI}{{\overline{I}}}

\nc{\eps}{\varepsilon}
\nc{\hrho}{{\hat{\rho}}}

\nc{\one}{{\mathbf{1}}}
\nc{\two}{{\mathbf{t}}}

\nc{\Rep}{{\mathop{\operatorname{\rm Rep}}}}

\nc{\Tot}{{\mathop{\operatorname{\rm Tot}}}}

\nc{\Ext}{{\mathop{\operatorname{\rm Ext}}}}

\nc{\CHom}{{\mathop{\operatorname{{\mathcal{H}}\it om}}}}

\nc{\Id}{{\mathop{\operatorname{\rm Id}}}}

\nc{\defi}{{\mathop{\operatorname{\rm def}}}}
\nc{\length}{{\mathop{\operatorname{\rm length}}}}

\nc{\Cliff}{{\mathsf{Cliff}}}
\nc{\Fl}{{\mathsf{Fl}}}
\nc{\Fr}{{\mathsf{Fr}}}
\nc{\fr}{{\text{\small\textsf{F}}}}
\nc{\Fib}{{\mathsf{Fib}}}
\nc{\Coh}{{\mathsf{Coh}}}
\nc{\FCoh}{{\mathsf{FCoh}}}

\nc{\cplus}{{\mathbf{C}_+}}
\nc{\cminus}{{\mathbf{C}_-}}
\nc{\cthree}{{\mathbf{C}_*}}
\nc{\Qbar}{{\bar{Q}}}

\nc{\bh}{{\bar{h}}}
\nc{\bOmega}{{\overline{\Omega}}}

\nc{\seq}[1]{\stackrel{#1}{\sim}}
\nc{\ul}{\underline}
\nc{\ol}{\overline}

\def\hilbb{{\operatorname{Hilb}}^{(1)}}
\begin{document}

\centerline{\Large{\bf Cherednik algebras  and Hilbert schemes}}
\vskip 4pt

\centerline{\Large{\bf in characteristic $p$}}
\vskip 10pt

\centerline{\sc Roman Bezrukavnikov, Michael Finkelberg and Victor Ginzburg}
\vskip 4pt

\centerline{(with Appendix by {\sc {Pavel Etingof}})}
\vskip 10mm

\hskip 80mm {\em To David Kazhdan with admiration} 

\begin{abstract}
We prove a {\em localization} theorem for the type $\mathbf{A}_{n-1}$
rational Cherednik 
algebra $\hh_c=\hh_{1,c}(\mathbf{A}_{n-1})$
over $\overline{\mathbb{F}}_p$, an algebraic closure of
the finite field.
In the most interesting special case where $c\in \mathbb{F}_p$,
we construct an Azumaya algebra $\ah_c$ on 
$\Hilb$, the Hilbert scheme of 
$n$ points in the plane, such that $\Gamma(\Hilb,
\,\ah_c)=\hh_c$.
Our localization theorem provides an equivalence between the
bounded derived categories of $\hh_c$-modules
and   sheaves of coherent   $\ah_c$-modules on $\Hilb$,
respectively. Furthermore, we show that the Azumaya algebra
splits on the formal neighborhood of each fiber of the Hilbert-Chow morphism.
This provides a link between
our results and those of Bridgeland-King-Reid and Haiman.
\end{abstract}

{\centerline{\bf Table of Contents}
\vskip -5mm
$\hspace{20mm}$ {\footnotesize \parbox[t]{115mm}{\,

\hp${}_{}$\hp1.{ $\;\,\,$} {\tt Introduction} \newline
\hp2.{ $\;\,\,$} {\tt Crystalline differential operators}\newline
\hp3.{ $\;\,\,$} {\tt Hamiltonian reduction in characteristic $p$}\newline
\hp4.{ $\;\,\,$} {\tt Azumaya algebras via Hamiltonian reduction}\newline
\hp5.{ $\;\,\,$} {\tt The rational Cherednik algebra of type $\mathbf{A}_{n-1}$}\newline
\hp6.{ $\;\,\,$} {\tt An Azumaya algebra on the Hilbert scheme}\newline
\hp7.{ $\;\,\,$} {\tt Localization functor for Cherednik algebras}\newline
\hp8.{ $\;\,\,$}  {\tt Induction functor and comparison with \cite{EG}}\newline
\hp9.{ $\;\,\,$} {\tt Appendix by Pavel Etingof:}\newline
\hp{ $\;\,\enspace\enspace$} {\tt The $p$-center of Symplectic
 reflection algebras}\newline
\newline
}}

\section{Introduction}
\subsection{} Let 
$c\in\BQ$ be 
a  rational number,
and $\hh_{1,c}(\mathbf{A}_{n-1})$
the rational Cherednik 
algebra of  type $\mathbf{A}_{n-1}$ with parameters $t=1$ and $c$
 that has been considered in
\cite{EG}
(over the ground field of complex numbers).

For all primes $p \gg n$,
we can reduce $c$ modulo $p$.
Thus, $c$ becomes an element of the finite
field $\BF_p$. We let $\k=\k_p$ be an
algebraic closure of $\BF_p$,
and let $\hh_c:=\hh_{1,c}(\mathbf{A}_{n-1},\k_p)$ 
be the Cherednik algebra,
viewed as an algebra over $\k_p$.
Unlike the case of characteristic zero,
the algebra $\hh_c$ has a large center,
called the $p$-center.
The spectrum of the $p$-center
is isomorphic to 
$[(\BA^2)^n/S_n]\tw$,
the Frobenius twist of
the $n$-th symmetric power of the plane $\BA^2$.

\subsection{} We consider $\Hilb$, the Hilbert scheme (over $\k_p$)
of $n$ points in the plane, see e.g. \cite{Na1}.
There is a canonical
 Hilbert-Chow map $\Ups: \Hilb\to (\BA^2)^n/S_n$
that induces an algebra isomorphism
\beq{hchow}
\G(\Hilb,\,\oo_{_{\Hilb}})\cong\k\big[(\BA^2)^n/S_n\big].
\eeq

Let $\hilbb$ denote the  Frobenius twist of $\Hilb,$
a scheme isomorphic to $\Hilb$ and equipped with a
canonical Frobenius morphism $\Fr:\Hilb\to\hilbb.$
We introduce an  Azumaya algebra $\ah_c$ on $\hilbb$
of degree $n!\cdot p^{n}$
(recall that  an Azumaya algebra has degree $r$ if each of its geometric
fibers is isomorphic
to the algebra of $r\times r$-matrices).
For all
sufficiently large primes $p$,
we construct a natural algebra isomorphism
(a version of the Harish-Chandra isomorphism from~\cite{EG})
\beq{intro1}
\G(\hilbb,\,\ah_c)\iso\hh_c.
\eeq
The restriction of
this isomorphism to the subalgebra $\G(\hilbb,\oo_{_{\hilbb}})$
yields, via \eqref{hchow}, the above mentioned   isomorphism between 
the algebra $\k\big[\bigl((\BA^2)^n/S_n\bigr)\tw\big]$
and  the $p$-center.
\begin{rem} More generally, for any $c\in \k$,
not necessarily an element of $\BF_p$,
there is an Azumaya algebra on the Calogero-Moser space
with parameter $c^p-c$ such that an analogue of isomorphism
\eqref{intro1} holds for the  Calogero-Moser space
instead of the Hilbert scheme.
This case is somewhat less interesting since
the Calogero-Moser space is affine while
 the Hilbert scheme is not.
\erem

The main idea used in the construction of isomorphism \eqref{intro1}
is to compare Nakajima's  description of $\Hilb$ by means of {\em Hamiltonian
reduction}, see \cite{Na1}, with (a refined version,
see \S\ref{hh}) of the construction introduced in \cite{EG} 
 describing  the
{\em spherical subalgebra} of $\hh_c$ as a
{\em quantum} Hamiltonian reduction of an algebra
of differential operators.

\subsection{}
We introduce the following set of rational numbers
\beq{qgood}
\QG=\{c\in \BQ\mid c\geq 0\enspace\&\enspace
c\not\in \frac{1}{2}+\Z\}.
\eeq

One of our main results (Theorem \ref{ah}) reads

\begin{theorem}\label{equiv_intro} Fix  $c\in\QG.$
 Then, there exists a constant
$d=d(c)$
such that for all primes $p> d(c)$,
the functor $\RG:\ D^b(\ah_c\dash\on{Mod})\to D^b(\hh_c\dash\on{Mod})$
is a triangulated equivalence between the bounded derived categories
of sheaves of coherent $\ah_c$-modules and finitely
generated $\hh_c$-modules, respectively, whose
inverse is the localisation functor
$M\mapsto{\ah_c\stackrel{_L}{\otimes}_{\hh_c}M}.$

Moreover, we
have
$H^i(\hilbb,\ah_c)=0,\,\forall i>0$.
\end{theorem}

\subsection{} Now, fix  $\xi\in [(\BA^2)^n/S_n]\tw$,
a point in the Frobenius twist of $(\BA^2)^n/S_n$.
We
write  $\op{Hilb}_\xi\tw=\Ups\inv(\xi)$ for the corresponding
 fiber  of the Frobenius twist of the 
Hilbert-Chow map, and let 
$\widehat{\op{Hilb}}_\xi\tw=\widehat{\Ups\inv(\xi)}$
denote its formal neighborhood, the completion of
$\hilbb$ along the subscheme $\op{Hilb}_\xi\tw$.

The theorem below, based on a similar result in \cite{bk},
says that the Azumaya algebra
$\ah_c$ {\em splits} on the formal  neighborhood of each fiber of the
Hilbert-Chow map, that is, we have
the following result (see Theorem \ref{splitting})

\begin{theorem}\label{spl_intro} For each $\xi\in [(\BA^2)^n/S_n]\tw$,
 there exists a vector bundle
$\CV_{c,\xi}$ on $\widehat{\op{Hilb}}_\xi\tw$ such that one has
$$\ah_c\big|_{\widehat{\op{Hilb}}_\xi\tw}\,\cong\,
(\eend\,\CV_{c,\xi})^{\op{opp}}.
$$
\end{theorem}
\begin{rem}\label{choice}
 The   splitting bundle is not unique;
it is only determined up to
twisting by an invertible sheaf.
\erem

Given $\xi$ as above, let
$\fm_\xi$ be the corresponding 
maximal ideal  in the $p$-center
of $\hh_c$.   Let $\wh{\hh}_{c,\xi}$,
resp. $\wh{\ah}_{c,\xi}=\ah_c\big|_{\widehat{\op{Hilb}}_\xi\tw}$, 
 be the $\fm_\xi$-adic completion
of $\hh_c$,
resp. of $\ah_c$. We write
$D^b(\wh{\hh}_{c,\xi}\dash\on{Mod})$,
resp. $D^b(\wh{\ah}_{c,\xi}\dash\on{Mod})$, for the bounded derived category
of finitely-generated complete topological
 $\wh{\hh}_{c,\xi}$-modules,
resp. $\wh{\ah}_{c,\xi}$-modules.
On the other hand, let $D^b(\Coh(\wh{\op{Hilb}}_\xi\tw))$
be the bounded derived category of coherent sheaves on 
the formal scheme $\widehat{\op{Hilb}}_\xi\tw$.

Fix $c\in\QG$. Then, for all primes $p>d(c)$,
Theorems \ref{equiv_intro}-\ref{spl_intro} imply the following

\begin{corollary}\label{intro2}
 The category $D^b(\wh{\hh}_{c,\xi}\dash\on{Mod})$ 
 is equivalent to $D^b(\Coh(\wh{\op{Hilb}}\tw_\xi))$.
\end{corollary}

Now let $\xi=0$ be the zero point in $[(\BA^2)^n/S_n]\tw$.
The fiber $\hilbb_0$ is isomorphic to the
(Frobenius twist of the) punctual Hilbert scheme.
 We will show that the algebra $\wh{\hh}_{c,0}$
contains a canonical dense
$\Z^2$-graded subalgebra $\hh^\circ=\oplus_{k,l\in\Z}\,\hh^{k,l}$,
see \S\ref{bigrading_sec}. 
The category of $\Z^2$-graded $\hh^\circ$-modules
may be thought of as a `{\em mixed version}'
of the category ${\wh{\hh}_{c,0}}\dash\on{Mod}$,
cf. \cite[Definition 4.3.1]{BGS}.

Recall
that the algebra  $\hh_c$ contains  a canonical 
 ${\mathfrak{s}\mathfrak{l}}_2$-triple, see  \cite{BEG}.
Let  ${\mathsf{h}}\in \hh_c$ denote the semisimple
element of that triple. 
We expect that the above mentioned $\Z^2$-grading 
has the property that,
for any $u\in \hh^{k,l}$, we have ${\mathsf{h}}\cdot u-
u\cdot {\mathsf{h}}= (k-l)\cdot u$.

\subsection{}\label{intro_conjectures}
Let $\k[S_n]$ denote the group
algebra of the Symmetric group on $n$ letters.
Write $\irr$ for the set of isomorphism classes
of simple $\k[S_n]$-modules. This set is labelled
by partitions of $n$, since by our assumptions $\charp\k>n$. 
In particular, we have the trivial 1-dimensional
representation $\triv$, and the sign representation $\sign$.

Let $\overline{\hh}:=\k[[\BA^{2n}]]\#S_n$ be the cross-product
of $S_n$ with  $\k[[\BA^{2n}]]$, the  algebra
of formal power series
in $2n$ variables
acted on by $S_n$ in a natural way.
We consider $D^b(\overline{\hh}\dash\on{Mod})$,
the bounded derived category of (finitely-generated)
complete topological
$\overline{\hh}$-modules.

Given a simple $\k[S_n]$-module $\tau$, write
$\tau_{_{\overline{\hh}}}$
for the corresponding  $\overline{\hh}$-module
obtained by pullback via the natural projection
$\overline{\hh}=\k[[\BA^{2n}]]\#S_n
\to \k[S_n],\, f\rtimes w\mapsto f(0)\cdot w$.
Similarly,
let
$L_\tau$ denote the corresponding simple  highest weight $\hh_c$-module,
the unique simple quotient of the standard
 $\hh_c$-module associated with $\tau$,
see \cite{DO},\cite{BEG}.

The results of Bridgeland-King-Reid
\cite{BKR} and Haiman \cite{H}, see also \cite{bk}, provide
an equivalence of categories
$$\op{BKR}:\ D^b(\Coh(\Hilb))\iso D^b(\k[\BA^{2n}]\#S_n\dash\on{Mod}),\quad
\scr F\mapsto \RG(\Hilb,\,\CP\stackrel{_L}{\otimes} \scr F),
$$
where $\CP$ denotes  
 the {\em Procesi 
bundle}, the `unusual' tautological rank $n!$  vector
bundle
on $\Hilb$ considered in \cite{H}.
Restricting this equivalence to the completion of the
zero fiber of the Hilbert-Chow
map, and using Corollary \ref{intro2},
one obtains the following composite equivalence
\beq{intro3}
\xymatrix{
{D^b({\wh{\hh}_{c,0}}\dash\on{Mod})}
\ar[rrr]_<>(0.5){\sim}^<>(0.5){\text{Corollary \ref{intro2}}}&&&
D^b(\Coh(\wh{\op{Hilb}}\tw_0))\ar[rr]_<>(0.5){\sim}^<>(0.5){\op{BKR}}&&
D^b({\overline{\hh}\dash\on{Mod}}).
}
\eeq
We recall that the equivalence of Corollary \ref{intro2}
involves a choice of splitting bundle $\CV_{c,0}$, cf. Remark 
\ref{choice}. 
This choice may be specified by the following

\begin{conj}\label{triv} Fix $c\in \QG$. Then, for all $p\gg 0$, we have

\vi One can choose the splitting bundle
$\CV_{c,0}$ in such a way that $\Gamma(\op{Hilb}\tw_0, \CV_{c,0}(-1))=L_\sign$.

\vii With this choice of $\CV_{c,0}$, the composite equivalence in \eqref{intro3}
preserves the natural $t$-structures,
in particular, induces an equivalence
${\wh{\hh}_{c,0}}\dash\on{Mod}\iso{\overline{\hh}\dash\on{Mod}},$
 of {\em abelian} categories, such that $L_\tau$ goes to
$\tau_{_{\overline{\hh}}}$, for any simple $S_n$-module $\tau$.
\end{conj}

For $c=\frac{1}{n}$, we expect that
$\Gamma(\op{Hilb}\tw_0, \CV_{c,0}\otimes
\op{BKR}\inv(\triv_{_{\overline{\hh}}}))$
is a  1-dimensional vector space that supports the
trivial representation of the group $S_n\sset\hh_c$.

\subsection{} In the special case of rank one,
i.e., for $n=2$, a complete classification and explicit
construction of 
simple $\hh_c$-modules (for $\charp\k>0$) has been
 obtained by Latour~\cite{La}.

 It seems certain that our results in characteristic
$p$ have their characteristic zero counterparts
for the double-affine Hecke algebra of type $\mathbf{A}_{n-1}$,
specialized 
at a root of unity, cf. \cite{Ch}; in that case 
one has to replace $\Hilb$ by $\Hilgm$, cf. \cite{Ob}.

Also, it is likely that the results of the present paper
can be generalized to the case of 
{\em symplectic reflection algebras}
associated with wreath products $\mathbf{\G}_n=S_n \ltimes\G^n,$
where $\G$ is a finite subgroup in $SL_2(\k)$, see \cite[\S11]{EG}.
More generally, we are going to study Azumaya algebras
arising via quantum Hamiltonian reduction from
the general Nakajima quiver varieties,
cf. \cite{Na2}, (of which
wreath-products are special cases). Our technique is ideally
suited for such a generalization, that
has been, in effect, suggested earlier
by  Nakajima and the first author. 

In another direction,
the general results of
\S4 below apply  verbatim to quantizations of Slodowy
slices considered in \cite{Pr}, see also \cite{GG}.

We are going to explore
these topics elsewhere.

\subsection*{Acknowledgments.}{\small
The authors are  grateful to I. Gordon for informing one of us
about his unpublished results included in~\cite{g}.
Also, we would like to thank A. Premet for pointing out
several inaccuracies involving restricted Lie
algebras that have occurred in the original version
of the paper.
The research of the first and third author
was partially supported by the NSF. The research
 of the second author
 was conducted for the Clay Mathematical
Institute and partially supported by the CRDF award RM1-2545-MO-03 and
RM1-2694,
and the ANR program "GIMP", contract number ANR-05-BLAN-0029-01.
He is grateful to V.~Vologodsky and V.~Dotsenko for patient
explanations,
and to the University of Chicago and Northwestern
University for hospitality and support.}

\section{Crystalline differential operators.}\label{crys}
\subsection{}\label{setup}
Unless specified otherwise,  we will be working
over the ground field $\k$, an algebraically closed field 
of characteristic $\charp\k=p>0$. We write
$\Fr: \k\to\k,\,k\mapsto k^p$
for the Frobenius automorphism.
Given a  $\k$-vector space $E$, it is convenient to introduce $E\tw$,
 a  vector space with the same underlying 
additive group
as $E$,  but with a `twisted'
$\k$-linear structure given by
$k\ccirc e:= \Fr\inv(k)\cdot e,\,\forall k\in \k,e\in E.$
Note that if $A$ is a $\k$-algebra, then
 the map
$A\to A, a\mapsto a^p$  is
an additive but {\em  not $\k$-linear}
map, that becomes {\em  $\k$-linear}
if considered as a map $A\tw\to A$.

Given an additive map $f: E\to F$ between 
two $\k$-vector spaces, we  say that $f$ is

\pb{$p$-{\em linear}, if $f(k\cdot e)= k^p\cdot f(e)$
for any $k\in k,e\in E$, i.e., if the
corresponding map $E\tw\to F$ is $\k$-linear;}

\pb{\gre if both vector spaces are equipped with $\Z$-gradings
$E=\bigoplus E(i),\,F=\bigoplus F(i)$,
and we have $f(E(i))\sset F(p\cdot i)$
for all $i$. }

\subsection{}
Let $X$ be a smooth algebraic veriety over $\k$
with structure sheaf  $\oo_X$. Write
$\k[X]$ for the corresponding algebra
 of global sections.
We let $X\tw$ denote an algebraic variety with the same
structure sheaf as $X$ but with the `twisted'
$\k$-linear structure.
Thus, $\oo_{X\tw}:=(\oo_X)\tw$,
and there is a canonical morphism  $\Fr: X \to X\tw$  called {\em Frobenius
morphism}, such that the map
$f\mapsto f^p$ on regular functions  becomes identified with
the natural sheaf imbedding $\Fr\hdot\oo_{X\tw} \into \oo_X$.

We write $\TT_X$ for the tangent
sheaf on $X$, and let $T^*X$ denote the total space of the cotangent bundle.
There is a canonical
isomorphism $T^*[X\tw]\cong [T^*X]\tw$,
and we will use the notation $T^*X\tw$ for these two isomorphic
varieties, and $\pi: T^*X\tw\to X\tw$ for the natural
projection.
The space $T^*X$, resp. $T^*X\tw$,  has a canonical symplectic structure,
which makes  $\k[T^*X]$
a Poisson algebra.

 Let $\D_X$ denote the
sheaf of
  crystalline
differential operators on $X$, that is, a sheaf of
algebras generated by $\oo_X$ and $\TT_X$. Let
 $\D(X):=\G(X,\D_X)$ denote the corresponding algebra of global sections.
More generally, given a locally-free coherent sheaf (= vector bundle)
$\CL$ on $X$, let $\D_X(\CL):=\CL\bigotimes_{\oo_X}
\D_X\bigotimes_{\oo_X}\CL^\vee$ be the
sheaf of
differential operators on $\CL$,
and $\D(X,\CL):=\Gamma(X,\D_X(\CL))$ the 
 algebra of its global
sections.

\subsection{}
The sheaf $\D_X$ is known to have a large center.
Specifically, for any vector field
$\xi\in \TT_X$, the $p$-th power of $\xi$
acts as a derivation, hence, gives rise to another
vector filed, $\xi\tp\in\TT_X.$
The assignment $\xi\mapsto \xi^p-\xi\tp$
extends to a canonical algebra 
imbedding
\beq{bzd}
\bzd:\
\xymatrix{
\sym \TT_{X\tw}\;\ar@{^{(}->}[rr]&&\Fr\idot\D_X
}, \quad\xi\mto \xi^p-\xi\tp,
\eeq 
(of sheaves
on $X\tw$) whose image, to be denoted $\z_{X\tw}\sset \Fr\idot\D_X$,
equals 
the center of $\Fr\idot\D_X$.
Therefore, the isomorphism $\pi\idot\oo_{T^*X\tw}\simeq \sym \TT_{X\tw}$
makes  $\Fr\idot\D_X$ a sheaf of $\pi\idot\oo_{T^*X\tw}$-algebras.
This way, 
we may (and will) view  $\Fr\idot\D_X$ as a coherent sheaf
on $T^*X\tw$, to be denoted $\D\tw$. 

 The sheaf $\D_X$ comes equipped with a standard increasing
filtration $\D^{\leq k}_X,$
$k=0,1,\ldots,$
by the order of differential operator. For the associated graded
sheaf, one has a graded algebra isomorphism $\gr\D_X\cong
\sym \TT_X=
\pi\idot\oo_{T^*X}$.
The filtration on $\D_X$ induces a filtration 
on  $\Fr\idot\D_X$ and also the filtration
 $\z^{\leq i}_{X\tw}:=\Fr\idot\D^{\leq i}_X\cap
\z_{X\tw}$ on  the central subalgebra $\z_{X\tw}$.
Observe that the latter algebra already has a natural grading
obtained, via the isomorphism $\z_{X\tw}=\sym \TT_{X\tw}$, from the standard
grading on the Symmetric algebra.
With this grading, one has a \gre algebra  isomorphism
$\sym\TT_{X\tw}\iso\gr\z_{X\tw},$ cf. Sect \ref{setup},
i.e. we have $\sym^i\TT_{X\tw}\iso
\gr^{p\cdot i}\z_{X\tw}
,\,\forall i$.

Now, 
view $\Fr\idot\D_X$ as a $\z_{X\tw}$-algebra. Thus,
 $\gr(\Fr\idot\D_X)$ becomes a $\gr\z_{X\tw}$-algebra
that may be viewed, by the isomorphism $\Spec(\gr\z_{X\tw})=T^*X\tw$,
as a $\Gm$-equivariant coherent sheaf on $T^*X\tw$.
On the other hand, consider the Frobenius morphism
$\Fr^{T^*X}: T^*X \to [T^*X]\tw$ and 
view $\Fr^{T^*X}\idot\oo_{T^*X}$ as a $\Gm$-equivariant coherent sheaf
of algebras on $T^*X\tw$, a $\Gm$-variety.
With this understood,
there is a natural $\Gm$-equivariant algebra  isomorphism
\beq{gr}
\gr\bigl(\Fr\idot\D_X\bigr)\simeq \pi\idot\Fr^{T^*X}\idot\oo_{T^*X}.
\eeq

\subsection{The Rees algebra.}\label{rees_sec} Let $\DD$ be an associative
algebra, and write  $\DD[t]:=\k[t]\otimes\DD,$
where $t$ is an indeterminate. We put a grading on
$\DD[t]$ by assigning $\DD$ grade degree zero,
and setting $\deg t=1$.
Recall that, given an increasing filtration
$0=\DD_{-1}\sset\DD_0\sset \DD_1\sset\ldots,$ on $\DD$, such that
$\DD_i\cdot\DD_j\sset\DD_{i+j}$
and $\bigcup_{i\geq 0}\,\DD_i=\DD$,
one defines the  {\em Rees algebra} of $\DD$ 
as the following graded subalgebra:
$$\rees\DD:=\sum\nolimits_{i\geq 0}\,t^i\cdot \DD_i
\;\sset\;
\DD[t].
$$

There are standard isomorphisms
\beq{rees_iso}
(\rees\DD)|_{\{0\}}\cong \gr\DD,\quad\text{and}\quad
(\rees\DD)_{(t)}\cong \k[t,t\inv]\otimes\DD,
\eeq
where, for any $\k[t]$-algebra $R$ 
we let $R_{(t)}:=\k[t,t\inv]\bigotimes_{\k[t]}R$ denote the localization of $R$,
and for any $s\in \k$
we use the notation $R|_{\{s\}}:=R/(t-s)R$.

Conversely, given a flat $\Z_{\geq 0}$-graded $\k[t]$-algebra
$\rr=\bigoplus_{i\geq 0}\,\rr(i)$, set $\DD:=\rr|_{\{1\}}=\rr/(t-1)\rr.$ This is
a $\k$-algebra equipped with a canonical
increasing filtration 
$\DD_i,\,i=0,1,\ldots,$ and with a canonical
graded algebra isomorphism $\gr\DD\cong \rr|_{\{0\}}=\rr/t\rr$.

The filtration on $\DD$ is defined in the following way.
Put $\DD[t,t\inv]:=\k[t,t\inv]\bigotimes_\k \DD$,
and view it as a $\Z$-graded  $\k[t,t\inv]$-algebra.
Further, inverting $t$, we get from $\rr$  a $\Z$-graded  $\k[t,t\inv]$-algebra
$\rr_{(t)}$
that contains $\rr$ as a $\k[t]$-subalgebra.
The definition of $\DD$ provides an
 isomorphism $\phi:$ ${\rr_{(t)}/(t-1)\rr_{(t)}}\iso \DD$ that
admits a unique lift to
the following graded $\k[t,t\inv]$-algebra isomorphism 
$${\phi}_{(t)}:
\rr_{(t)}\iso\DD[t,t\inv]=\k[t,t\inv]\otimes\frac{\rr}{(t-1)\rr},\quad
\rr(i)\ni u\mapsto t^i\otimes \bigl(u\,\text{mod}(t-1)\rr\bigr).
$$

The above mentioned increasing filtration on $\DD$ is defined by
\beq{filt}
\DD_i:= \DD\cap {\phi}_{(t)}(t^{-i}\cd\rr),\quad
i=0,1,\ldots\;.
\eeq

Assume next  that
$\DD=\bigoplus_{i\ge 0}\,\DD(i)$ is  a {\em graded}
algebra,
and view it as a filtered algebra with filtration
being induced by the grading,
that is, defined by $\DD_i:=\bigoplus_{j\leq i}\,\DD(j)$.
Then, we have
$$\rees\DD =\sum\nolimits_{i\geq 0}\,t^i\cdot \bigl(\oplus_{j\leq i}\,
\DD(j)\bigr)= \bigl(\sum\nolimits_{i\geq 0}\, t^i\cdot\DD(i)\bigr)[t].$$
We see that, for a graded algebra $\DD$,
one has the following  graded
algebra isomorphism
\beq{rees_gr}
\k[t]\otimes\DD\iso \rees\DD
,\quad
\k[t]\bigotimes\bigl(\oplus_{_i}\DD(i)\bigr)
\ni\; f\otimes (\sum_i\, u_i)\;\mto\;
f\cdot \sum_i\, t^i\cd u _i.
\eeq

\subsection{The sheaf $\rr\D\tw$.}\label{rrDtw}
We apply the Rees algebra construction
 to  $\Fr\idot\D_X,$ viewed as a sheaf of filtered algebras.
Thus, we get a sheaf $\rees\Fr\idot\D_X$
of graded $\oo_{X\tw}$-algebras.
In $\rees\Fr\idot\D_X$,
we also have a central
subalgebra $\rees\z_{X\tw}\sset\rees\Fr\idot\D_X.$

By \eqref{rees_gr}, the canonical grading on $\z_{X\tw}$
provides a
\gre algebra isomorphism $\rees\z_{X\tw}
\simeq \k[t]\otimes \z_{X\tw}.$ 
Thus, we obtain the following canonical algebra maps
\beq{central_rees}
\xymatrix{
\sym \TT_{X\tw}\ar[r]^<>(0.5){\eqref{bzd}}_<>(0.5){\sim}&
\z_{X\tw}
\ar[rr]^<>(0.5){z\mapsto 1\otimes z}&&
\k[t]\otimes \z_{X\tw}\ar@{=}[r]^<>(0.5){\eqref{rees_gr}}&\rees\z_{X\tw}
\ar@{^{(}->}[r]&
\rees\Fr\idot\D_X.
}
\eeq 
The composite map in \eqref{central_rees} is a \gre map to be denoted
$\bzdr$. This map  specializes
at $t=0$ to the map $\xi\mapsto \xi^p$,
and  at $t=1$ to   \eqref{bzd}.

Further, the algebra isomorphism  $\rees\z_{X\tw}
\simeq \k[t]\otimes \z_{X\tw}$ yields 
a direct product decomposition
\beq{prod_spec}
\Spec\bigl(\rees\z_{X\tw}\bigr)\;\cong\;
\BA^1\times \Spec\z_{X\tw}=
\BA^1\times T^*X\tw.
\eeq
We will often identify $\rees\Fr\idot\D_X$, a
graded $\rees\z_{X\tw}$-algebra,
  with a  coherent sheaf
of algebras on
$\Spec(\rees\z_{X\tw})$,
that is, on $\BA^1\times T^*X\tw$.
The resulting sheaf on  $\BA^1\times T^*X\tw$,
to be denoted
 $\rr\D\tw$, is easily seen to be {\em flat} with respect to
 the first factor $\BA^1$.
Moreover, \eqref{rees_iso} yields, in view of
\eqref{gr}, the  following
isomorphisms
of sheaves of algebras on $T^*X\tw$, resp.,  on $(\BA^1\sminus\{0\})\times T^*X\tw$:
\begin{align}\label{rees_D}
(\rr\D\tw)\big|_{\{0\}\times
T^*X\tw}\,\cong\,\Fr^{T^*X}\idot\oo_{T^*X},\quad
(\rr\D\tw)\big|_{(\BA^1\sminus\{0\})\times T^*X\tw}\,\cong\,\op{pr}_2^*(\D\tw),
\end{align}
where $\pr_2: (\BA^1\sminus\{0\})\times T^*X\tw\to T^*X\tw$ is the second projection.

Further, the grading on the Rees algebra $\rees\z_{X\tw}$
makes $\Spec(\rees\z_{X\tw})$
a $\Gm$-variety.
It follows from formula \eqref{rees_gr} that
the natural $\Gm$-action on $\Spec\bigl(\rees\z_{X\tw}\bigr)$
corresponds, via \eqref{prod_spec},
to the $\Gm$-{\em diagonal} action on 
$\BA^1\times T^*X\tw$.
The sheaf $\rr\D\tw$ on $\BA^1\times T^*X\tw$
comes equipped with a canonical
 $\Gm$-equivariant structure.

\section{Hamiltonian reduction in characteristic $p$}
\subsection{Lie algebras in characteristic $p$.}\label{not_sect}
Let $A$ be a connected linear algebraic group
  over $\k$.
Write
 $A^{(1)}$ for the {\em Frobenius twist} of $A$, cf. \S\ref{crys},
an algebraic group isomorphic to $A$ and equipped with
an algebraic group morphism $\Fr: A\to A\tw$, called the
Frobenius morphism. The kernel of this morphism is an infinitesimal
group scheme $A_1\sset A$, called {\em  Frobenius kernel}.
By definition, one has an  exact sequence:
$$
1\too A_1 \too A \stackrel{\Fr}{\too}  A^{(1)}\too 1.
$$

The Lie algebra  $\aa:=\Lie A$ 
may be viewed as the vector space of left invariant vector 
fields on $A$. This  vector space comes equipped
with a natural structure
of $p$-Lie algebra, i.e., we have a
$p$-{\em power}  map
$\aa\to\aa,\,x\mapsto x^{[p]}$, see \cite{J} or  \cite{Ja}.

Let $\sym\aa$, resp. $\Ua$, be the symmetric, resp. 
 enveloping, algebra of $\aa$. The group $A$ acts
on $\sym\aa$ and $\Ua$ by algebra automorphisms via the adjoint action.

The standard increasing filtration
$\U\idot\aa$ on the enveloping algebra
gives rise to a graded algebra $\rees\Ua=\sum_{i\geq 0}\,t^i\cdot\U_i\aa$.
Jacobson's argument \cite[ch. V,\S7, (60)-(64)]{J} shows that
the  following assignment
\beq{pcenter}
\bzur:\  
\xymatrix{
\sym\aa\tw\;\ar@{^{(}->}[r]&
\rees\Ua,
}\quad\aa\tw\ni x\mto x^p-t^{p-1}\cd x^{[p]}
\eeq
gives a well-defined  $\Ad A$-equivariant  injective
 \gre algebra homomorphism, cf. also \cite{PS}.
The  image of the map $\bzur$ is  an $\Ad A$-stable
subalgebra
contained in the center of $\rees\Ua$.

Specialization of  the map  $\bzur$ at $t=0$ 
reduces to
the \gre algebra map
$\sym\aa\tw\into\sym\aa,\,x\mapsto x^p.$
On the other hand, specializing the map $\bzur$
at $t=1$, one obtains an algebra imbedding
$\bzu:=\bzur|_{t=1}\,: \sym\aa\tw\into \Ua$. The 
image $\z(\aa):=\bzu(\sym\aa\tw)$ of this
imbedding is a central subalgebra in
$\Ua$ generated by the elements
$\{x^p-x\tp\}_{x\in\aa},$ usually referred to as the
$p$-{\em center} of $\Ua$, cf. e.g. \cite{Ja}.
Thus, the map in \eqref{pcenter}
may be identified with the composite of the
following chain of 
 algebra homomorphisms, completely
analogous to those in \eqref{central_rees}:
$$
\xymatrix{
\sym\aa\tw\;\ar[r]_<>(0.5){\sim}^<>(0.5){\bzu}&
\z(\aa)
\ar[rr]^<>(0.5){z\mapsto 1\otimes z}&&
\k[t]\otimes \z(\aa)\ar@{=}[r]^<>(0.5){\eqref{rees_gr}}&\rees\z(\aa)\;
\ar@{^{(}->}[r]&
\;\rees\Ua
}.
$$

The adjoint action on $\z(\aa)$ of the Frobenius kernel
$A_1\sset A$ is trivial, hence, the $A$-action on $\z(\aa)$
factors through
$A\tw$.

\subsection{The Artin-Schreier map.}\label{twists}
Let $\aa^*$
denote the $\k$-linear dual of $\aa$,
and write $\atw:=(\aa\tw)^*$
for the  $\k$-linear dual of $\aa\tw$.
 For any linear function
$\lambda\in \aa^*$, the assignment
$x\mapsto \la(x)^p$ gives a $p$-linear map $\aa\to\k$,
that is, 
 $\k$-linear function on
$\aa\tw$, to be denoted $\la\tw.$
This way, one obtains a   $p$-linear map
$\aa^*\to (\aa\tw)^*,\,\la\mapsto\la\tw.$
The latter map gives a
canonical $\k$-vector space isomorphism
$(\aa^*)\tw\iso (\aa\tw)^*,\,\la\mapsto\la\tw.$

Let $\BX(\aa)\sset\aa^*$ denote the subspace
of fixed points of the coadjoint action of $A$ on $\aa^*$.
Such a fixed point may be viewed as an $\Ad A$-invariant
Lie algebra homomorphism $\aa\to\k$ (note that a
Lie algebra homomorphism  $\aa\to\k$ need not necessarily be
 $\Ad A$-invariant).

Given $\chi\in \BX(\aa)$,
 we 
 write $\chi\tpp$ for the function
$\aa\to \k,\,x\mapsto   \chi(x^{[p]})$.

\begin{lemma}\label{le} For any $\chi\in \BX(\aa)$,
the function $\chi\tpp: \aa\to\k$ is  a $p$-linear map,
that is, $\chi\tpp\in\atw$. 
\end{lemma}

\begin{proof} The above mentioned Jacobson's formula implies that,
in $\Ua$, one has an equality, cf. also \cite[Lemma 2.1]{Ja}:
$$ (x+y)^{[p]}-x^{[p]}-y^{[p]}=(x+y)^p-x^p-y^p.
$$
Now, 
extend $\chi$ to an algebra
homomorphism $\Ua\to\k$, and apply the resulting map to the equation
above.
We find
\begin{align*}
\chi\bigl((x+y)^{[p]}-x^{[p]}-y^{[p]}\bigr)&=
\chi\bigl((x+y)^p-x^p-y^p\bigr)=\chi(x+y)^p-\chi(x)^p-\chi(y)^p\\
&=
\bigl(\chi(x)+\chi(y)\bigr)^p-\chi(x)^p-\chi(y)^p=0.
\end{align*}
The Lemma follows.
\end{proof}

Lemma \ref{le} shows that the assignment $\chi\mapsto \chi\tpp$,
as well as the following assignment
\beq{as}
\vkb:\ \BX(\aa)\map\atw,\quad
\vkb(\chi)=\chi\tw-\chi\tpp:\ 
x\mto \chi(x)^p-\chi(x\tp),
\eeq
 gives a well-defined
$p$-linear map $\BX(\aa)\map\atw$.
 The map \eqref{as} will play an important role later in this paper,
it may be thought of as a Lie algebra analogue of the
{\em Artin-Schreier map}.

Let $\I_\varphi\sset\sym\aa\tw=\k[\atw]$
denote the
  maximal ideal  corresponding to a point~$\varphi\in\atw.$ 

Now, given $\chi\in  \BX(\aa)$,
 extend it
  to an algebra map $\chi: \Ua\to\k,$ as in the proof of Lemma \ref{le}.

\begin{corollary}\label{compositeIJ} Let $\chi\in  \BX(\aa)$. For
the  composite homomorphism below we have
$$
\Ker\big[\xymatrix
{\sym\aa\tw
\;\ar[rr]^<>(0.5){\bz_{_\U}}_<>(0.5){\eqref{pcenter}}&&
\Ua\ar[r]^<>(0.5){\chi}&\k}\big]=\I_{\as(\chi)}.
$$ 
\end{corollary}
\begin{proof} For $x\in\aa$, we compute
$\dis \chi(\bz_{_\U}(x))=\chi(x^p-x\tp)=\chi(x)^p-\chi(x\tp)=\as(\chi)(x).$
\end{proof}

Let
$\BX(A):=\Hom(A,\Gm)$ be the character lattice
of the algebraic group $A$. 
The differential of a character
${f}: A\to \Gm$ at $1\in A$  is
a  linear function $x\mapsto x({f})(1)$ on the
Lie algebra $\aa$, which is clearly an element
of $\BX(\aa).$ We denote this linear function
by $\dlog{f}$, so that the assignment ${f}\mapsto\dlog{f}$
yields
 an additive group homomorphism $\dlog: \BX(A)\to \BX(\aa).$

Fix $f\in  \BX(A)$ and put $\phi:=\dlog{f}\in\BX(\aa).$
Observe that, for any $x\in\aa$ viewed as a left invariant vector field
on $A$,  we have
$x({f})=\phi(x)\cdot{f}$.
It follows that the $p$-th power of $x$,
 viewed as a left invariant differential operator on
$A$ of order $p$, acts on a character
${f}\in \BX(A)$ as multiplication by the constant
$\phi(x)^p\in\k$.
On the other hand, this differential operator
is a derivation, which corresponds to
the  left invariant vector field $x^{[p]}$.
Thus, we also have $x^{[p]}({f})=\phi(x^{[p]})\cdot{f}$.
Combining together the equations above,  we deduce
$\phi(x\tp)=\phi(x)^p,\,\forall x\in\aa,$ that is,
$\vkb(\phi)=0$. Thus, we have proved
\beq{BX}
\dlog(\BX(A))\sset \Ker\big[\as:\ \BX(\aa) \map \atw\big].
\eeq

\subsection{Restricted enveloping algebras.} Fix $\chi\in\BX(\aa)$,
and let $\chi:\Ua\to\k$ be the corresponding algebra homomorphism.
\begin{defn}\label{IJ}
Let $I_\chi:=\Ker(\Ua\to\k)$ 
denote the kernel
of $\chi$, the
two-sided ideal in $\Ua$ generated by the elements
$\{x-\chi(x)\}_{x\in\aa}.$ 
Also, in $\z(\aa)$, consider the following ideal
$$I\tw_\chi:=I_\chi\cap \z(\aa)
\,\sset\,\z(\aa)\;\sset \Ua,\quad\text{and set}\quad
\uchi:=\Ua/\Ua\cd I\tw_\chi.
$$
\end{defn}
Here, $\Ua\cd I\tw_\chi\sset I_\chi$, is 
 an $\Ad A$-stable two-sided ideal
in $\Ua$, and the quotient $\uchi$ is an associative
 algebra of dimension $\dim\uchi=p^{\dim\aa}$, 
 called $\chi$-{\em restricted}
enveloping algebra. By definition
there is an exact sequence
\beq{ug}
0\too\Ua\cd I\tw_\chi\too\Ua\stackrel{r}\too \uchi\too 0.
\eeq

Corollary \ref{compositeIJ} shows that $I\tw_\chi$ is a maximal
ideal in $\z(\aa)$ that goes, under the isomorphisms
$\z(\aa)\cong\sym\aa\tw\cong\k[\atw]$, to
the  maximal ideal in $\k[\atw]$ corresponding to the point
$\vkb(\chi)\in\atw,$ that is, to the ideal
$\I_{\vkb(\chi)}$. Thus, we have
\beq{ass}
I\tw_\chi=\bzu(\I_{\vkb(\chi)}),\quad\text{hence}\quad
\uchi=\Ua/\Ua\cd \bzu(\I_{\vkb(\chi)}).
\eeq

Observe further that,
the $A_1$-action on $\z(\aa)$ being trivial,
it  preserves the ideal $\Ua\cd I\tw_\chi$,
hence induces a well-defined $A_1$-action on
$\uchi$ by algebra automorphisms.
We set 
\beq{fi}
\fri_\chi\,:=\,r(I_\chi)\,=\,I_\chi/\Ua\cd I\tw_\chi\;\sset\;\uchi.
\eeq
Thus, $\fri_\chi$ is an $A_1$-stable two-sided ideal
in $\uchi$ generated by the elements
$\{x-\chi(x)\}_{x\in\aa}.$

In the special case $\chi=0$,
the restricted enveloping algebra $\fu_0(\aa):=\Ua/\Ua\cdot I\tw_0$
inherits from $\Ua$ the structure of a Hopf algebra. This Hopf
algebra is  dual to
$\k[A_1]$, the coordinate ring of  the Frobenius kernel $A_1$.

  For any $\chi$, the map $\ad x: u\mapsto x\cdot u - u\cdot x,
\,x\in\aa,u\in\uchi,$
extends to a well-defined $\fu_0(\aa)$-action on $\uchi$,
that is, to an algebra map
$\ad: \fu_0(\aa)\map \End_\k\bigl(\uchi\bigr)$. This  $\fu_0(\aa)$-action
corresponds to the adjoint action on $\uchi$ of the
Frobenius kernel $A_1$.

\subsection{Quantum Hamiltonian reduction.}\label{quant_red}
Let $\DD$ be any associative, not necessarily commutative,
$\k$-algebra
equipped with an algebraic
action of the group $A$ by algebra automorphisms
and with an $A$-equivariant algebra map $\rho: \Ua\to \DD$
such that the adjoint $\aa$-action on $\DD$, given by
$\ad x: u\mapsto \rho(x)\cdot u - u\cdot \rho(x),\, x\in\aa, u\in\DD,$ is equal
to the differential of the $A$-action.

Let  $I\sset \Ua$ be an  $\Ad A$-stable
 two-sided
ideal.
Then,
 $\DD\cd \rho(I)$ is an $A$-stable
left ideal  in $\DD$. 
 It is easy to verify that
 multiplication in $\DD$
descends to a well-defined associative
algebra structure on $(\DD/\DD\cd \rho(I))^A$, the space
of $A$-invariants in $\DD/\DD\cd \rho(I)$.

Abusing notation, from now on we will write
$\DD\cd I$ instead of $\DD\cd \rho(I)$.

\begin{rem} The algebra  $(\DD/\DD\cd I)^A$ may be thought of
as a `Hamiltonian reduction' of $\DD$ with respect to
$I$.
\erem

Observe also that, if $u\in\DD$ is such that
$u\,\text{mod}(\DD\cd I)\,\in (\DD/\DD\cd I)^A$, then
the operator of right multiplication by $u$
descends to a well-defined map
$R_u: \DD/\DD\cdot I\to\DD/\DD\cdot I$.
Moreover, the assignment $u\mapsto R_u$
induces an algebra isomorphism
$\dis (\DD/\DD\cd I)^A\iso \bigl(\End_\DD(\DD/\DD\cd I)\bigr)^{\op{opp}}.
$

More generally, let $M$ be a left $\DD$-module equipped
with an $A$-equivariant structure (i.e., such that the
action map $\DD\otimes M\to M$ is  $A$-equivariant).
The algebra map $\Ua\to\DD$ makes $M$ 
an $\Ua$-module. The space
$(M/I\cdot M)^A$  acquires a natural left
$(\DD/\DD\cdot I)^A$-module structure,
to be called a Hamiltonian reduction of $M$. 
Similar construction applies to right $\DD$-modules.

Next, fix $\chi\in\BX(\aa)$, and let
$\uchi$ be the corresponding $\chi$-restricted enveloping algebra.
Recall that this algebra comes equipped with the
adjoint action of $A_1$,  the Frobenius kernel.
Let $\DD$ be  an associative algebra equipped with
$A_1$-action, and with $\rho: \uchi\to\DD$,
an $A_1$-equivariant algebra morphism.
One shows similarly that, given an
$A_1$-stable two-sided ideal $\fri\sset \uchi$,
there is a natural associative algebra
structure on $(\DD/\DD\cdot\fri)^{A_1}$.
In the special case $\fri=\fri_\chi$, see \eqref{fi}, the algebra
$(\DD/\DD\cdot\fri_\chi)^{A_1}$
may be thought of as a quantum hamiltonian reduction
of $\DD$ with respect to the action of
$A_1$, an `infinitesimal' group-scheme.

Later on, we will be interested in the following
special
case of this construction.
Let $E$ be a finite dimensional
$A$-module such that the induced action map $\rho: \Ua\to\End E$
descends to the algebra $\uchi$, i.e., vanishes on
the ideal $I\tw_\chi\sset\z(\aa)$.
Write
${E_\chi:=\{e\in E\mid} {x(e)=\chi(x)\cdot e,\,\forall x\in \aa\}}$
for the  $\chi$-weight space of $\Ua$.

We put  $\DD:=\End E\, (=\End_\k E)$.
We claim that  the  quantum Hamiltonian reduction
of the algebra $\DD$ with respect to the $\uchi$-action
is canonically isomorphic
to $\End_\k E_\chi$.

In more detail, we form an associative 
algebra $\DD_\chi:= (\DD/\DD\cdot\fri_\chi)^{A_1}$.
The  natural action of $\DD$ on $E$
descends to a well-defined $\DD_\chi$-action on the
weight space $E_\chi$. On the other hand,
 we have 
 a right $\DD_\chi$-action on
$\DD/\DD\cdot\fri_\chi$.

We are going to establish
canonical algebra isomorphisms
\beq{chipsi}\End_\DD(\DD/\DD\cdot\fri_\chi)^{\op{opp}}
\,\stackrel{\mathbf{l}}\cong\,
\DD_\chi\,
\stackrel{\mathbf{r}}\cong\,
\End_\k E_\chi.
\eeq

The isomorphisms above follow from a more general result below
that involves two characters $\chi,\psi\in\BX(\aa)$, such that
$\vkb(\chi)=\vkb(\psi)$. 
In this case,  in $\Ua$ we have $I\tw_\chi= I\tw_\psi$.
Hence there is a canonical identification
$\uchi=\upsi$, and we may view the $A$-representation $E$
either as an  $\uchi$- or as an $\upsi$-module.
We set $\DD:=\End_\k E$, as above,
and consider the Hamiltonian reductions
$\DD_\chi$ and $\DD_\psi$.
We have a $\DD_\chi$-action on the
weight space 
$E_\chi\sset E$,
and a  similar $\DD_\psi$-action on the
weight  space $E_\psi$.

Observe that the natural $\DD$-bimodule structure
on the algebra $\DD$, via left and right multiplication,
descends to a $\DD\dash\DD_\chi$-bimodule structure
on the vector space
$\Hom_\k(E_\chi,E)$. Further,
the right $\DD_\chi$-action on $\DD/\DD\cdot\fri_\chi$,
resp., $\DD_\psi$-action on $\DD/\DD\cdot\fri_\psi$,
gives  the following space
a natural $\DD_\chi\dash\DD_\psi$-bimodule structure
\beq{bimod}
{_\psi\DD_\chi}:=\Hom_{\DD}\bigl(\DD/\DD\cd\fri_\chi,\,
\DD/\DD\cd\fri_\psi\bigr)^{\psi-\chi},
\eeq
where the superscript `$\psi-\chi$' denotes
the $(\psi-\chi)$-weight component of the
natural (adjoint) $A$-action.

\begin{lemma}\label{weight} \vi The restriction
map
$\End E\to \Hom(E_\chi,E),$ resp. $\End E\to \Hom(E_\psi,E),$
induces a  $\DD\dash\DD_\chi$-bimodule,  resp. 
$\DD\dash\DD_\psi$-bimodule, isomorphism
$$\DD\big/\DD\cd\fri_\chi\overset{\op{res}}\iso\Hom_\k(E_\chi,E),\quad
\text{resp.},\enspace \DD\big/\DD\cd\fri_\psi\overset{\op{res}}
\iso\Hom_\k(E_\psi,E).
$$

\vii We have
the following $\DD_\chi\dash\DD_\psi$-bimodule
isomorphisms
\beq{compo}
{_\psi\DD_\chi}\overset{\op{res}}\iso \Hom_{\DD}\bigl(\Hom_\k(E_\chi,E),
\Hom_\k(E_\psi,E)\bigr)^{\psi-\chi}
\,\stackrel{_\sim}{\longleftarrow}\,
\Hom_\k(E_\chi, E_\psi).
\eeq

\noindent
For $\chi=\psi$, it follows that the  map $\mathbf{l}$
in \eqref{chipsi} is an algebra isomorphism, furthermore,
the composite map in \eqref{compo} induces
an {\sl algebra} isomorphism
$\DD_\chi={_\chi\DD_\chi}\iso \End_\k(E_\chi)$,
which gives the  isomorphism $\mathbf{r}$ in \eqref{chipsi}.
\end{lemma}

\begin{proof}[Sketch of Proof.] Associated to a vector subspace
$F\sset E$, one has a left ideal
$J_F\sset\DD$, defined by
$J_F:=\{f\in \DD=\End_\k E\enspace\big|\enspace f|_F=0\}$. Moreover, any
left ideal $J\sset \DD$  has the form
$J=J_F$ where the corresponding
 subspace $F\sset E$ can be recovered from $J$
by the formula $F=\bigcap_{f\in J}\,\Ker f$.
Applying this to the left ideal
$J=\DD\cdot\fri_\chi$ we get
$\DD\cdot\fri_\chi=\{f\in \End_\k E\mid f|_F=0\}$,
where $F=\bigcap_{x\in \fri_\chi}\,\Ker \rho(x)$.
The latter space equals $E_\chi$, by definition.
Thus, we deduce 
$\DD\cdot\fri_\chi=\{f\in \End_\k E\mid f|_{E_\chi}=0\}=\Hom_\k(E/E_\chi,E)$,
hence $\DD/\DD\cdot\fri_\chi\cong\Hom_\k(E_\chi,E)$.

The rest of the proof is an elementary exercise
which we leave for the reader.
\end{proof}

\subsection{Moment maps.}\label{group actions}
Let $A$ be a linear algebraic group as in sect.\ref{not_sect},
 and let $A\times X \to X$ be an algebraic action
on $X$, a smooth $\k$-variety.
Any element $x\in \aa$
gives rise to an algebraic vector  field $\xi_x$  on $X$.
We may view $\xi_x$ as a regular function on $T^*X$.
This way, the assignment $x\mapsto \xi_x$ extends
uniquely to an $A$-equivariant
 Poisson algebra map $\mu_{_{\op{alg}}}: \sym\aa\to \k[T^*X]$.
Since $\sym\aa\cong\k[\aa^*]$, this algebra map
induces an $A$-equivariant
 morphism $\mu: T^*X\to\aa^*$, called  {\em moment map},
such that the algebra map $\mu_{_{\op{alg}}}$ becomes the pull-back via  $\mu$.

There is also a noncommutative analogue of
the  Poisson algebra map $\mu_{_{\op{alg}}}.$
Specifically, the   Lie algebra morphism $x\mapsto \xi_x$
 extends uniquely
 to an $A$-equivariant associative algebra homomorphism
$\mu_{_{\U}}: \Ua\to \D_X$
(more generally, given  an $A$-equivariant
vector bundle $\CL$,
one defines similarly
 an associative algebra homomorphism
$\mu_{_{\U}} :\Ua\to \D(X,\CL)$).
The morphism $\mu_{_{\U}}$ is compatible with  natural filtrations,
hence,  induces a canonical
{\em graded} algebra  homomorphism $\mu_{_\rr}: \rees\Ua\map\rees\D_X$.
We may view the later homomorphism
as a map $\rees\Ua\map\pi\idot\rr\D\tw,$
cf. Sect. \ref{rrDtw}.

We have the following commutative diagram:
\beq{U_square}
\xymatrix{
\sym\aa\tw
\ar@{^{(}->}[dd]_<>(0.5){\eqref{pcenter}}^<>(0.5){\bzur:\,
x\mapsto x^p-t^{p-1}\cdot x\tp}
\ar[rrr]^<>(0.5){\mu\tw_{_{\op{alg}}}}&&&\z_{X\tw}=\sym \TT_{X\tw}=
\pi\idot\oo_{_{T^*X\tw}}
\ar@{^{(}->}[dd]_<>(0.5){\eqref{central_rees}}^<>(0.5){\bzdr:\,
\xi\mapsto \xi^p-t^{p-1}\cdot \xi\tp}\\
\\
\rees\Ua\ar[rrr]^<>(0.5){\mu_{_\rr}}&&&\rees(\Fr\idot\D_X)=\pi\idot\rr\D\tw.
}
\eeq 
The map $\mu_{_\rr}$ in the bottom row specializes  at $t=1$
to the map $\mu_{_\rr}|_{t=1}=\muu: \Ua\to \D_X$ considered earlier,
and specializes at $t=0$ to the map $\mu_{_\rr}|_{t=0}=
\mu_{_{\op{alg}}}$.
The vertical maps in  diagram \eqref{U_square} 
 are the \gre algebra morphisms considered earlier.

Further, fix  $\chi\in\BX(\aa)$,
let   $I\tw_\chi\sset\z(\aa)$  be the corresponding ideal,
cf. Definition \ref{IJ},
and $\muu(I\tw_\chi)\sset \z_{X\tw}$ its image
in $\Fr\idot\D_X=\pi\idot\D\tw.$
 From commutativity of diagram \eqref{U_square}
for $t=1$ and formula \eqref{ass}, we deduce that
 the canonical isomorphism
$\z_{X\tw}\iso\pi\idot\oo_{_{T^*X\tw}}$
takes the ideal $\muu(I\tw_\chi)$ 
to  the ideal  $\mu_{_{\op{alg}}}(\I_{\vkb(\chi)})\sset \pi\idot\oo_{_{T^*X\tw}}$.
On the other hand, it follows from
Poincar\'e-Birkhoff-Witt theorem that the associated graded
ideal
$\gr I\tw_\chi$ equals the
 the augmentation ideal in $\gr\z(\aa)=\sym\aa\tw$. Thus, specializing
diagram  \eqref{U_square} 
at $t=0$ and $t=1$, respectively, we find
\beq{gr_mu}
\muu(I\tw_\chi)=\mua(\I_{\as(\chi)}),\quad\text{resp.},
\quad
\gr\bigl(\muu(I\tw_\chi)\bigr)=
\mua(\I_0).
\eeq

We introduce the following subscheme in $T^*X\tw$:
\beq{mutw}
\ttw:=\text{zero scheme of $\muu(I\tw_\chi)$}=
\text{zero scheme of $\mu_{_{\op{alg}}}(\I_{\vkb(\chi)})$}
=\big[\mu\tw\big]\inv(\vkb(\chi)),
\eeq
the scheme-theoretic fiber 
 of the
moment map $\mu\tw: T^*X\tw\to\atw$
over the point $\vkb(\chi)$, cf. \eqref{as}.
For example, if $f: A\to\Gm$ is an algebraic group character
and $\chi:=\dlog f$, then we have $\as(\chi)=0,$ see
\eqref{BX}. Hence, $\ttw
=\big[\mu\inv(0)\big]\tw.$

\subsection{Hamiltonian reduction of differential operators.}
\label{ham_diff}
We keep the setup of Sect.~\ref{group actions}.
Thus, we have an $\Ad A$-invariant  homomorphism $\chi: \aa\to\k,$
 the corresponding ideal $I\tw_\chi\sset\z(\aa)$,
and  its image  $\muu(I\tw_\chi)\sset \pi\idot\D\tw$
under  the  map
$\muu: \Ua\to \pi\idot\D\tw$,
induced by the $A$-action.
This image is a central subalgebra
in  the  Azumaya algebra $\D\tw$
 on $T^*X\tw$, hence,
$\D\tw\cdot \pi\hdot \muu(I\tw_\chi)\sset \D\tw$
is an $A$-stable two-sided ideal.

We put $
\D\tw_\chi:=
\D\tw/\D\tw\cdot\pi\hdot \muu(I\tw_\chi).
$ Since
$\pi\hdot \muu(I\tw_\chi)=\pi\hdot\mua(\I_{\vkb(\chi)}),
$ by formula \eqref{gr_mu}, we have 
\beq{mtw}\D\tw_\chi=\D\tw/\D\tw\cdot\pi\hdot \muu(I\tw_\chi)
=\D\tw/\D\tw\cdot\pi\hdot\mua(\I_{\vkb(\chi)})
=\D\tw\big|_{\ttw}.
\eeq
We see that $\D\tw_\chi$ is a restriction of the sheaf $\D\tw$ to $\ttw$,
the scheme-theoretic fiber of the moment map, see
\eqref{mutw}.
Thus, $\D\tw_\chi$
is a coherent sheaf of associative 
algebras on the subscheme $\ttw$.
By construction, the 
map $\Ua\to\D(X)=\G(T^*X\tw,\D\tw)$ descends, in view of
 exact sequence  \eqref{ug}, to
an $A_1$-equivariant algebra homomorphism
\beq{act_chi}
\rho_\chi: \ \uchi\too \G\bigl(T^*X\tw,\,
\D\tw_\chi\bigr).
\eeq

Recall the two-sided ideal $\fri_\chi\sset\uchi,$ see
\eqref{fi}, and put
\beq{Dchi}
\ee_\chi\,:=\,\bigl({\D\tw_\chi}\big/\D\tw_\chi \cd\fri_\chi\bigr)^{A_1}.
\eeq
This is a sheaf  on  $\ttw$
that may be thought of as a
 Hamiltonian reduction of the algebra $\D\tw_\chi$ with respect to the
action of $A_1$.
The construction of Sect. \ref{quant_red} applied
to $\DD:=\D\tw_\chi$ and to the homomorphism \eqref{act_chi},
gives $\ee_\chi$  
the natural structure
of a coherent sheaf of 
associative algebras
on the scheme 
${\ttw}\sset T^*X\tw$.
Observe further that the $A$-action 
on $\D\tw$ factors, when restricted to $A_1$-invariants,
through the quotient  $A\tw=A/A_1$. 
Thus, the  sheaf ${\ee_\chi}$
acquires an $A\tw$-equivariant structure.

On the other hand,  rather than 
 performing the Hamiltonian reduction
of $\D\tw_\chi$ with respect to the $A_1$-action,
one may   perform the Hamiltonian reduction
of  $\D\tw$, a larger object, with respect to the
action of $A$, a larger group, that is,
to consider $A$-invariants in
$\D\tw\big/\D\tw\cd  I_\chi$
(abusing notation,
 we will write $\D\tw\cd I_\chi$ instead of $\D\tw\cd \pi\hdot\mu_{_\U}(I_\chi)$, and use similar
notation in other cases).

The  elementary result below is a manifestation of the general principle
saying that Hamiltonian reduction can be performed in stages:
to make a  Hamiltonian reduction with respect to
 $A$, one can first perform   Hamiltonian reduction with
respect to
$A_1$, and then make reduction with respect to 
$A\tw=A/A_1$.

\begin{lemma}\label{stages} There is a 
canonical algebra isomorphism
$$
\G(X,\,\D_X/\D_X\cd  I_\chi)^A
\cong\,\G\bigl(\ttw,\,{\ee_\chi}\bigr)^{A\tw}.
$$
\end{lemma}
\begin{proof} It is clear that in $\D\tw/\D\tw\cd I\tw_\chi$, 
one has an equality
$\D\tw\cdot I_\chi\big/\D\tw\cd I\tw_\chi
=\D\tw_\chi\cdot \fri_\chi$.
Thus, we obtain
$
\D\tw\big/\D\tw\cd I_\chi\cong
{\D\tw_\chi}\big/{\D\tw_\chi}\cd\fri_\chi.
$
Taking $A_1$-invariants on both sides, we deduce
an isomorphism (of sheaves of associative algebras on
$\ttw$):
\beq{ham_tw}
\bigl(\D\tw\big/\D\tw\cd I_\chi\bigr)^{A_1}
\,\cong\,\bigl(\D\tw_\chi\big/{\D\tw_\chi}\cd\fri_\chi\bigr)^{A_1}=
{\ee_\chi}.
\eeq
Applying the functor $\G(T^*X\tw,-)^{A\tw}$ to 
 $A\tw$-equivariant sheaves in
 \eqref{ham_tw},  we obtain a chain of
 canonical algebra isomorphisms
 \begin{align*}
\G(X,\,\D_X/\D_X\cd I_\chi)^A
&\cong \G\bigl(T^*X\tw,\,\D\tw/\D\tw\cd I_\chi\bigr)^A\cong
\G\bigl(T^*X\tw, \,(\D\tw/\D\tw\cd I_\chi)^{A_1}\bigr)^{A\tw}_{_{}}\\
\text{(by \eqref{ham_tw})}\quad&\cong
\G\bigl(T^*X\tw, \ee_\chi\bigr)^{A\tw}
\cong\G\bigl(\ttw,\,
{\ee_\chi}\bigr)^{A\tw}.
\end{align*}
The Lemma follows.
\end{proof}

\subsection{The case of free $A$-action.}\label{free}
Keep the notation of \S\ref{ham_diff} and let $\chi=0$,
hence $I_\chi=I_+\sset \Ua$
is  the
augmentation ideal. Write $\D_X\cdot I_+=\D_X\cdot \aa$ for the left
ideal generated by  the image of $I_+$ under the homomorphism
$\mu_{_{\U}}: \Ua\to \D(X)$.
Applying the construction of section \ref{quant_red}
to the algebra $\DD:=\D_X$ and 
the two-ideal $I_+\sset\Ua$
one gets an associative 
 algebra 
$\bigl(\D_X/\D_X\cdot I_+\bigr)^A$.

Assume now that the $A$-action on
$X$ is free and, moreover,  there is a smooth 
variety $Y$,
and a smooth {\em universal geometric quotient} morphism $\pr_Y: X\onto
Y$ 
(whose fibers
are exactly the $A$-orbits), see \cite[Definition 0.7]{GIT}.
It is well-known that the algebra of differential operators
on $Y$ can be expressed in terms of differential 
operators on $X$ as follows
\beq{XY}\D_Y\cong 
\bigl((\pr_Y)\idot(\D_X/\D_X\cdot I_+)\bigr)^A=
\bigl((\pr_Y)\idot(\D_X/\D_X\cdot\mu_{_{\U}}(\aa))\bigr)^A.
\eeq

More generally, fix an algebraic homomorphism
$\chi: A\to\Gm$.
Given a free $A$-action on $X$,
let $\oo_Y(\chi)$
be an invertible sheaf  on $Y$
defined as the subsheaf of $(\pr_Y)\idot\oo_X$
formed by the functions $f$ such that
$a^*(f)=\chi(a)\cdot f,\,\forall a\in A.$
Let $\D_Y(\oo_Y(\chi))$ be the corresponding sheaf
of {\em twisted differential operators},
and
 $\D(Y,\chi):=\G(Y,\,\D_Y(\oo_Y(\chi))$ the
algebra of its global sections.

There is a $\chi$-twisted version of formula \eqref{XY} that
provides a canonical isomorphism
$\dis{\bigl((\pr_Y)\idot(\D_X/\D_X\cd I_\chi)\bigr)^A}
\iso \D_Y(\chi)$
(isomorphism of sheaves of algebras on $Y$).
Taking global sections on each side of
the isomorphism,
we get algebra isomorphisms:
\beq{XA}
\xymatrix{
\G\bigl({\ttw},\,{\ee_\chi}\bigr)^{A\tw}\;
\ar[rrr]^<>(0.5){{\text{Lemma \ref{stages}}}}_<>(0.5){\sim}&&&
\;\G\bigl(X,\,\D_X/\D_X\cd I_\chi\bigr)^A}\cong\,
\D(Y,\chi).
\eeq

The isomorphism above makes sense, in effect, not only for
$\chi\in\BX(A),$
but also
in a slightly more general setting where
$\chi\in\BX(\aa)$  is an $\Ad A$-invariant Lie algebra character
that does not necessarily  exponentiate
to an algebraic group homomorphism
$A\to\Gm$.
Although, generally, the sheaf  $\oo_Y(\chi)$
is not defined in such a case,
 the corresponding sheaf $\D_Y(\chi)$
of {\em twisted differential operators} 
is always well-defined, cf. [BB],
and the isomorphism in \eqref{XA} still holds.

\begin{rem}
The algebra $\D(Y,\chi)$ may be thought of as a quantization
of the commutative algebra $\k[\mu\inv(\chi)]^A$, 
the coordinate ring of the Hamiltonian reduction
of $T^*X$ with respect to the 1-point
orbit $\{\chi\}\sset\aa^*$ and the moment map
$\mu: T^*X\to\aa^*$.
\erem

\section{Azumaya algebras via Hamiltonian reduction}
\subsection{The main result.}
\label{some}
Let $X$ be a smooth $A$-variety. Below, we are going to
extend considerations of  section \ref{free}
to a
more general case
where the $A$-action on $X$ is not necessarily free,
but the corresponding Hamiltonian $A$-action on
$T^*X$ is free on an open subset of $T^*X$.

There is a natural action of the multiplicative group $\Gm$
on  the vector space $\aa^*$ and also on the vector
bundle $T^*X$, by dilations.
The moment map $\mu: T^*X\to\aa^*$ is clearly
compatible with these two actions.
It is also compatible with the $A$-actions, and the latter
 commute with the $\Gm$-actions.
Thus, $\mu$ is an
equivariant morphism between  $\Gm\times A$-varieties.

Let $\mu\inv(0)\sset T^*X$ be
 the scheme-theoretic zero fiber  of the moment map.
This is clearly
 a $\Gm\times A$-stable subscheme in $T^*X$.

  From now on, we make the following

\begin{assumption}\label{assume}
There is a Zariski open $\Gm\times A$-stable 
subscheme $\M\sset \mu\inv(0)$
which is a reduced smooth locally-closed connected subvariety
in $T^*X$ such that\hfill\newline
\pb{The differential of the moment map $\mu: T^*X\to\aa^*$
is surjective at any point of $\M$;}\hfill\newline
\pb{The
$A$-action on $\M$ is free,  moreover,
there is a smooth  variety $\mm$ and 
a smooth universal geometric quotient morphism $\M\to\mm$
(in particular, it is a
principal $A$-bundle 
whose fibers are precisely the $A$-orbits in $\M$), see \cite{GIT};}\hfill\newline
\pb{The natural $\Gm$-action 
 on the algebra $\k[\mm]$ (arising from the $\Gm$-action on
$\M$) has no
negative weights, more geometrically,
the induced  $\Gm$-action on the scheme
$\mm_{\op{aff}}:=\Spec\k[\mm]$, the {\em affinization}
of $\mm$, is an  {\em attraction}.}
\pb{The canonical projection
$f: \mm\to \mm_{\op{aff}}$ is a proper morphism.}
\end{assumption}

The assumptions above insure that the standard symplectic
structure on  $T^*X$ induces a   symplectic
structure on  $\mm$. Thus, the manifold $\mm$ may be thought of as
a Hamiltonian reduction of $T^*X$ at $0$.

Now,
 let $\chi\in\BX(\aa)$ be  
such that $\vkb(\chi)=0$. Recall the
notation $\ttw:=[\mu\tw]\inv(\vkb(\chi))$.
The equation  $\vkb(\chi)=0$
implies that $\ttw=[\mu\inv(0)]\tw$, which 
 is clearly a $\Gm\times A\tw$-stable
 subscheme in $T^*X\tw$. Further, by
the Basic Assumptions, the scheme
$[\mu\inv(0)]\tw$ contains $\M\tw$,
the Frobenius twist of $\M$, as  
 an open subscheme. Thus, for any $\chi\in\BX(\aa)$ 
such that $\vkb(\chi)=0$, the Basic Assumptions yield
 the following diagram
\beq{in_to}
\xymatrix{
\ttw=[\mu\inv(0)]\tw\enspace&&
\enspace\M\tw\enspace\ar@{_{(}->}[ll]
\ar@{->>}[rrrr]^<>(0.5){\varpi}_<>(0.5){\text{principal $A^{(1)}$-bundle}}
&&&&
\enspace\mm\tw.
}\eeq

We are going to define a coherent sheaf $\az_\chi$ on  $\mm\tw$
that  will be an Azumaya
$\oo_{\!\!_{\mm\tw}}$-algebra of degree $p^{1/2\dim \mm}$, to be
called the {\em quantum hamiltonian reduction of $\D_X$ at $\vkb(\chi)$}.
To this end, 
we restrict  $\ee_\chi=(\D\tw_\chi/\D\tw_\chi\cdot \fri_\chi)^{A_1}$, an
$A^{(1)}$-equivariant 
sheaf on
$\ttw\sset T^*X^{(1)}$, cf. \eqref{Dchi},  to the  open
subset $\M^{(1)}$, and consider the push-forward of that restriction
under the map
$\varpi:\M\tw\map\mm\tw$, cf. \eqref{in_to}.

\begin{defn}\label{az_defn} We define the following
coherent
sheaf
of associative algebras on $\mm\tw$:
$$\az_\chi:=\,
\varpi\idot\bigl(\ee_\chi\big|_{\M\tw}\bigr)^{A\tw},
\quad\text{and put}\quad \A_\chi:=\G(\mm\tw,\az_\chi).
$$
\end{defn}

Assume next that we are given two  points
$\chi, \psi\in\BX(\aa),$ such that
 the character ${\chi-\psi:}$
$\aa\to\k$ can be exponentiated 
to a group homomorphism $A\to\Gm$, i.e., such that  $\chi-\psi\in \dlog(\BX(A)$).
Then, formula \eqref{BX} yields $\vkb(\chi)=\vkb(\psi)$.
Hence we have
$\ttw=\ttwp$, and we may view
the set $\mm$ as a geometric quotient of an open subset of either
$\ttw$ or~$\ttwp$.

The result below that will play a key role
in subsequent sections,
is   a generalization of 
\cite[Proposition 4.8]{bk}.

\begin{theorem}\label{roma} Let $\chi\in\aa^*$ be an $A$-fixed point
 such that $\vkb(\chi)=0$
and such that the Basic Assumptions \ref{assume} hold.  Then we have

\vi The sheaf
$\az_\chi$ is a sheaf of Azumaya algebras on $\mm\tw$
equipped with a canonical algebra morphism
$\dis\Xi_\chi:\
\G(X,\,\D_X/\D_X\cdot I_\chi)^A\map\G(\mm\tw,\az_\chi).$

\vii For  all $i>0,$ we have
$ H^i(\mm\tw, \az_\chi)=0$.

\viii If the algebra $\A_\chi=\G(\mm\tw,\az_\chi)$ has finite
homological
dimension then the  functors below give mutually inverse
equivalences 
of bounded derived categories of  sheaves of coherent $\az_\chi$-modules and finitely
generated $\A_\chi$-modules, respectively:
$$
\xymatrix{
D^b(\az_\chi\dash\on{Mod})\enspace
\ar@<1ex>[rrrr]^<>(0.5){{\scr L}\,\longrightarrow\, \RG(\mm\tw, {\scr L})}&&&&
\enspace
D^b(\A_\chi\dash\on{Mod}).
\ar@<1ex>[llll]^<>(0.5){{\az_\chi\stackrel{_L}{\otimes}_{\A_\chi} L\, \longleftarrow\, L}}
}
$$

\iv Let $\psi \in\BX(\aa)$ be another  point satisfying all the assumptions
above
and
such that  $\chi-\psi\in \dlog(\BX(A))$. 
Then the corresponding
Azumaya algebras $\az_\chi$ and $\az_\psi$ are Morita equivalent
(but not necessarily isomorphic).
\end{theorem}

The rest of this section is devoted to the
proof of the Theorem.

As will be explained  in Sect. \ref{proof_roma} below, part (iii) of Theorem
\ref{roma} is entirely due to
\cite[Proposition~2.2]{bk}. A result similar to
part (iv) of the theorem is also contained in an
updated version of \cite[\S2.3.1]{BMR}.

\begin{rem} \vi The assumption of the Theorem that $\vkb(\chi)=0$
may be relaxed, as will be explained elsewhere.

\vii We will show, in the course of the proof of Theorem
\ref{roma}, that the sheaf $\ee_\chi\big|_{\M\tw}$
is also an Azumaya algebra, specifically, we have
an Azumaya algebra isomorphism:
\beq{Efree}
\ee_\chi\big|_{\M\tw}=\varpi^*\az_\chi.
\eeq

\viii It will also follow from the proof that
the sheaf $\az_\chi$, viewed
as a vector bundle on $\mm\tw$, is a deformation of the
vector bundle $\Fr\idot\oo_\mm$; in particular,
in the Grothendieck group $K(\Coh(\mm\tw))$
on has $[\az_\chi]=[\Fr\idot\oo_\mm]$,
furthermore,
 $R^if\idot\tw\az_\chi=0$ for any $i>0$,
where
$f\tw: \mm\tw\to \mm\tw_{\op{aff}}$ denotes the affinization morphism.
\end{rem}

\subsection{Deformation construction.}\label{deform}
We are going to apply the Rees algebra formation
 to all the objects involved
in the construction of the algebra $\az_\chi$.

In more detail, the standard filtration on $\Ua$ induces a filtration on
the ideal $I_\chi\sset\Ua$, and  we  form a graded ideal 
$\rees I_\chi\sset \rees\Ua$, which is generated by
the set $\{{x-}{t\cdot\chi(x)\}_{x\in\aa}}$.
Further, let $\rr\z(\aa)$ be a
$\k[t]$-subalgebra in $\rees\Ua$ generated by the
image of the homomorphism $\bzur:\sym\aa\tw\to\rees\Ua,$
see
\eqref{pcenter}. Thus,  $\rr\z(\aa)\sset \rees\Ua$ is a graded central 
subalgebra, and $\rr{I\tw_\chi}:=\rr\z(\aa)\cap\rees I_\chi$ is
  a graded ideal in  $\rr\z(\aa)$. For any $x\in \aa$, 
we have $x^p-{t^{p-1}\cdot x\tp}$
$-t^p\cdot\as(\chi)(x)=\bzur(x)-t^p\cdot\as(\chi)(x)\in \rr\z(\aa).$
On the other hand, the following equations
\begin{align*}
x^p-t^{p-1}\cd x\tp&-t^p\cd\as(\chi)(x)
=\bigl(x^p-t^p\cd\chi(x)^p\bigr)-\bigl(t^{p-1}\cd x\tp-
t^p\cd\chi\tpp(x)\bigr)\\
&=
\bigl(x-t\cd\chi(x)\bigr)^p-t^{p-1}\cd\bigl(x\tp-t\cd\chi(x\tp)\bigr)\in\rees I_\chi
\end{align*}
show that $x^p-t^{p-1}\cdot
x\tp-t^p\cd\as(\chi)(x)\in \rr{I\tw_\chi}.$
Moreover, it is easy to verify that
the elements of this form generate  $\rr{I\tw_\chi}$
as an ideal.

Next, we apply the
Rees algebra construction to the filtered sheaf
$\Fr\idot\D_X$. Inside $\Fr\idot\D_X$,
we have two left ideals $\Fr\idot\D_X\cdot \muu(I\tw_\chi)\sset
\Fr\idot\D_X\cdot \muu(I_\chi),$ generated by
the images of the sets
$I\tw_\chi\sset I_\chi\sset\Ua$,
respectively, under the
moment map $\muu: \Ua\to\Fr\idot\D_X$.
The filtration on 
 $\Fr\idot\D_X$ induces by restriction natural
 filtrations on  $\Fr\idot\D_X\cdot \muu(I\tw_\chi)$
and $\Fr\idot\D_X\cdot \muu(I_\chi).$
Thus, we obtain
 graded ideals $\rees(\Fr\idot\D_X\cdot \muu(I\tw_\chi))\sset
\rees(\Fr\idot\D_X\cdot \muu(I_\chi))
\sset\rees(\Fr\idot\D_X).$

On the other hand, we have a
moment map $\mu_{_\rr}:\rees\Ua\to\rees(\Fr\idot\D_X),$
which takes the central subalgebra $\rr\z(\aa)\sset\rees\Ua$
into the central subalgebra $\rees\z_{X\tw}$
$\sset\rees(\Fr\idot\D_X),$
cf. diagram \eqref{U_square}.
Hence, the image of  $\rr{I\tw_\chi}$
is a subalgebra $\mu_{_\rr}(\rr{I\tw_\chi})\sset \rees\z_{X\tw},$
and we have
$\mu_{_\rr}(\rr{I\tw_\chi})\cd
\rees(\Fr\idot\D_X)\sset$
 $\rees(\Fr\idot\D_X\cd\muu(I\tw_\chi)),$
where the inclusion is strict, in general.

Further, we have a $\Gm$-equivariant sheaf $\rr\D\tw$ on
$\BA^1\times T^*X\tw$ corresponding
to the graded algebra $\rees(\Fr\idot\D_X)$.
The above constructed  graded ideals in the algebra  $\rees(\Fr\idot\D_X)$
give rise to the following three $\Gm$-equivariant sheaves of 
left ideals in 
 $\rr\D\tw$:
\begin{align}\label{ideals}
\rr\D\tw\cd\{x^p-t^{p-1}\cdot
x\tp-t^p\cdot\as(\chi)(x)\}_{x\in\aa}&=
\rr\D\tw\cdot\rr{I\tw_\chi}\\
&\sset
\rr(\D\tw\cdot I\tw_\chi)\sset
\rr(\D\tw\cdot I_\chi),\nonumber
\end{align}
where we follow our usual convention to drop the symbols $\muu$ and
$\mu_{_\rr}$ from the notation.

 By  \eqref{rees_gr}, we have a
 graded
algebra isomorphism
$\rees\z_{X\tw}\cong\oo_{_{\BA^1\times T^*X\tw}}.$
Thus,  $\mu_{_\rr}(\rr{I\tw_\chi})$ may be viewed as a
subset  in $\oo_{_{\BA^1\times T^*X\tw}}$, and 
from commutativity of diagram \eqref{U_square}
we deduce that the set
$\big\{\mu\tw_{_{\op{alg}}}(x)-t^p\cdot\as(\chi)(x)\big\}_{x\in\aa}$
 generates the ideal $\mu_{_\rr}(\rr{I\tw_\chi})\cdot \oo_{_{\BA^1\times T^*X\tw}}.$
Thus, in $\BA^1\times T^*X\tw$, we have
\beq{setI}
\text{Zero-scheme of}\enspace\;\mu_{_\rr}(\rr{I\tw_\chi})\,=\,
\{(t,\xi)\in \BA^1\times T^*X\tw\mid
\mu\tw(\xi)=t^p\cd\as(\chi)\}.
\eeq

The sheaf $\rr\D\tw/\rr\D\tw\cd\rr{I\tw_\chi}$ is clearly supported
on the subscheme \eqref{setI}.
We conclude that 
its quotient $\rr\D\tw/\rr(\D\tw\cd I_\chi),$
cf. \eqref{ideals}, 
is  supported on
the subscheme \eqref{setI} as well.

Now let $\chi\in\BX(\aa)$ be such that  $\vkb(\chi)=0$.
Then the set 
in \eqref{setI} reduces to a direct product
$\BA^1\times
[\mu\inv(0)]\tw$. 
One checks further, going through the identifications
used above, that
the $\Gm\times A\tw$-action on $\BA^1\times\ttw=\BA^1\times
[\mu\inv(0)]\tw$ arising from the
natural grading and  from the $A\tw$-action on  $\rr{I\tw_\chi}$,
respectively,
is the one
where the group
$\Gm$ 
acts diagonally, and the group $A\tw$ acts only on the factor
$[\mu\inv(0)]\tw$.

Recall the sheaf $\rr\D\tw/\rr(\D\tw\cdot I_\chi)$,
which  is supported on
 $\BA^1\times
[\mu\inv(0)]\tw$, since $\as(\chi)=0$. We restrict this  sheaf
 to the open subscheme
$\BA^1\times\M\tw
\sset \BA^1\times
[\mu\inv(0)]\tw$. 
We have the principal $\Gm$-equivariant $A\tw$-bundle
$\id_{_{\BA^1}}\boxtimes \varpi:
\BA^1\times\M\tw\map\BA^1\times\mm\tw.$
We  define
\beq{CF}\rr\CF:=(\id_{_{\BA^1}}\boxtimes \varpi)\idot
\bigl(\rr\D\tw/\rr(\D\tw\cdot I_\chi)\bigr).
\eeq
This 
is a $\Gm$-equivariant quasi-coherent sheaf of  $A$-modules on
$\BA^1\times\mm\tw$ which is
flat over the $\BA^1$-factor. Let $\rf^A$ denote the subsheaf of
its $A$-invariant
sections.

\begin{lemma}\label{commute_gr}
There is a  $\Gm$-equivariant sheaf  isomorphism
$\rf^A|_{\{0\}\times
\mm\tw}\cong\Fr\idot\oo_\mm.$
\end{lemma}

In the course of the proof below,
we will repeatedly use
 the following elementary result
\begin{lemma}\label{proj} Let $\DD$ be a graded algebra,
and $M$ a $t$-torsion free, graded $\DD[t]$-module ($\deg t=1$) such
that $M/tM$ is a rank $m$ free, resp. projective,  graded $\DD$-module.
Then, $M$ is a  rank $m$ free, resp. projective,  $\DD[t]$-module.\qed
\end{lemma}

\begin{proof}[Proof of Lemma \ref{commute_gr}.]
Recall that every point in
any $\Gm$-variety is known to have a
$\Gm$-stable {\em affine} Zariski open neighborhood.
Applying this to the $\Gm$-action
on $\mm$, we may replace $\mm$ by
a $\Gm$-stable affine Zariski open subset
 $\CY\sset \mm$. Let $Y$ be the inverse image
of $\CY$ under the bundle map $\M\to\mm$.
Thus, $Y$ is a $\Gm\times A$-stable  affine Zariski open subset
in $\M$, and we put 
$$\RF:=\G(\BA^1\times\CY\tw,\rf)=\G\bigl(\BA^1\times Y\tw,
\,\rr\D\tw/\rr(\D\tw\cdot I_\chi)\bigr).
$$
Thus, $\RF$ is a graded flat $\k[t]$-module,
where $t$ stands for the coordinate on $\BA^1$.

With these notations, the statement of the Lemma
amounts to the claim that, for all sufficiently small $\Gm$-stable
Zariski open affine subsets $\CY\sset\mm$, there is a natural
graded space  isomorphism
\beq{gr_space}
\RF^A/t\cd\RF^A\cong\G(\CY\tw,\Fr\idot\oo_\CY)\,\Bigl(=\k[\CY]=\k[Y]^A\Bigr).
\eeq

To prove \eqref{gr_space}, let $\ct$ be the Zariski open (possibly empty) subset
in $T^*X$ formed by the points $\xi\in T^*X$
such that the differential
of $\mu: T^*X\to \aa^*$ is surjective at $\xi$.
Further, let $\J\sset \sym\aa=\k[\aa^*]$ denote the augmentation
ideal. It is clear that $\mu\inv(0)$
is the 
 zero scheme of the ideal $\oo_{T^*X}\cdot \mua(\J)
\sset\oo_{T^*X}$, and that this ideal
is reduced at any point of
$\mu\inv(0)\,\cap\,\ct.$ 
It follows,
since   $\gr  I_\chi=\J$, that
on the Frobenius twist of
$\ct$, one has:
$$\gr\bigl(\Fr\idot(\D_X\cd I_\chi)\bigr)\big|_{\ct\tw}=
\bigl((\gr\Fr\idot\D_X)\cd(\gr \Fr\idot I_\chi)\bigr)\big|_{\ct\tw}
=\Fr\idot\bigl(\oo_{T^*X}\cd\mua(\J)\bigr)\big|_{\ct\tw},
$$
where we 
 identify $\gr(\Fr\idot\D_X)$ 
with the corresponding $\Gm$-equivariant sheaf on $T^*X\tw$.

We now use our Basic Assumptions saying that  the differential
of $\mu$ is surjective at any point of the open subset
$\M\sset\mu\inv(0)$. Hence,
$\M\sset\ct\cap \mu\inv(0)$ is a non-empty Zariski-open
subset in $\mu\inv(0)$. 
Since the sheaf $\gr(\Fr\idot\D_X)/\gr(\Fr\idot\D_X\cd I_\chi)$
is  supported on $[\mu\inv(0)]\tw$, we obtain
\begin{align}\label{long1}
\bigl(\gr(\Fr\idot\D_X)/\gr(\Fr\idot\D_X\cd
I_\chi)\bigr)\big|_{\M\tw}
&=\bigl(\gr(\Fr\idot\D_X)/\gr(\Fr\idot\D_X\cd
I_\chi)\bigr)\big|_{\ct\tw\cap\M\tw}\\
&=\bigl((\Fr\idot\oo_{T^*X})\big/
\Fr\idot\bigl(\oo_{T^*X}\cd\mua(\J)\bigr)\bigr)\big|_{\ct\tw\cap\M\tw}\nonumber\\
&=\Fr\idot(\oo_{T^*X}\big|_{\ct\cap\M})
=\Fr\idot\oo_{\M}.\nonumber
\end{align}

The definition of Rees
algebra  implies readily, see \eqref{rees_iso} and \eqref{rees_D},
that the restriction to $\{0\}\times \M\tw$ of the
sheaf $\rr\D\tw/\rr\D\tw\cd\rr{I\tw_\chi}$ is isomorphic to
$\bigl(\gr(\Fr\idot\D_X)/\gr(\Fr\idot\D_X\cd
I_\chi)\bigr)\big|_{\M\tw}$,
the sheaf in the top line of
\eqref{long1}.
 Hence, equations \eqref{long1} and a flat base change
yield 
$\rf\big|_{\{0\}\times\mm\tw}\cong\varpi\idot\Fr\idot\oo_{\M}.$
Restricting this sheaf isomorphism to our 
affine open subset $\CY\tw$ and taking global sections,
we obtain canonical graded space isomorphisms
\beq{kY}
\RF/t\cd\RF \cong
\G\bigl(\CY\tw,\,\rf|_{\{0\}\times\CY\tw}\bigr)
\cong \G(\CY\tw,\,\varpi\idot\Fr\idot\oo_{Y})\cong\G(Y\tw,
\Fr\idot\oo_{Y})=\k[Y].
\eeq

We conclude, comparing the above isomorphisms with those in
\eqref{gr_space} that
proving the Lemma
reduces to  the following
result: {\em The canonical 
map below gives, for $\CY$ sufficiently small, an isomorphism}
\beq{reduction2}
\RF^A/t\cd\RF^A\iso (\RF/t\cd\RF)^A.
\eeq

To prove this, we may assume,
shrinking $\CY$ if necessary,
that there is an \'etale map $\theta:\widetilde{\CY}\to\CY$
such that  the principal $A$-bundle
$Y\to \CY$  becomes trivial
after pull-back via $\theta$, i.e.,
 that there is
an $\Gm\times A$-equivariant isomorphism
 $Y\times_\CY\widetilde{\CY} \cong A\times \widetilde{\CY}$.
Using \'etale
base change for the Cartesian square 
$$
\xymatrix{
A\times \widetilde{\CY}\ar@{=}[r]&Y\times_{\CY}\widetilde{\CY}
\ar[r]^<>(0.5){\theta}\ar[d]&Y\ar[d]\\
&\widetilde{\CY}\ar[r]^<>(0.5){\theta}&\CY}
$$
we obtain  $\Gm\times A$-equivariant graded algebra isomorphisms
\beq{kkk}
\k[\widetilde{Y}]=
\k[\widetilde{\CY}]\otimes_{\k[\CY]}\,\k[Y]\cong
\k[Y\times_\CY\widetilde{\CY}]\cong\k[A]\otimes \k[\widetilde{\CY}].
\eeq
In the leftmost term of this formula,
we have used the notation $\widetilde{Y}:=Y\times_\CY\widetilde{\CY}$,
and in the rightmost term of the formula,
the  group $A$ acts trivially on 
the factor $\k[\widetilde{\CY}]$.

Next, we set
$$
\widetilde{\RF}:=\G\Big(\BA^1\times \widetilde{Y}\tw,
\,(\id_{_{\BA^1}}\times\theta\tw)^*\bigl(\rr\D\tw/\rr(\D\tw\cdot I_\chi)\bigr)|_{\BA^1\times
Y\tw}\Big)=\k[\BA^1\times \widetilde{Y}\tw]\bigotimes_{\k[\BA^1\times \CY\tw]}\,\RF.
$$
Since $\k[\widetilde{Y}\tw]$ is flat
over $\k[Y\tw]$, from \eqref{kY} and \eqref{kkk}
we find 
$$
\widetilde{\RF}/t\cd\widetilde{\RF}=
\k[\widetilde{Y}\tw]\otimes_{\k[\CY\tw]}\,(\RF/t\cd\RF)=
\k[\widetilde{Y}\tw]\otimes_{\k[\CY\tw]}\,\k[Y]=\k[\widetilde{Y}]
\cong\k[A]\otimes \k[\widetilde{\CY}].
$$
Thus, we see that $\widetilde{\RF}/t\cd\widetilde{\RF}$
is a  rank one  free
$\k[A]\otimes \k[\widetilde{\CY}]$-module.
We deduce from Lemma \ref{proj} that
$\widetilde{\RF}$ is  isomorphic
to a rank one free
$\k[A]\otimes \k[\BA^1\times \widetilde{\CY}]$-module.
Furthermore, it is easy to show that
this isomorphism can be chosen to be $A$-equivariant.

Now, the functor $M\mto M^A$, of $A$-invariants,
 takes short exact sequences of $A$-modules
of the form $\k[A]\otimes E$ to short exact sequences.
Hence applying this functor
to the  short exact sequence
$$0\too\widetilde{\RF}\stackrel{t}\too\widetilde{\RF}\too\widetilde{\RF}/
t\cd\widetilde{\RF}
\too 0,$$
we deduce the isomorphism
$$
\widetilde{\RF}^A/t\cd\widetilde{\RF}^A\iso (\widetilde{\RF}/t\cd\widetilde{\RF})^A.
$$
The latter  isomorphism yields \eqref{reduction2} since
the morphism $\widetilde{Y}\tw\to Y\tw$ is faithfully flat.
The Lemma is proved.
\end{proof}

\subsection{Deformation of the algebra $\az_\chi$.}\label{mimick}
 We mimick formulas
\eqref{mtw} and \eqref{Dchi} and put
\begin{align}\label{rchi}
\rr\D\tw_\chi\,&:=\,\rr\D\tw/\rr(\D\tw\cdot I\tw_\chi),\quad\text{and}\\
&\rr\ee_\chi\,:=\,\bigl(\rr\D\tw/\rr(\D\tw\cdot I_\chi)\bigr)^{A_1}
\cong\bigl({\rr\D\tw_\chi}\big/\rr(\D\tw_\chi
\cd\fri_\chi)\bigr)^{A_1}.\nonumber
\end{align}
The sheaves
$\rr\D\tw_\chi$ and $\rr\ee_\chi$ are both supported on
 $\BA^1\times
[\mu\inv(0)]\tw$, see \eqref{setI}. 
Following the same strategy as has been used in the
construction of the sheaf $\az_\chi$,
we define
$$ \rr\az_\chi:=
\bvp\idot\bigl(\rr\ee_\chi\big|_{\BA^1\times 
\M\tw}\bigr)^{A\tw}.
$$
This is a $\Gm$-equivariant sheaf of associative 
algebras on $\BA^1\times \mm\tw$, 
viewed as a $\Gm$-variety with diagonal action.
Write $\op{pr}_2: \BA^1\times \mm\tw\to\mm\tw$ for the second projection.

\begin{lemma}\label{rees_az} The sheaf  $\rr\az_\chi$ is  flat over the
$\BA^1$-factor, and we have
$$
\rr\az_\chi\big|_{\{0\}\times\mm\tw}\,\cong\,\Fr\idot\oo_\mm,\quad
\text{and}\quad\rr\az_\chi\big|_{(\BA^1\sminus\{0\})\times \mm\tw}
\,\cong\,\op{pr}_2^*\az_\chi.
$$
\end{lemma}

\begin{proof}
 Both the flatness statement and the isomorphism on the right
are immediate from the
corresponding properties of the Rees algebra.

It remains to study the restriction of the sheaf
$\rr\az_\chi$ to the special divisor
$\{0\}\times\mm\tw$. From the definition, we find
\begin{align*}
\rr\az_\chi\big|_{\{0\}\times \mm\tw}\,&\cong\,
\bigl((\id_{_{\BA^1}}\boxtimes \varpi)\idot\rr\ee_\chi\bigr)^{A\tw}
\big|_{\{0\}\times \mm\tw}\\
&\cong\,
(\id_{_{\BA^1}}\boxtimes \varpi)\idot
\bigl(\rr\D\tw/\rr(\D\tw\cd I_\chi)^{A_1}\bigr)^{A\tw}\big|_{\{0\}\times
\mm\tw}\\
&\cong\,(\id_{_{\BA^1}}\boxtimes \varpi)\idot
\bigl(\rr\D\tw/\rr(\D\tw\cd I_\chi)\bigr)^{A}\big|_{\{0\}\times
\mm\tw}=\rf^{A}\big|_{\{0\}\times
\mm\tw}.
\end{align*}
But $\rf^{A}\big|_{\{0\}\times
\mm\tw}\cong \Fr\idot\oo_{\mm}$ by Lemma \ref{commute_gr},
and we are done.
\end{proof}

Next, we set $\BR:= \G(\BA^1\times \mm\tw,\,\rr\az_\chi)$,
 a graded flat $\k[t]$-algebra such that
$\BR/(t-1)\BR$
$=\G(\mm\tw,\az_\chi)$.
Applying  formula
\eqref{filt} to the algebra  $\rr:=\BR$
 and using  Lemma \ref{rees_az}, we get
a natural increasing filtration
on
the algebra $\G(\mm\tw,\az_\chi)$
such that for the associated graded algebra,
to be denoted
$\gr^\rr\G(\mm\tw,\az_\chi),$
we have $\gr^\rr\G(\mm\tw,\az_\chi)=\BR/t\BR.$
On the other hand,
the $\Gm$-action  induces a grading on
the algebra $\G(\mm, \oo_\mm)=\G(\mm\tw,\Fr\idot\oo_\mm)$.

Recall the affinization morphism $f: \mm\to\mm_{\op{aff}}$.
Given $\xi\in \mm\tw_{\op{aff}}$,
let $\mm\tw_\xi:=[f\tw]\inv(\xi)
\sset\mm\tw$ be the fiber of $f\tw: \mm\tw\to\mm\tw_{\op{aff}}$
over $\xi$.

\begin{prop}\label{gr_az} \vi The sheaf $\az_\chi$ is
locally free, and $H^i(\mm\tw,\az_\chi)=0$ for all $i>0$.

\vii There is a graded algebra isomorphism
$$
\gr^\rr\G(\mm\tw,\az_\chi)\cong 
\G(\mm,\oo_\mm).
$$

\viii For any $q>0$ we have $R^qf\tw_*\az_\chi=0$ and
$H^q(\mm\tw_\xi,\az_\chi)=0,\,\forall \xi\in \mm\tw_{\op{aff}}.$
\end{prop}

\begin{proof}
 We consider the following diagram
\beq{cart}
\xymatrix{
\mm\ar[r]^<>(0.5){\Fr}\ar[d]^<>(0.5){f}&{\mm\tw_{}}={\{0\}\times\mm\tw_{}}\ar[d]^<>(0.5){f\tw}
\ar@{^{(}->}[rr]^<>(0.5){\ti}&&\BA^1\times \mm\tw
\ar[d]^<>(0.5){\tilde{f}:=\id_{\BA^1}\times f\tw}\\
\mm_{\op{aff}}\ar[r]^<>(0.5){\Fr}&\mm\tw_{\op{aff}}=\{0\}\times\mm\tw_{\op{aff}}
\ar@{^{(}->}[rr]^<>(0.5){i}&&\BA^1\times \mm\tw_{\op{aff}}.
}
\eeq

Lemma \eqref{rees_az} says that
$\ti^*\rr\az_\chi\cong \Fr\idot\oo_\mm$.
Thus, $\rr\az_\chi$ is a $\Gm$-equivariant sheaf on $\BA^1\times\mm\tw$
such that its restriction to the subvariety
$\{0\}\times\mm\tw$ is a locally free sheaf.
It follows  that the sheaf $\rr\az_\chi$ must be itself
locally free. Indeed, every   point in  $\mm\tw$
has  a
$\Gm$-stable {\em affine} Zariski open neighborhood.
Taking global sections
of  $\rr\az_\chi$  over such a neighborhood, we see that our claim
reduces to Lemma \ref{proj}.
Thus, we have proved that the sheaf  $\rr\az_\chi$,
 hence its restriction
to $\{1\}\times\mm\tw$, is a locally free sheaf.
But $\rr\az_\chi\big|_{\{1\}\times\mm\tw}=\az_\chi$,
thus, the first  claim of part (i) of the Lemma is proved.

Recall further that the variety $\mm$ is symplectic, hence
it has  trivial canonical bundle.
Therefore, by the Grauert-Riemenschneider theorem
(see \cite{EV} for $\charp\k >0$ case),
 the higher direct image sheaves $R^qf_*\oo_\mm$ vanish for all
$q>0$. Therefore, from the commutative square on the left
of \eqref{cart}
we deduce 
\beq{grauert}
R^qf\tw_*(\ti^*\rr\az_\chi)=R^qf\tw_*\Fr\idot\oo_\mm=0,
\quad\text{for all}\quad q>0.
\eeq

We are going to use \eqref{grauert} to prove part (iii) of the Proposition.
To this end, let $\wh{D}:=\Spec\k[[t]]\sset \BA^1$
denote the completion of the line $\BA^1$ at the origin.
Thus,  $\wh{D}\times\mm\tw$, resp., 
 $\wh{D}\times\mm\tw_{\op{aff}}$, 
is  the
formal completion of the scheme $\BA^1\times\mm\tw$
 along the  closed subscheme $\{0\}\times \mm\tw$,
 resp., 
formal completion of
$\BA^1\times \mm\tw_{\op{aff}}$
  along the subscheme $\{0\}\times \mm\tw_{\op{aff}}$.
We have  $\wh{D}\times\mm\tw=\tilde{f}\inv(\wh{D}\times\mm\tw_{\op{aff}})$,
and the restriction of $\tilde{f}$ gives
a morphism
$\wh{f}:\ \wh{D}\times\mm\tw\map \wh{D}\times\mm\tw_{\op{aff}}$.

Let $\rr\az_\chi\big|_{\wh{D}\times\mm\tw}$ denote the
restriction of $\rr\az_\chi$ to $\wh{D}\times\mm\tw$.
We will view  $\rr\az_\chi\big|_{\wh{D}\times\mm\tw}$ as a sheaf
of abelian groups on the closed fiber $\mm\tw=\{0\}\times \mm\tw$.
This sheaf is clearly isomorphic to an inverse limit
of sheaves which are iterated extensions of the
sheaf $\ti^*\rr\az_\chi=\Fr\idot\oo_\mm$.
Hence, formula \eqref{grauert} implies that
$$
R^q\wh{f}_*(\rr\az_\chi\big|_{\wh{D}\times\mm\tw})
=0,
\quad\text{for all}\quad q>0.
$$

On the other hand, let $\wh{i}: \wh{D}\times\mm\tw_{\op{aff}}\into
\BA^1\times \mm\tw_{\op{aff}}$ denote the imbedding.
For any $q=0,1,\ldots,$
one has $\wh{i}^*R^q\tilde{f}_*\rr\az_\chi=
R^q\wh{f}_*(\rr\az_\chi\big|_{\wh{D}\times\mm\tw}),$
by the Formal Functions
 Theorem, cf. \cite[III, 11.1]{Har} or
[EGA III, Sect. 4].
Thus, we have proved 
\beq{grauert2}
\wh{i}^*R^q\tilde{f}_*\rr\az_\chi=
R^q\wh{f}_*(\rr\az_\chi\big|_{\wh{D}\times\mm\tw})
=0,
\quad\text{for all}\quad q>0.
\eeq

Observe next that, for each $q$, the sheaf $R^q\tilde{f}_*\rr\az_\chi$
is a $\Gm$-equivariant coherent
sheaf on $\BA^1\times \mm\tw_{\op{aff}}$.
Therefore,
the  support of $R^q\tilde{f}_*\rr\az_\chi$
is a  $\Gm$-stable closed subscheme in $\BA^1\times \mm\tw_{\op{aff}}$.
The assumption that  the $\Gm$-action on $\mm$ be  attracting
implies that any non-empty  $\Gm$-stable closed
subscheme in $\BA^1\times \mm\tw_{\op{aff}}$ has a non-empty
intersection with  $\{0\}\times \mm\tw_{\op{aff}}$.
Hence, if the sheaf $R^q\tilde{f}_*\rr\az_\chi$ is non-zero,
it must have a non-zero restriction
to the subscheme  $\{0\}\times \mm\tw_{\op{aff}}$.
Since $\{0\}\times \mm\tw_{\op{aff}}\sset
\wh{D}\times\mm\tw_{\op{aff}}$, this would yield,
in particular, that $\wh{i}^*R^q\tilde{f}_*\rr\az_\chi\neq 0,$
contradicting \eqref{grauert2}.
Part (iii) of the Proposition follows.

Taking  global sections over the
{\em affine} open set $(\BA^1\sminus\{0\})\times\mm\tw_{\op{aff}}$,
from the vanishing result of
part (iii) and  the second isomorphism of Lemma \ref{rees_az},
for any $q>0$, we find
\begin{align*}0=\G\bigl((\BA^1\sminus\{0\})&\times\mm\tw_{\op{aff}},\,
R^q\tilde{f}_*\rr\az_\chi\bigr)=
H^q\bigl((\BA^1\sminus\{0\})\times\mm\tw,\,\rr\az_\chi\bigr)\\
&=H^q\bigl((\BA^1\sminus\{0\})\times\mm\tw,\,\pr_2^*\az_\chi\bigr)
=
\k[t,t\inv]\otimes H^q(\mm\tw,\az_\chi),
\end{align*}
This completes the proof  of part (i) 
of the  Proposition.

Now, the (ordinary) direct image sheaf $\tilde{f}_*\rr\az_\chi$
is by construction flat over $\BA^1$,
hence, $L^qi^*(\tilde{f}_*\rr\az_\chi)=0$
for all $q>0$. Therefore, the
vanishing of the higher direct images implies that
the Proper Base Change theorem for the Cartesian square
on the right of diagram \eqref{cart} involves no higher
derived functors. Thus,  using  Base Change and
the definition of affinization we  obtain
\beq{base_change}
i^*\tf_*\rr\az_\chi=
f\tw_*\ti^*\rr\az_\chi=f\tw_*\Fr_*\oo_\mm=\Fr\idot
f_*\oo_\mm=\Fr\idot\oo_{\mm_{\op{aff}}}.
\eeq
Observe that we have
 $\G(\BA^1\times \mm\tw_{\op{aff}},\,\tilde{f}_*\rr\az_\chi)=
\G(\BA^1\times \mm\tw, \rr\az_\chi)=\BR$, the graded
 $\k[t]$-algebra involved in the definition of filtration on
$\G(\mm\tw,\,\az_\chi)$.
Therefore, applying the  global sections functor
${\G(\BA^1\times \mm\tw_{\op{aff}}, -)}$ to both sides in  \eqref{base_change}
we obtain
$$
\gr^\rr\G(\mm\tw,\az_\chi)=\BR/t\BR\cong
\G(\mm\tw_{\op{aff}},\,\Fr\idot\oo_{\mm_{\op{aff}}})=
\G(\mm\tw,\,\Fr\idot\oo_\mm)=
\G(\mm,\oo_\mm).$$ 
This proves part (ii) of the  Proposition.
\end{proof}

\subsection{Comparison of characters.}\label{comp_char_sec}
Fix two characters  $\chi,\psi$  such that $\chi-\psi\in\dlog\BX(A)$.
Let $\oo_\mm(\psi-\chi)$
be the $(\psi-\chi)$-weight subsheaf
in the push-forward of $\oo_\M$ under the bundle map
$\M\to\mm$. The $(\psi-\chi)$-isotypic component of the
regular representation $\k[A]$ being 1-dimensional,
we conclude that  $\oo_\mm(\psi-\chi)$ is a  rank 1 locally free  sheaf on $\mm$.

Let $\D\tw_\M$ denote the restriction of the Azumaya algebra
$\D\tw$ to the 
subset  $\M^{(1)}\sset \ttw=\ttwp$. We put
${_\chi\ee_\psi}:=\on{Hom}_{{\D\tw_\M}}(\D\tw_\M/\D\tw_\M\cdot\fri_\chi,
\D\tw_\M/\D\tw_\M\cdot\fri_\psi)$.
This is an  $A\tw$-equivariant
 sheaf  on $\M\tw$. Therefore,
$\varpi\idot({_\chi\ee_\psi})$ is
 a sheaf
on $\mm\tw$ with fiberwise $A\tw$-action.
Let ${_\chi\az_\psi}$ be the $(\psi-\chi)$-weight
component of $\varpi\idot({_\chi\ee_\psi})$.
It is clear that  ${_\chi\az_\psi}$ is a coherent
sheaf 
 of $\az_\chi\dash\az_\psi$-bimodules.

Next, we
 mimic the argument in Sect. \ref{deform} and
 observe that the standard increasing
filtration on $\varpi\idot\D\tw_\M$
induces a natural  increasing
filtration on $\varpi\idot({_\chi\ee_\psi})$. Therefore,
using the Rees algebra construction,
we may form a $\Gm$-equivariant sheaf $\rr({_\chi\ee_\psi})$
on $\BA^1\times\M\tw$. This gives an increasing
filtration on
 $\G(\mm\tw,\,{_\chi\az_\psi})$, an ${\mathsf{A}}_\chi\dash{\mathsf{A}}_\psi$-bimodule,
and we write $\gr^\rr\G(\mm\tw,\,{_\chi\az_\psi})$
for the associated graded $(\gr{\mathsf{A}}_\chi)\dash(\gr{\mathsf{A}}_\psi)$-bimodule.

\begin{lemma}\label{morita_az}
\vi The sheaf ${_\chi\az_\psi}$ is
locally free, moreover, we have 
a natural $\ee_\chi\dash\ee_\psi$-bimodule
isomorphism ${_\chi\ee_\psi}=\varpi^*({_\chi\az_\psi})$.

\vii Assume that $R^qf_*\oo_\mm(\chi-\psi)=0$ holds for all $q>0$.
Then, there is a graded
$(\gr{\mathsf{A}}_\chi)\dash(\gr{\mathsf{A}}_\psi)$-bimodule 
isomorphism
$\gr^\rr\G(\mm\tw,\,{_\chi\az_\psi})\cong 
\G(\mm,\oo_\mm(\chi-\psi)).
$ Moreover, we have
$$R^qf\tw_*({_\chi\az_\psi})=0\quad
\text{and}\quad
H^q(\mm\tw,\,{_\chi\az_\psi})=0,\quad\forall 
q>0.
$$
\end{lemma}

\begin{proof}  We argue as in 
 Proposition \ref{gr_az},
and show first that
\beq{psiphi_def}
\bigl((\id_{\BA^1}\times\varpi)\idot\rr({_\chi\ee_\psi})\bigr)^{\psi-\chi}
\Big|_{\{0\}\times\mm\tw}\,\cong\,\Fr\idot\oo_\mm(\psi-\chi).
\eeq
Now, the sheaf 
$\Fr\idot\oo_\mm(\psi-\chi)$ is clearly locally free.
Hence, arguing as in  \eqref{Efree} we deduce 
that the sheaf $\rr({_\chi\ee_\psi})\big|_{\{0\}\times\M\tw}$
is a locally free sheaf
on $\M\tw$ which is isomorphic to $\varpi^*\Fr\idot\oo_\mm(\psi-\chi)$.
Using Lemma \ref{proj} one shows that $\rr({_\chi\ee_\psi})$
is a  locally free sheaf
on $\BA^1\times\M\tw$. Thus,
${_\chi\ee_\psi}=\rr({_\chi\ee_\psi})\big|_{\{1\}\times\M\tw}$
is a  locally free sheaf
on $\M\tw$ and, moreover, we have 
$
{_\chi\ee_\psi}=\varpi^*({_\chi\az_\psi}).$
\end{proof}

\begin{rem}\label{T^c} The above argument shows that there is a
flat family (over $\BA^1$) of  coherent sheaves such that
nonzero members of the family are all isomorphic to 
${_\chi\az_\psi}$ 
and the fiber over $0\in\BA^1$ is isomorphic
to $\Fr\idot\oo_\mm(\psi-\chi)$. In particular,
in the Grothendieck group of coherent sheaves
on $\mm\tw$ one has an equality
 $[{_\chi\az_\psi}]=[\Fr\idot\oo_\mm(\psi-\chi)]$.
\end{rem}
\subsection{Proof of Theorem \ref{roma}.}\label{proof_roma}
We have $\az_\chi=
\varpi\idot\bigl(\ee_\chi\big|_{\M\tw}\bigr)^{A\tw}$.
The sheaf $\az_\chi$ in the left hand side of this equality is 
locally free, by Proposition~\ref{gr_az}(i).
It follows that the sheaf $\ee_\chi\big|_{\M\tw}$
is also locally free, moreover, 
we have $\ee_\chi\big|_{\M\tw}\cong\varpi^*\az_\chi,$ see
\eqref{Efree}.

Now, let  $x\in\M^{(1)}$ and write
 $\bar{x}:=\varpi(x)\in\mm\tw$ for its image.
Let $\ee_x$, resp. $\az_{\bar{x}}$, denote  the
 geometric  fiber at $x$, resp. at
$\bar{x}$, of the corresponding  locally free
sheaf.
We deduce from  $\ee_\chi\big|_{\M\tw}\cong\varpi^*\az_\chi$ that
there is an algebra isomorphism
\beq{stalks}
\ee_x\cong \az_{\bar{x}},\quad\text{for all}\quad x\in\M^{(1)}\quad
\text{and}\quad
\bar{x}:=\varpi(x)\in\mm\tw.
\eeq

Thus, to prove that  $\az_\chi$ is an Azumaya algebra,
it suffices to show that $\ee_x$ is a matrix algebra,
for any $x\in\M^{(1)}.$
 By definition,
we have $\ee_x= (\D\tw_x/\D\tw_x\cdot\fri_\chi)^{A_1},$
where $\D\tw_x$ is the geometric  fiber at $x$
of the sheaf $\D\tw$. 

We know that $\D\tw$
is an Azumaya algebra on $T^*X\tw$.
Hence, there is a vector space $E$ and an 
algebra isomorphism $\D\tw_x\cong \End_\k E$.
Using this, from the last statement of Lemma \ref{weight}
we deduce the following algebra isomorphisms
$$\ee_x= (\D\tw_x/\D\tw_x\cdot\fri_\chi)^{A_1}\cong
\bigl(\End E\big/\End E\cd\fri_\chi\bigr)^{A_1}=
\DD_\chi\cong \End_\k(E_\chi).
$$
Thus,  $\ee_x$ is a matrix algebra, as claimed.
This proves that $\az_\chi$ is an Azumaya algebra.

To complete the proof of part (i) of the Theorem,
consider the following chain
of canonical algebra maps
\begin{align}\label{gl_az}
\G(X,\,\D_X/\D_X\cdot&I_\chi)^A\;
\underset{^\sim}{\stackrel{\text{Lemma \ref{stages}}}\tooo}\;
\G\bigl(\ttw,\,\ee_\chi\bigr)^{A\tw}
\stackrel{\text{restriction}}{-\!\!\!-\!\!\!-\!\!\!-\!\!\!\longrightarrow}\;
\G\bigl(\M\tw,\,\ee_\chi\big|_{\M\tw}\bigr)^{A\tw}\\
&\cong
\G\bigl(\mm\tw,
\,\varpi\idot(\ee_\chi\big|_{\M\tw})\bigr)^{A\tw}\cong
\G\bigl(\mm\tw,
\,\varpi\idot(\ee_\chi\big|_{\M\tw})^{A\tw}\bigr)\cong
\G(\mm\tw,\az_\chi).\nonumber
\end{align}
The composite map provides the algebra map claimed in  part (i) of Theorem
\ref{roma}.

Part (ii) of the Theorem follows directly from the cohomology
vanishing in Proposition~\ref{gr_az}(i); 
To prove   (iv), we 
 fix a point
$x\in\M\tw$, and write $\D\tw_x\cong \End_\k E$.
We know
that for the geometric fibers (at $x$) of the Azumaya
algebras $\ee_\chi,$ resp., $\ee_\psi$,
one has the following formulas
$$(\ee_\chi)_x=\End_\k E_\chi,\quad\text{resp.}\quad
(\ee_\psi)_x=\End_\k E_\psi.$$
 Further, the sheaf ${_\chi\ee_\psi}$
is locally free by Lemma \ref{morita_az}(i).
The fiber of  that sheaf at $x$ is
an $(\ee_\chi)_x\dash(\ee_\psi)_x$-bimodule,
and Lemma \ref{weight}(ii) yields 
 the following $(\ee_\chi)_x\dash(\ee_\psi)_x$-bimodule
isomorphisms, 
cf. \eqref{bimod}: 
\begin{align*}
({_\chi\ee_\psi})_x&\cong
\on{Hom}_{{\D\tw_x}}\bigl(\D\tw_x/\D\tw_x\cd\fri_\chi,\,
\D\tw_x/\D\tw_x\cd\fri_\psi\bigr)^{\psi-\chi}\\
&\cong
\Hom_{\End E}\bigl(\End E\big/\End E\cd\fri_\chi,\,
\End E\big/\End E\cd\fri_\psi\bigl)^{\psi-\chi}
={_\psi\DD_\chi}\cong\Hom_\k(E_\chi, E_\psi).
\end{align*}
Thus,  the sheaf ${_\chi\ee_\psi}$  is a 
sheaf 
 of locally-projective  $\ee_\chi\dash\ee_\psi$-bimodules.
In particular, the Azumaya 
 algebras  $\ee_\chi$ and $\ee_\psi$
are  Morita equivalent.

Now, for   $\bar{x}:=\varpi(x)\in\mm\tw$, we have
$({_\chi\az_\psi})_{_{\bar{x}}}=({_\chi\ee_\psi})_x$.
Therefore,
we see from \eqref{stalks} and Lemma \ref{morita_az}(i)
that the sheaf ${_\chi\az_\psi}$  is a 
sheaf 
 of locally-projective $\az_\chi\dash\az_\psi$-bimodules.
Hence, it 
provides the required  Morita equivalence
between the   Azumaya algebras  $\az_\chi$ and $\az_\psi$.

Part (iii) of Theorem \ref{roma}
is a special case of the following  more general result
due to \cite[Proposition 2.2]{bk}.

\begin{prop}\label{bk_equiv}
Let $\mm$ be a smooth connected variety over $\k$
with the trivial canonical class, and such that
the morphism $\mm\to\mm_{\op{aff}}$ is proper.
Let
$\az$ be an Azumaya algebra on $\mm$ such that $H^i(\mm,\az)=0,
\,\forall i>0$ and, moreover,
the   algebra $\A:=\G(\mm, \az)$  has finite homological
dimension.

Then, the  functor below provides an equivalence
beteen the bounded derived categories
of sheaves of coherent $\az$-modules and finitely-generated
$\A$-modules, respectively:
$$D^b(\az\dash\op{Mod})\too
D^b(\A\dash\on{Mod}^{\text{fin. gen.}}),\quad
{\scr F}\mto \op{RHom}_{\az\dash\op{Mod}}(\az,{\scr F}).
$$
\end{prop}

The proof of this Proposition exploits the technique
of Serre functors, and is similar in spirit to the
proof of \cite[Theorem 2.4]{BKR}.

Part (iii) of our Theorem  follows from the Proposition since
for any $\az$-module ${\scr F}$, one has
$\op{RHom}_{\az\dash\op{Mod}}(\az,{\scr F})\cong
\op{RHom}_{\oo_{_\mm}\dash\op{Mod}}(\oo_{_\mm}, {\scr F})=
\op{\RG}(\mm,{\scr F}).$

This completes the proof of the Theorem \ref{roma}.\qed

\section{The rational Cherednik algebra of type $\mathbf{A}_{n-1}$.}\label{hh}
\subsection{Basic definitions.}
Let $W:=S_n$ denote the Symmetric group and 
 $\Z[W]$ denote the group algebra of $W$.
Write $s_{ij}\in W$ for the transposition $i\leftrightarrow j$.
We consider two sets of variables
$x_1,\ldots,x_n,$ and $y_1,\ldots,y_n,$ and let
$W=S_n$ act on the polynomial algebras
$\Z[x_1,\ldots,x_n]$ and
$\Z[y_1,\ldots,y_n]$ by permutation
of the variables.

Let $\mathbf{c}$ be an indeterminate.
We define  the rational Cherednik algebra
of type $\mathbf{A}_{n-1}$ as a $\Z[\mathbf{c}]$-algebra, $\hh$,
with generators
$x_1,...,x_n,y_1,...,y_n$ and  $\Z[W]$, and
the following defining relations, see~\cite{EG}:
\begin{equation}\label{a_n}
\begin{array}{lll}\displaystyle
&{}_{_{\vphantom{x}}}s_{ij}\cdot x_i=x_j\cdot s_{ij}\quad,\quad
 s_{ij}\cdot y_i=y_j\cdot s_{ij}\,,&
\forall i,j\in\{1,2,\ldots,n\}\;,\;i\neq j\break\medskip\\
&{}_{_{\vphantom{x}}}{}^{^{\vphantom{x}}}
[y_i,x_j]= {\mathbf{c}}\cdot s_{ij}\quad,\quad[x_i,x_j]=0=[y_i,y_j]
\,,&
\forall i,j\in\{1,2,\ldots,n\}\;,\;i\neq j\break\medskip\\
&{}^{^{\vphantom{x}}} [y_k,x_k]= 1-{\mathbf{c}}\cdot\sum_{i\neq k}\;s_{ik}\,.
\end{array}
\end{equation}

Given a field $\k$ and  $c\in \k$, we let
$\hh_c:= \k\bigotimes_{\Z[\mathbf{c}]}\,\hh$ be the $\k$-algebra 
obtained from
$\hh$ by extension of scalars via the homomorphism
$\Z[\mathbf{c}]\to\k,\,f\mto f(c).$

We keep our standing assumption $\charp \k > n$, and
write  $\e=\frac{1}{n!}\sum_{g\in
W}\,g\in \k[W]\sset\hh_c\,$
 for the symmetrizer idempotent.
  Let $\ehe\sset \hh_c$
be the
 {\em Spherical
subalgebra}, see~[EG].

Let $\h:=\k^n$ be the tautological permutation representation
of $W$. We  identify the variables
$x_1,\ldots,x_n,$ resp. $y_1,\ldots,y_n,$
with coordinates on $\h$, resp. on $\h^*$.
The algebras $\hh_c$ and $\ehe$ come equipped
with compatible increasing filtrations such that
all elements of $W$ and $x_i\in \h^*\sset\hh_c$ have
filtration degree zero, and elements 
$y_i\in\h\sset\hh_c$ have
filtration degree $1$. 

The {\em Poincar\'e-Birkhoff-Witt theorem}
for  rational Cherednik algebras, cf. \cite{EG}, yields graded algebra isomorphisms
\beq{pbw}
\gr\hh_c\,\cong\,\k[\h^*\times\h]\#W,
\quad\text{and}\quad
\gr(\ehe)\,\cong\,\k[\h^*\times\h]^W.
\eeq

\subsection{Dunkl representation.}\label{B}
Write $\hreg$ for
an affine Zariski open dense subset of $\h$
formed by points with pairwise distinct coordinates.
The group  $W$ acts  naturally on the algebra
$\D(\hreg)$ of crystalline differential operators on 
$\hreg$
and we let $\D(\hreg)^W\sset \D(\hreg)$ be the subalgebra
of $W$-invariant differential operators.
The standard increasing filtration on the algebra of differential
operators induces an increasing filtration on the
subalgebra $\D(\hreg)^W$, and we have $\gr\D(\hreg)^W\cong
\k[T^*\hreg]^W=\k[\h^*\times\hreg]^W.$

According to Cherednik, see also \cite{EG},\cite{DO},
there is an injective algebra homomorphism
\beq{dunkl}
\Th_c:\ \ehe\into \D(\hreg)^W,
\eeq
called {\em Dunkl representation} of the algebra $\ehe$.

Let $\B_c:=\Th_c(\ehe)\sset \D(\hreg)^W$ be the image of the
Dunkl representation.
 We equip the algebra
$\B_c$ with increasing filtration induced from the
standard one on $\D(\hreg)^W$. 
Then, $\gr\B_c$ becomes a graded subalgebra in
$\gr\D(\hreg)^W=\k[\h^*\times\hreg]^W,$
which is known to be equal to
$\k[\h^*\times\h]^W\sset\k[\h^*\times\hreg]^W,$
see \cite{EG}.
Further, the map $\Th_c$ in \eqref{dunkl}
is known to be filtration preserving,
and it was proved in \cite{EG} that
the associated graded map, $\gr\Th_c$, induces graded
algebra isomorphisms
\beq{gr_dunkl}
{\xymatrix{
\k[\h^*\times\h]^W
\ar@/^2pc/[rrrr]|-{\;\Id\;}
\ar@{=}[r]
&\gr(\ehe)
\ar[rr]^<>(0.5){\gr\Th_c}_<>(0.5){\sim}&&
{\gr\B_c}\ar@{=}[r]&
\k[\h^*\times\h]^W}}
\sset \k[\h^*\times\hreg]^W.
\eeq

\subsection{The `radial part' construction.}\label{hc_sec}
In this section, we let  $\k$ be an arbitrary 
algebraically closed field,
either of characteristic zero or of characteristic $p$.

Let $V$ be an $n$-dimensional vector space over $\k$.
In case the field $\k$ has finite characteristic
we assume throughout that $\charp\k > n\geq2$.

We put   $G=\GL(V)$ and let $\g=\Lie G=\gl(V)$ be the Lie algebra of
$G$. 
We consider the vector space $\EE:=\g\times V$.

\begin{defn}\label{eer} Let $\eer\sset\EE=\g\times V$ be
a Zariski open dense subset formed by
the pairs $(x,v)$ such that $v$ is a 
{\em cyclic vector} for the operator $x: V\to V$.
\end{defn}

We recall that the endomorphism $x\in\g$ admits a cyclic vector
if and only if $x$ is {\em regular} (not necessarily
semisimple), i.e., 
the centralizer of $x$ in $\g$ has dimension $n$.

Fix a nonzero volume element ${\mathsf{vol}}\in\wedge^n V^*$. We
introduce the following polynomial 
function on~$\EE$:
\beq{ffs}
(x,v)\mto \ffs(x,v):=\langle\vol, v\wedge x(v)\wedge\ldots
\wedge x^{n-1}(v)\rangle.
\eeq
It is clear that we have $\eer=\EE\sminus\ffs\inv(0)$,
in particular, $\eer$ is an affine
variety.

The group
$G$ acts  on  $\g$ via the adjoint action, and acts
naturally on $V$. This gives a $G$-diagonal action on $\EE$
such that $\eer$ is
a $G$-stable subset of $\EE$.
We compose  the first projection
$\EE=\g\times V\to\g$ with the adjoint
quotient map $\g\to\g/\Ad G=\h/W$,
and restrict the resulting morphism
to the subset $\eer\sset\EE$. This way we get
a morphism $p: \eer\to \h/W$.
The group $G$
clearly acts along the fibers of $p$, and we have
the following  well-known result.

\begin{lemma}\label{U}
\vi The $G$-action on $\eer$ is free and
each fiber of $p$ is a single $G$-orbit; 

\vii Furthermore, the map
$p: \eer\to \h/W$ is a universal geometric
quotient morphism. \qed
\end{lemma}

It follows from the Lemma that $\eer$ is a
principal $G$-bundle over $\h/W$.

Given an integer $c\in\Z$, we put
$$\oo(\eer,c):=\{f\in \k[\eer]\mid g^*(f)=(\det g)^c\cdot f,
\enspace\forall g\in G\}.
$$
It is clear that pull-back via the bundle projection $p$
makes $\oo(\eer,c)$ a $\k[\h/W]$-module.
Also,
observe that $\ffs\in \oo(\eer,1)$.

\begin{corollary}\label{U_cor} For any $c\in \Z$, the
space $\oo(\eer,c)$ is a rank one free  $\k[\h/W]$-module
with generator $\ffs^c$.\qed
\end{corollary}

\begin{notation}\label{IJc} For $c\in\k$, we consider a Lie algebra
homomorphism $\chi_c: \g\to\k,\, x\mto c\cdot\Tr(x)$.
Let
$I_c:=I_{\chi_c}\sset\Ug,$ 
denote
 the
two-sided ideal generated by the elements
$\{x-\chi_c(x)\}_{x\in\g},$ cf. Definition~\ref{IJ}.
\end{notation}

The action of $G$ on $\eer$ induces
an algebra map $\Ug\to \D(\eer)$.
We fix $c\in\k$, and
perform the Hamiltonian reduction of the sheaf $\D_\eer$,
of crystalline differential operators on $\eer$, at the
point $\chi_c$. This way, we get an associative algebra
$[\D(\eer)/\D(\eer)\cd I_c]^G$.

 From Lemma \ref{U} and the isomorphism on the right of formula
\eqref{XA} we deduce

\begin{prop}\label{U_prop} For any $c\in\k$, there is a natural 
algebra isomorphism
$$[\D(\eer)/\D(\eer)\cd I_c]^G\cong\D(\h/W).$$
\end{prop}

More explicitly, if $c$ is an integer, then the isomorphism of
the Proposition is obtained by transporting
the action of differential operators on $\eer$
via the
bijection $\k[\h/W]\iso \oo(\eer,c),\, f\mto \ffs^c\cdot p^*(f)$,
provided by Corollary \ref{U_cor}. 
\begin{rem} The explicit construction
of the isomorphism
shows in particular
that
the algebra $\D(\h/W,\chi_c)$ of  {\em twisted}  differential operators
coming from  the right hand side of the general
formula \eqref{XA} turns out to be canonically isomorphic, in our 
case, to 
the algebra $\D(\h/W)$ of {\em ordinary} differential operators.
Thus,  we have put
 $\D(\h/W)$ on the right hand side of the
isomorphism of Proposition \ref{U_prop} (although we have only justified
this for integral values of $c$, the same holds for arbitrary
 values of $c$ as well).
\erem

Observe further that the algebras
$[\D(\eer)/\D(\eer)\cd I_c]^G$ and $\D(\h/W)$
both come equipped with natural
increasing filtrations and the isomorphism
 of
the Proposition is filtration preserving.

The isomorphism of Proposition \ref{U_prop}
may be viewed as a refined version of the
`radial part' construction considered in \cite{EG}.

\begin{rem} The action of differential operators on
$\k[\hreg]^W$ gives rise to the following natural 
algebra inclusions: $\D(\h)^W\sset\D(\h/W)\sset\D(\hreg)^W$.
We also remark  that the space $\D(\h)^W$ has {\em infinite} codimension
in $\D(\h/W)$.
\end{rem}
\subsection{A Harish-Chandra homomorphism.}\label{hc2_sec}
Let $\grs\sset\g$
denote the Zariski open dense subset
of semisimple regular elements.
Observe that the eigen-spaces of an
element $x\in\grs$ give a direct
sum decomposition $V=\ell_1\bigoplus\ldots\bigoplus
\ell_n$. Hence, any $v\in V$ we can be uniquely
written as $v=v_1+\ldots+v_n$
where $v_i\in\ell_i,\,i=1,\ldots,n.$
Such a vector $v$ is a cyclic vector for $x$
if and only if none of the $v_i$'s vanish.

We put $U:= \{(x,v)\in \eer\mid x\in\grs\}.$
Thus, $U$ is an affine $G$-stable Zariski open dense
subset in $\eer$, and the geometric quotient morphism
$p: \eer\to\h/W$ restricts
to a geometric quotient morphism
$p: U\to \hreg/W$.

Now, the group $G$ acts naturally on
the algebra $\ke$.
 We observe that
the first projection
$\EE=\g\times V\to\g$ induces an isomorphism
of $G$-invariants $\k[\g]^G\iso \ke^G$,
since the center of $G$ acts trivially
on $\k[\g]$ and nontrivially
on any homogeneous polynomial $f\in\k[V]$ such that $\deg f >0$.

Let $\Delta_\g$ denote  the second order
Laplacian on $\g$ associated to a
nondegenerate invariant bilinear form.
We will identify $\Delta_\g$ with the operator
 $\Delta_\g\otimes 1\in \D(\g)\otimes\D(V)=\D(\EE)$ acting
 trivially in the $V$-direction. Restricting the latter
differential operator to $U$, we may view
 $\Delta_\g$ as an element of the
algebra $\G(U,\,\D_U/\D_U\cd I_c)^G$.

Write $x_1,\ldots,x_n$ for coordinates in $\h=\k^n$,

\begin{prop}\label{prop}
For any $c\in \k$, there is a natural filtration preserving
algebra
isomorphism 
$\dis\Psi_c:\G(U,\,\D_U/\D_U\cd{I}_c)^G\iso \D(\hreg)^W,$
that
reduces to the `Chevalley restriction' map:
$$f\mto f\big|_{\h\times\{0\}},\quad\k[\g\times V]^G\too
\k[\hreg\times\{0\}]^W=\k[\hreg]^W,
$$
on polynomial\footnote{notice that although
$\hreg\times\{0\}$ is {\em not} a subset of $U$
the restriction to $\hreg\times\{0\}$ is well-defined
for a polynomial on $\EE=\g\times V$.} zero order differential operators,
and such that
 $\,\Psi_c(\Delta_\g)=\cm_c$, where
\begin{equation}\label{H}
\cm_c=\sum\nolimits_j\;\frac{\partial^2}{\partial x_j^2}\;-\;
\sum\nolimits_{i\ne j}\;\frac{c(c+1)}{(x_i-x_j)^2}\,
\end{equation}
is the Calogero-Moser operator with rational potential,
corresponding to the parameter~$c$.
\end{prop}
\begin{proof}[Sketch of Proof.]
We restrict the isomorphism of Proposition
\ref{U_prop} to $U\sset\eer$, equivalently, we
apply formula \eqref{XA} to the geometric quotient
morphism $p: U\to\hreg/W$
and to  the character $\chi_c: \g\to\k$.
This way, we deduce an  algebra isomorphism
\beq{double}
\G(U,\,\D_U/\D_U\cd{I}_c)^G\iso \D(\hreg/W)
\eeq

Now, the natural projection $\hreg\to\hreg/W$
 is
a Galois covering with Galois group $W$.
Therefore, pull-back via 
the projection gives rise to  a canonical
isomorphism $\D(\hreg/W)\cong\D(\hreg)^W$.
Thus, composing with \eqref{double}
yields an   algebra isomorphism
$\Psi'_c: \G(U,\,\D_U/\D_U\cd{I}_c)^G\iso\D(\hreg)^W$.

Finally, let  $R_{+}$ be the
set of  positive roots
of our root system $R\sset \h^*$ of type $\mathbf{A}_{n-1}$,
 and set
$\delta:=\prod_{\alpha\in R_{+}}\alpha$.
We conjugate the map $\Psi'_c$ by $\delta$.
That is, for any  
$u\in \G(U,\,\D_U/\D_U\cdot{I}_c)^G,$
let $\Psi_c(u)$ be a differential operator
on $\hreg$ given by $\Psi_c(u):=M_\delta\ccirc\Psi_c'(u)\ccirc
M_{1/\delta}$, where $M_f$ denotes
the  operator
of multiplication by a 
function $f\in\k[\hreg].$
The map $u\mapsto  \Psi_c(u)$ thus defined
gives the isomorphism $\Psi_c$ of the Proposition.

The equation $\Psi_c(\Delta_\g)=\cm_c$ is verified by a direct
computation similar
to one  in the proof of \cite[Proposition 6.2]{EG}.
We leave details to the reader.
\end{proof}

We will need the following analogue of 
the surjectivity part of
\cite[Corollary~7.4]{EG}.

\begin{prop}\label{surject} Let $\k={\mathbb{Q}}$ be the field
of rational numbers. Then, for all $c\in\k$,
the algebra $\B_c$
is contained in the image of the following composite
map, cf. Proposition \ref{prop}:
$$\overline{\Psi}_c:\
\xymatrix{
\G\bigl(\EE,\,\D_\EE/\D_\EE\cd I_c\bigr)^G
\ar[rrr]^<>(0.5){\text{restriction}}&&&
\G(U,\,\D_U/\D_U\cd{I}_c)^G\ar[r]^<>(0.5){\Psi_c}_<>(0.5){\sim}&
\D(\hreg)^W}.
$$
\end{prop}

\begin{proof} We repeat the argument used in \cite{EG}, which
is quite standard. 

Specifically, the algebra $\B_c$ contains a
subalgebra $\mathsf{C}_c\sset\B_c$ formed by so-called
{\em Calogero-Moser integrals}. The algebra
$\mathsf{C}_c$ is a commutative algebra
containing the Calogero-Moser operator $\cm_c$, and isomorphic
to $(\sym\h)^W$, due to a result by Opdam.
Moreover,
the associated graded map corresponding to
the imbedding $\mathsf{C}_c\into\B_c$ induces
an isomorphism, cf. \eqref{gr_dunkl}:
$$
\gr\mathsf{C}_c\cong(\sym\h)^W\into
(\sym(\h\oplus\h^*))^W=\k[\h^*\times\h]^W.
$$

Observe next that the imbedding
$\g= \g\times\{0\}\into\g\times V=\EE$
induces an isomorphism
$(\sym\g)^G\iso (\sym\EE)^G,$ very similar to the
isomorphism $\k[\g]^G\iso\ke^G$ explained in Sect.~\ref{hc_sec}.
Thus, we identify $(\sym\g)^G$ with $(\sym\EE)^G$,
and view the latter as a subalgebra
in $\D(\EE)^G$ formed by constant coefficient
differential operators.
Clearly, this is a commutative subalgebra
that contains $\Delta_\g=\Delta_\g\otimes 1\in\D(\EE)^G$.

The homomorphism $\Psi_c$ of Proposition \ref{prop},
hence the composite map $\overline{\Psi}_c$
of Proposition \ref{surject}, takes the algebra
$(\sym\g)^G$, viewed as  subalgebra
in $\G(\EE,\,\D_\EE/\D_\EE\cd{I}_c)^G$,
to a commutative subalgebra of $\D(\hreg)^W$ containing
$\cm_c$. Further, one proves by a standard argument
that $\overline{\Psi}_c\bigl((\sym\g)^G\bigr)\sset \mathsf{C}_c$,
cf. e.g. \cite{BEG},
moreover, the induced map
$\gr\overline{\Psi}_c:\
 \gr(\sym\g)^G\map \gr\mathsf{C}_c=(\sym\h)^W$
is the Chevalley isomorphism $(\sym\g)^G\iso(\sym\h)^W$.
It follows that the map $\overline{\Psi}_c$ induces
an isomorphism 
$(\sym\g)^G\iso \mathsf{C}_c$.

Now the algebra  $\B_c$ is 
known to be generated by  the two
subalgebras $\k[\h^*]^W$ and $\mathsf{C}_c$
(since 
$\gr\B_c=\k[\h^*\times\h]^W$, viewed
as a Poisson algebra with respect to the
natural Poisson structure on $\h^*\times\h=T^*\h$,
is known to be generated by the two
subalgebras $\k[\h^*]^W$ and $\k[\h]^W$).
By Proposition \ref{prop}, we have
$\overline{\Psi}_c(\ke^G)=\k[\h^*]^W$ and as
we have explained above, one also has
$\overline{\Psi}_c\bigl((\sym\g)^G\bigr)=\mathsf{C}_c$.
We conclude that $\B_c$ is 
equal to the subalgebra in $\D(\hreg)^W$
generated by $\overline{\Psi}_c(\ke^G)$ and $\mathsf{C}_c,$
hence, is contained in the image of the map
$\overline{\Psi}_c$.
\end{proof}

\section{An Azumaya algebra on the Hilbert scheme}\label{hilb_sec}
\subsection{Nakajima construction reviewed.}\label{calo}
We keep the notations of \S\ref{hc_sec}. In 
particular, we have a vector space $V$ over $\k$,
such that $\charp\k> n\geq 2,$ where $n=\dim V.$
We put $G:=\GL(V)$ and $\g:=\Lie G=\gl(V).$
We will freely identify $\g^*$
with $\g$
via the   pairing
 $\g\times\g\to\k,\,(x,y)\mto\frac{1}{n}\Tr(x\cdot y)$.
 
The group $G$ acts naturally on $V$
and also on  $\g$, via the adjoint action.
We consider the $G$-diagonal action on 
the vector space $\EE=\g\times V$,
and the corresponding
 Hamiltonian $G$-action on the cotangent
bundle:
$T^*\EE= \EE^*\times\EE
\cong\g\times\g\times V^*\times V.$
The  moment map for this action
is given by the formula
\beq{moment}
 \mu:\ T^*\EE=\g\times\g\times V^*\times V\too \g^*\cong \g,\quad
(x,y,\vc,v) \mto [x,y]+\vc\otimes v\in \g.
\eeq

Recall the notation introduced in \ref{IJc}.
Observe that the Lie algebra
homomorphism  $\chi_c=c\cdot\Tr\in\g^*$ corresponds, under the identification
$\g^*\cong\g$, to the element
$c\cdot\Id_V\in\g$.

Following Nakajima, we introduce the set
$$
\fM_c:=\mu\inv(\chi_c)=\big\{(x,y,\vc,v)\in \g\times\g\times V^*\times V
\mid [x,y]+\vc\otimes v=c\cdot\Id_V\big\}.
$$
This is an  affine algebraic variety
 equipped with a natural $\GL(V)$-action.

If  $c\neq 0$, then $\fM_c$ is
 known, see \cite{Na1}, \cite{Wi}, to be smooth,
 moreover, the $G$-action on
$\fM_c$ is free. The quotient
$\mm_c:=\fM_c/G$ is a well-defined smooth affine algebraic
variety of dimension $2\dim V$, called
{\em Calogero-Moser space}. It was first considered
in \cite{KKS}, and studied in \cite{Wi}, cf. also \cite{Na1}.
By definition, $\mm_c$ is
the Hamiltonian reduction of $T^*\EE$ with respect to
the 1-point $G$-orbit $\chi_c\in\g^*$.
The standard symplectic structure on the cotangent
bundle thus induces a symplectic structure
on $\mm_c$.

If $c=0$, then the set $\fM_0$ is not smooth,
and $G$-action on  $\fM_0$ is not free.
Let $\mss$ be the subset of `stable points' formed by quadruples
$(x,y,\vc,v)\in \g\times\g\times V^*\times V$ such that $v\in V$ is a cyclic vector for 
$(x,y)$, i.e., such that there is no nonzero
proper subspace $V'\sset V$ that contains $v$
and that is  both $x$- and $y$-stable.
Then, $\mss$ is known to be a smooth 
Zariski open $G$-stable subset in $\fM_0$.
Moreover, the differential of
the moment map $\mu$, see \eqref{moment},
is known, cf. \cite{Na1}, to be surjective at any point of  $\mss$,
and the $G$-action on $\mss$ is free.

The following  description of $\Hilb$,
 the  Hilbert scheme of 
zero-dimensional length $n$ subschemes in the
affine plane $\af^2$,
is essentially due to Nakajima \cite{Na1}.

\begin{prop}\label{nak1}
There exists a smooth geometric quotient morphism
$\mss$ $\to \Hilb$.
\end{prop}

\begin{rem}
It is known that $\Hilb$ is a smooth connected (non-affine) 
algebraic variety
of dimension $2\dim V$.
\erem

Thus, the Hilbert scheme $\Hilb$
may be viewed as a `Hamiltonian reduction'
of $T^*\EE$ at
the 1-point $G$-orbit $\{0\}\sset\g^*$.
In particular, $\Hilb$ has a natural symplectic structure.

It is known that for any quadruple
$(x,y,\vc,v)\in \fM_0$, the operators $x,y$ can be
put simultaneously in the upper-triangular form.
Hence, the diagonal components
of these two operators give a pair of
elements $\diag x,\diag y\in\h$, well defined up to simultaneous 
action of $W=S_n$. The assignment
$(x,y,\vc,v)\mto (\diag x,\diag y)$ clearly descends
to a morphism $\Upsilon :\Hilb\map$ ${(\h\oplus\h)/W}$, called 
{\em Hilbert-Chow morphism}.
It is known that the  Hilbert-Chow morphism 
induces an algebra isomorphism
\beq{weyl}
\G\bigl(\Hilb,\,\oo_{_{\Hilb}}\bigr) \stackrel{\Upsilon_*}\iso 
\G\bigl((\h\times\h)/W,\,\oo_{_{(\h\times\h)/W}}\bigr)=
\k[\h\oplus\h]^W.
\eeq

\subsection{The Azumaya algebra.} Let
$\as: \BA^1\to\BA^1,\,c\mapsto\varkappa(c)=c^p-c$
be the classical Artin-Schreier map. This map is related to
the map $\as: \BX(\g)\to \g\ttt$ defined  in  \eqref{as}
by the formula
 $\varkappa(\chi_c)=\chi_{\varkappa(c)}=
\varkappa(c)\cdot \Tr\tw.$

We introduce the following simplified notation for the
 scheme-theoretic fiber of the moment
map $\mu\tw: T^*X\tw\to \g\ttt$ over the point $\varkappa(\chi_c)$:
$$T\ttt_{\varkappa(c)}:=T\ttt_{\varkappa(\chi_c)}=
[\mu\tw]\inv(\chi_{\varkappa(c)})=\fM_{\varkappa(c)}\tw.
$$

We are going to apply the general Hamiltonian
reduction procedure of Sect. \ref{some} to the
algebraic group $A:=G$, the Lie algebra character $\chi=\chi_{\varkappa(c)},$
and  the natural $G$-action on the variety $X:=\EE$.

If $\varkappa(c)\neq 0$, then the variety
 $T\ttt_{\as(c)}=\fM_{\varkappa(c)}\tw$
is smooth, the $G$-action on this  variety
is free, and there is a smooth geometric quotient
map $\fM_{\varkappa(c)}\too
\mm_{\varkappa(c)}$, where $\mm_{\varkappa(c)}$
is  the Calogero-Moser
variety with parameter $\varkappa(c)=c^p-c$.
Thus, the construction  of Sect. \ref{some},
applied to $\M=T\ttt_{\as(c)}$, produces 
an Azumaya algebra $\az_c:=\az_{\chi_c}$ on  $\mm_{\as(c)}$.

Assume now that $\varkappa(c)=c^p-c=0$, that is,
the element $c\in\k$ is contained in the finite
subfield $\BF_p\sset\k$. Then,
we apply the construction
 of Sect. \ref{some} to the open
subset $\M:=\mss\sset\mu\inv(0)$,
of stable points.
By Proposition \ref{nak1}, we  obtain
an Azumaya algebra $\az_c:=\az_{\chi_c}$ on  the Frobenius
twist of $\Hilb$, to be denoted $\hilbb$.

For any $c\in\k$, the Azumaya algebra $\az_c$
comes equipped with canonical algebra
homomorphism
\beq{DtoA}
\Xi_c: \ \G\bigl(\EE,\,\D_\EE/\D_\EE\cd I_c\bigr)^G
\stackrel{\eqref{gl_az}}\tooo \G(\mm_{\varkappa(c)}\tw,\,\az_c)
\eeq

\subsection{A Harish-Chandra homomorphism for the Azumaya algebra.}
We are now going to construct 
a Harish-Chandra homomorphism for the Azumaya algebra 
$\az_c$.

Recall the open subset $U\sset\EE$ formed
by the pairs $(x,v)\in\grs\times V$  such that $v$ is a cyclic vector
for $x$.
\begin{prop}\label{phi_az} For any $c\in\k$,
there is an
algebra homomorphism $\Psi_c^\az$ making the following
diagram commute:
$$
\xymatrix{
\G\bigl(\EE,\,\D_\EE/\D_\EE\cd I_c\bigr)^G
\ar[rrrr]_<>(0.5){\eqref{DtoA}}^<>(0.5){\Xi_c}
\ar[d]_<>(0.5){\text{restriction}}
&&&&
\G(\mm\tw_{\varkappa(c)},\,\az_c)
\ar[d]^<>(0.5){\Psi_c^\az}
\\
{\G\bigl(U,\,\D_U/\D_U\cd I_c\bigr)^G}
\ar[rrrr]_<>(0.5){{\text{Proposition \ref{prop}}}}^<>(0.5){\Psi_c}
&&&&\D(\hreg)^W.
}
$$
\end{prop}
\begin{proof}
We only consider the most interesting case $\as(c)=0$.

The action of $G$ on $U$ induces a Hamiltonian $G$-action on $T^*U$.
Let $\mu_U: T^*U\to \g^*$ be the corresponding moment map,
and set $\UU:= \mu_U\inv(0)$.
 Since $U$ is an open subset of $\EE$,
the map $\mu_U$  clearly equals
 the restriction
of the moment map $\mu: T^*\EE\map\g^*$ to
the open subset  $T^*U=\EE^*\times U\sset
\EE^*\times \EE$.
Thus,
we have $\UU=\mu_U\inv(0)=(\EE^*\times U)\,\bigcap\,\mu\inv(0).$

It is crucial  for us that one has an open inclusion
$$\UU\;\sset\;
\mss=\M.$$
This trivially follows from definitions since
a vector $v\in V$ which is cyclic for $x\in \g$
is necessarily also cyclic for any pair of the form $(x,y)\in \g\times\g$.

Thus, we have the diagram
$$
\xymatrix
{
{\UU\tw_{_{}}\enspace}\ar@{^{(}->}[rrr]_<>(0.5){\text{open imbedding}}&&&
{\enspace\M\tw=[\mss]\tw_{_{}}\enspace}
\ar@{->>}[rrrr]^<>(0.5){\varpi}_<>(0.5){\text{geom. quotient map}}&&&&
{\enspace\mm\tw_{_{}}=\hilbb}.
}
$$

Restricting the sheaf $\ee_c:=\ee_{\chi_c}$, cf. \eqref{Dchi},
from $[\mss]\tw$ to $\UU\tw$ yields a
 $G\tw$-equivariant
algebra map 
$\G([\mss]\tw,\,
\,{\ee_c})\map\G(\UU\tw,\,
\,{\ee_c}).
$
Further, applying Lemma \ref{stages}
to $X:=U$ and $A:=G$, we get an algebra isomorphism 
\beq{stages_U}
\G(U,\,\D_U/\D_U\cd  I_c)^G
\cong\,\G\bigl([\mu_U\inv(0)]\tw,\,
\,{\ee_c}\bigr)^{G\tw}=\G\bigl(\UU\tw,\,
\,{\ee_c}\bigr)^{G\tw}.
\eeq

By definition, we have
$\dis\az_c:=\,
\varpi\idot\bigl(\ee_c\big|_{[\mss]\tw}\bigr)^{G\tw}.$
We obtain the following 
chain of algebra  homomorphisms
\begin{align*}
\G(\mm\tw,\az_c)&=\G\left(\mm\tw,
\,\varpi\idot\bigl(\ee_c\big|_{[\mss]\tw}\bigr)^{G\tw}\right)\\
&=
\G\left(\mm\tw,
\,\varpi\idot\bigl(\ee_c\big|_{[\mss]\tw}\bigr)\right)^{G\tw}
=\G\bigl([\mss]\tw,\ee_c\bigr)^{G\tw}\\
&\stackrel{\text{restriction}}\tooo
\G\bigl(\UU\tw,\ee_c\bigr)^{G\tw}
\underset{^{\eqref{stages_U}}}{\;\iso\;}
\G(U,\,\D_U/\D_U\cd I_c)^G\enspace
\underset{^{\text{Prop. \ref{prop}}}}{\stackrel{\Psi_c}\tooo}\enspace
\D(\hreg)^W.
\end{align*}

We let $\Psi^\az_c$ be the composite homomorphism.
Commutativity of the diagram of the Proposition is
immediate from the construction above.
\end{proof}

\section{Localization functor for Cherednik algebras}
\subsection{From  characteristic zero to characteristic $p$.}\label{mod_p}
We begin by reminding the general technique of transferring
various results valid over fields of characteristic zero to
similar results in characteristic $p$, provided
$p$ is sufficiently large.

We fix $c=a/b\in\BQ$ with $b>0$. For any prime $p > b$,
reducing modulo $p$,
we may (and will) treat $c=a/b$ as an element of
${\mathbb{F}}_p$. We let $\k_p\supset{\mathbb{F}}_p$
denote an algebraic closure of ${\mathbb{F}}_p$, and consider the
corresponding  $\k_p$-algebras 
$\hh_c$ and $\e\hh_c\e$.

We begin with the following characteristic $p$ analogue of
 Proposition \ref{surject}.
\begin{lemma}
\label{surject_p}
 Fix $c=a/b\in\BQ$.
For all sufficiently large primes $p$,
we have an inclusion of
$\k_p$-algebras $\B_c\sseq \im(\overline{\Psi}_c)$.
\end{lemma} 
\begin{proof} For $c=a/b\in \BQ$, we  consider
the ring $\BZ[\frac{1}{n!},c]=\BZ[\frac{1}{n!},\frac{1}{b}]$
obtained by inverting $n!$ and $b$.
The algebras
$\G(\EE,\,\D_\EE/\D_\EE\cd{I}_c)^G$ and $\D(\hreg)^W$,
are both defined over $\Z$, hence have 
natural  $\BZ[\frac{1}{n!},c]$-integral structures.
We will denote the corresponding  $\BZ[\frac{1}{n!},c]$-algebras
by the same symbols.
Therefore, we may (and will) consider  the $\BZ[\frac{1}{n!},c]$-integral
version
$\overline{\Psi}^\Z_c:\G(\EE,\,\D_\EE/\D_\EE\cd{I}_c)^G\map \D(\hreg)^W$,
of the homomorphism
of Proposition \ref{prop}.
Thus, $\overline{\Psi}^\Z_c$ is a
 homomorphism of $\BZ[\frac{1}{n!},c]$-algebras.

Further, we may choose 
a finite set of generators in the $\BQ$-algebra  $\B_c$ 
in such a way that the $\BZ[\frac{1}{n!},c]$-subalgebra
generated by this set is contained in the
$\BZ[\frac{1}{n!},c]$-integral form of  $\D(\hreg)^W$.
Denote this $\BZ[\frac{1}{n!},c]$-subalgebra
 by  $\B^\Z_c$. Thus, we get
a diagram of $\BZ[\frac{1}{n!},c]$-algebra maps
$$
\xymatrix{
\G\bigl(\EE,\,\D_\EE/\D_\EE\cd I_c\bigr)^G
\ar[rr]^<>(0.5){\overline{\Psi}^\Z_c}&&\D(\hreg)^W&&
\ar@{_{(}->}[ll]_<>(0.5){j}\B_c^\Z,
}
$$
where $j$ denotes the inclusion.

Proposition \ref{surject} says that 
$
\BQ\bigotimes_{\BZ[\frac{1}{n!},c]}\B_c\sseq
\BQ\bigotimes_{\BZ[\frac{1}{n!},c]}\im(\overline{\Psi}^\Z_c).
$
Since all the algebras involved are finitely generated,
it follows that there exists an integer
$q\in \BZ$ such that
$\B_c^\Z\sset\im(\overline{\Psi}^\Z_c)[\frac{1}{q}].$
Thus, for all primes $p> n$ which do not divide
$q$, reducing the above inclusion modulo $p$,
we get
$\k\bigotimes_{\BZ[\frac{1}{n!},c]}\B_c\sseq
\k\bigotimes_{\BZ[\frac{1}{n!},c]}\im(\overline{\Psi}^\Z_c).$
The Lemma is proved.
\end{proof}

We now consider the algebras  $\hh_c$ and $\ehe$
over the ground field $\BQ$ of the rational numbers,
and let $\hh_c\e\hh_c$ be the two-sided ideal in $\hh_c$ generated by 
the idempotent $\e$.

We will use the following 
result from \cite{g}.
\begin{prop}\label{gordon}
For any $c\in\QG,$ see \eqref{qgood},
we have $\hh_c=\hh_c\e\hh_c$. Thus,
the $\ehe\dash\hh_c$-bimodule $\e\hh_c$ provides a
Morita equivalence between
the algebras $\hh_c$ and $\ehe$.\qed
\end{prop}

We are going to deduce a similar result in characteristic $p$,
which reads

\begin{corollary}
\label{gordonp} Given $c\in\QG$, there exists a constant
$d=d(c)$
such that for all primes $p> d(c)$,
the $\ehe\dash\hh_c$-bimodule $\e\hh_c$ provides a
Morita equivalence between the $\k_p$-algebras 
 $\hh_c$ and~$\ehe$.
\end{corollary}
\begin{rem} We emphasize that, in this Corollary and in various
other results below, {\em a rational value}
of the parameter $c$ must be fixed first. The choice of $c$
dictates a lower bound $d(c)$ for allowed primes $p$,
and only after that one considers the corresponding
Cherednik algebras over $\k_p$. Thus, if $c\in \BQ$ and
$p$ have been chosen as above, and $c'\in\BQ$ is such that
$c'=c \mod p$, then $\hh_c\cong \hh_{c'}$ as
$\k_p$-algebras; yet, it is quite possible that we have
$d(c)< p < d(c')$, hence, the results of this section
do not apply for $\hh_{c'}$ viewed as a $\k_p$-algebra.
\erem

The proof of the Corollary
 will exploit the following standard result of commutative algebra,
\cite{Gr}, Expose IV, Lemma 6.7.
\smallskip

\noindent
{\bf Generic Flatness Lemma.}
 {\em Let $A$ be a commutative noetherian integral
domain,
$B$ a (commutative) $A$-algebra of finite type, and
$M$ a finitely generated $B$-module. Then, there is a nonzero element
$f\in A$ such that $M_{(f)}$, the localization of $M$,  is a free
$A_{(f)}$-module.}\qed
\smallskip

\begin{proof}[Proof of Corollary \ref{gordonp}.]
Let
${\mathbf{c}}$ denote an independent variable, 
and let $\BZ[\frac{1}{n!}][{\mathbf{c}}]$, 
be the localization of the
polynomial ring $\BZ[{\mathbf{c}}]$ at the number $n!\in\BZ$.
Given a nonzero element $f\in \BZ[\frac{1}{n!}][{\mathbf{c}}]$
and a $\BZ[\frac{1}{n!}][{\mathbf{c}}]$-module
$M$, we write $M_{(f)}$ for the
localization of $M$ at $f$, a
module over the localized ring $\BZ[\frac{1}{n!}][{\mathbf{c}}]_{(f)}$.

Let  $\hh$ and $\e\hh\e$
be the  {\em universal}
 Cherednik algebras,
viewed as algebras  over the ground ring
$\BZ[\frac{1}{n!}][{\mathbf{c}}]$.
We first establish the following

\begin{claim}
\label{Quillen}
There exists a polynomial $f\in\BZ[\frac{1}{n!}][{\mathbf{c}}]$ such that 
$\hh_{(f)}=(\heh)_{(f)}$ and, moreover,
such that $f(c)\neq 0$ for any $c\in \QG$.
\end{claim}

To prove the Claim, consider
the standard increasing 
filtration $F\idot$ on $\hh$ and the induced  filtration
$F\idot(\hh/\heh)$ on
the quotient algebra  $\hh/\heh$.\footnote{The argument
below is similar to \cite{Q}.}
The associated graded $\gr^F(\hh/\heh)$,
is a finitely generated  $\BZ[\frac{1}{n!}][{\mathbf{c}}]$-algebra.
Clearly, it suffices to show that
this algebra 
vanishes generically over $\Spec\BZ[\frac{1}{n!}][{\mathbf{c}}]$. 
To this end, observe that  $\gr(\hh/\heh)$ is a finitely generated
module over
the graded algebra $\gr(\hh)$, which is a quotient of the smash-product algebra
$$\BZ[\frac{1}{n!}][{\mathbf{c}}][\h\times \h^*]\ltimes\BZ[\frac{1}{n!}][{\mathbf{c}}][S_n]$$
(here we regard $\h$ as a free rank $n$ module over $\BZ[\frac{1}{n!}][{\mathbf{c}}]$;
its dual $\h^*$ is also free of rank $n$; by the Poincare-Birkhoff-Witt
theorem proved in \cite{EG},
the algebra $\gr(\hh)$ and the above smash-product become
isomorphic after tensoring with $\BQ$). 
We deduce that
 $\gr(\hh/\heh)$ is a finitely generated module over
the {\em commutative} algebra $\BZ[\frac{1}{n!}][{\mathbf{c}}][\h\times \h^*]$.

The Generic Flatness Lemma implies that there exists
a nonzero polynomial 
$f\in\BZ[\frac{1}{n!}][{\mathbf{c}}]$ such
 that $\gr(\hh/\heh)_{(f)}$
is free over $\BZ[\frac{1}{n!}][{\mathbf{c}}]_{(f)}$. 
On the other hand, 
for $c\in\QG$,
by Proposition \ref{gordon} we have 
$\BQ\otimes\hh_c=\BQ\otimes(\hh_c\e\hh_c)$;
that is, the fiber
of $\BQ\otimes\gr(\hh/\heh)$ over $c\in\QG$
vanishes. 
We conclude that $\hh=\hh\e\hh$ holds over a nonempty Zariski
open subset of
$\on{Spec}(\BZ[\frac{1}{n!}][{\mathbf{c}}])$ that has a nontrivial
intersection with any closed subscheme
${\{{\mathbf{c}}=c\}},$
$c\in\QG$.
This proves Claim \ref{Quillen}.

We complete the proof of Corollary \ref{gordonp} as follows.
By Claim \ref{Quillen}, there exist 
 a  polynomial $f\in\BZ[\frac{1}{n!}][{\mathbf{c}}]$ 
 and elements $h'_i,h''_i\in\hh,\, i=1,\ldots, m,$
such that 
\beq{gor}
f\cdot \boldsymbol{1}_\hh=\sum_{i=1}^m\,
h'_i\cd\e\cd h''_i\quad\text{holds in}\quad\hh,
\eeq
where $\boldsymbol{1}_\hh$ denotes the unit of the 
$\BZ[\frac{1}{n!}][{\mathbf{c}}]$-algebra
$\hh$. We may specialize this equation
at any rational value $c=a/b\in\BQ$ to obtain
a similar  equation for the corresponding
$\BZ[\frac{1}{n!},\frac{1}{b}]$-algebras.
If $c\in\QG$  then, according to  Claim \ref{Quillen}, we may 
further assume that $f(c)=k/l \neq 0$.

Now, let $p$ be a prime such that $p> \operatorname{max}\{n,k,l\}$.
Reducing (the specialization at $c$ of) equation \eqref{gor} modulo
$p$,  for the corresponding
$\BF_p$-algebras
we get $f(c)\cdot \boldsymbol{1}_{\hh_c}=$
${\sum_{i=1}^m\,
h'_i\cdot\e\cdot h''_i.}$
Thus $f(c)=k/l$ is a nonzero, hence, invertible element
in  $\BF_p$ and,
since $\BF_p\sset\k_p$, in the $\k_p$-algebra $\hh_c$
 we obtain $\boldsymbol{1}_{\hh_c}=\frac{1}{f(c)}\sum_{i=1}^m\,
h'_i\cdot\e\cdot h''_i$. 

Thus, we have proved that $\boldsymbol{1}_{\hh_c}\in
\hh_c\e\hh_c$, and the first statement of the Corollary
follows. It is well-known that this implies
the last statement of the Corollary as well.
\end{proof}

\subsection{Localization of the Spherical subalgebra.}
We have the following 
$\az_c$-version
of Proposition \ref{surject}.

\begin{theorem}\label{rad_p}   Fix $c=a/b\in\BQ$. 
For all sufficiently large primes $p$, we have:

The 
image of the Harish-Chandra
homomorphism $\Psi_c^\az$ of Proposition \ref{phi_az}
is equal to the subalgebra $\B_c\sset\D(\hreg)^W$
(algebras over $\k_p$).
Moreover, the resulting map gives an algebra
 isomomorphism
$\Psi_c^\az:\G(\hilbb,\az_c){\iso\B_c.}$
\end{theorem}

Proof of this Theorem will be given later in this section.

\begin{rem} A similar construction also produces
an isomorphism $\Psi_c^\az:\ \G(\mm\tw_{\varkappa(c)},\az_c)\iso
\B_c,$ for {\em all} $c\in\k_p$ (not only for $c\in \mathbb{F}_p$).
The proof of this generalization is similar to the
proof of Theorem \ref{rad_p}, but involves {\em twisted}
differential operators and {\em twisted cotangent bundles.}
It will be presented elsewhere.
\erem

Composing the isomorphism  of Theorem \ref{rad_p}
with the inverse of the Dunkl representation \eqref{dunkl},
we obtain the following Azumaya version of the {\em Spherical
Harish-Chandra isomorphism} considered in 
\cite{EG}
\beq{spher}
\Phi_c^\az:\ \G(\hilbb,\,\az_c)
\stackrel{\Psi^\az_c}{\isto}
\B_c
\stackrel{(\Th_c)\inv}{\isto}
\ehe.
\eeq

Thus, we have proved part (i) of the following theorem,
which is one
 of the main results of the paper
\begin{theorem}\label{Ehe}  Fix $c=a/b\in\BQ.$
Then there exists a constant
$d=d(c)$
such that for all primes $p> d(c)$, we have:

\vi  The composite morphism in \eqref{spher} yields a $\k_p$-algebra isomorphism
$\Phi_c^\az:$
$\ehe\iso\G(\hilbb,\,\az_c).$

\vii  $\;H^i(\hilbb,\az_c)=0$ for all $i>0$.

\viii  If  $c\in\QG$, then the (derived)
global sections 
functor $\RG:\ D^b(\az_c\dash\on{Mod})\to D^b(\e\hh_c\e\dash\on{Mod})$
is an equivalence of bounded derived categories.
\end{theorem}

\begin{proof} Part (i) follows from
Theorem \ref{rad_p} (to be proved below),
and part (ii) is a consequence of
Theorem \ref{roma}(ii).
We are going to  deduce part (iii) of Theorem \ref{Ehe} from
Theorem \ref{roma}(iii). To do so,
we need to know that the algebra
$\ehe\cong\G(\hilbb,\,\az_c)$  has
finite homological dimension. But this follows from
the Morita equivalence of Corollary ~\ref{gordonp}, since
the algebra $\hh_c$ is known to have
finite homological dimension, which is equal to $2n$,
cf.~\cite{EG}.
\end{proof}

\begin{proof}[Proof of Theorem \ref{rad_p}]\label{proof_rad}
Fix $c\in \BQ$, and let $p\gg 0$ be such that Lemma \ref{surject_p}
holds for $p$. 

Recall an increasing
filtration on the algebra $\G(\hilbb,\,\az_c)$,
introduced in section~\ref{deform},
such that for the associated graded
algebra we have
$\gr^\rr\G(\hilbb,\,\az_c)=\G(\Hilb,\,\oo_{_{\Hilb}}),$
see {Proposition \ref{gr_az}(ii).

We consider the commutative diagram of Proposition \ref{phi_az}.
All maps in that diagram are filtration preserving,
and the corresponding commutative diagram of associated
graded  maps reads
\beq{incl1}
\xymatrix
{
\gr\G\bigl(\EE,\,\D_{\EE}/\D_{\EE}\cd I_c\bigr)^G
\ar[rrrr]_<>(0.5){\eqref{DtoA}}^<>(0.5){\gr\Xi_c}
\ar[d]_<>(0.5){{\text{restriction}}}^<>(0.5){{\op{res}}^{\EE}_U}
&&&&
\gr^\rr\G(\hilbb,\,\az_c)
\ar[d]^<>(0.5){\gr\Psi_c^\az}
\\
\gr\G\bigl(U,\,\D_U/\D_U\cd I_c\bigr)^G
\ar[rrrr]_<>(0.5){\sim}^<>(0.5){\gr\Psi_c}
&&&&
\gr\D(\hreg)^W=\k[\h^*\times\hreg]^W.}
\eeq

Further, we have the following graded algebra isomorphisms:
\beq{chow_iso}
\xymatrix{\gr^\rr\G\bigl(\hilbb,\,\az_c\bigr)
\ar[rrr]_<>(0.5){\sim}^<>(0.5){{\text{Proposition \ref{gr_az}(ii)}}}&&&
\G\bigl(\Hilb,\oo_{_\Hilb}\bigr)\;\ar[r]_<>(0.5){\sim}^<>(0.5){{\eqref{weyl}}}&
\;\k[\h^*\times\h]^W
}.
\eeq

Moreover,  by going through definitions, it is easy to verify
that the composite map in \eqref{chow_iso}
is equal to the map  $\gr\Psi^\az_c$ in \eqref{incl1}.
It follows, in particular, that $\Psi^\az_c$ is
an injective morphism.

On the other hand,  Lemma  \ref{surject_p} and Proposition
\ref{phi_az}
yield
\beq{incl_main}
\B_c\;\sset\; \im(\overline{\Psi}_c):=\Psi_c\left(
{\op{res}}^{\EE}_U\bigl(\G(\EE,\,\D_{\EE}/\D_{\EE}\cd I_c)^G\bigr)\right)\;\sseq\;
\im(\Psi_c^\az).
\eeq
Hence, using commutativity of diagram \eqref{incl1}, we
obtain the following commutative diagram
 of graded algebra morphisms
\beq{incl2}
\xymatrix{
\k[\h^*\times\h]^W\ar@{=}[d]^<>(0.5){\eqref{chow_iso}}
\ar@{=}[r]&\gr\B_c\ar@{^{(}->}[r]^<>(0.5){\imath}&
\im(\gr{\overline{\Psi}}_c)\ar@{^{(}->}[d]^<>(0.5){\jmath}\\
\gr^\rr\G\bigl(\hilbb,\,\az_c\bigr)\ar@{->>}[rr]^<>(0.5){\gr\Psi_c^\az}
&&
\im(\gr\Psi_c^\az)
.
}
\eeq
We deduce from commutativity of diagram \eqref{incl2} that both
$\imath$ and $\jmath$ must be surjective and,
therefore,
$\im\bigl(\gr\Psi_c^\az\bigr)=\gr\B_c$.
Now, the inclusions in \eqref{incl_main} show that
we must have $\im(\Psi_c^\az)=\B_c$ and,
moreover, the map $\Psi_c^\az$ gives an isomorphism
$\G\bigl(\hilbb,\,\az_c\bigr)$
${\iso\B_c.}$}
\end{proof}

\subsection{Localization of the algebra $\hh_c$.}
\label{ah_main} We introduce the following
{\em  localization functor}
 $\loc :\ N\mapsto\az_c
\stackrel{L}{\otimes}_{\e\hh_c\e}N$, which is
the left adjoint  to the functor 
$\RG:\ D^b(\az_c\dash\on{Mod})\to D^b(\e\hh_c\e\dash\on{Mod})$.
Since $\RG(\az_c)=\e\hh_c\e$,
and the functor $\RG$ is an equivalence
by Theorem \ref{Ehe}, we conclude that
 $\loc (\e\hh_c\e)=\az_c$ is also an equivalence which is a
quasi-inverse
to $\RG(-).$

Observe next that $\e\hh_c$ is a projective $\e\hh_c\e$-module, by
Corollary ~\ref{gordonp}. Hence
we conclude that
$\CR_c:=\loc (\e\hh_c)=\az_c\otimes_{\e\hh_c\e}\e\hh_c$
 is a locally free sheaf of $\az_c$-modules.
Moreover, it is easy to see by looking
at the  restrictions of the associated graded modules
to the
generic locus of $(\h\times\h)/W$ that
the rank of $\e\hh_c$ viewed as a projective
 $\e\hh_c\e$-module  equals $n!$.
Therefore, we deduce that
$\CR_c$ is a vector bundle on $\hilbb$
of rank $n!\cdot p^{2n}$.

We put $\ah_c:=\eend_{\az_c}(\CR_c)$. This is clearly
an Azumaya algebra on $\hilbb$ again,
and the degree of this  Azumaya algebra
is equal to $n!\cdot p^{2n}$.
Further, the left $\hh_c$-action
on each fiber of the sheaf $\CR_c$ induces a natural algebra
map
\beq{ah_map}
\hh_c\too\G(\hilbb,\ah_c).
\eeq

The second main result of the paper reads
\begin{theorem}\label{ah} Fix $c\in\QG$. Then, there exists
a constant $d=d(c)$ such that for all
 primes $p>d(c)$,
we
have

\vi The map \eqref{ah_map}
is an algebra isomorphism, moreover,
$\RG^i(\hilbb,\ah_c)=0,\,\forall i>0$.

\vii The functor $\RG:\ D^b(\ah_c\dash\on{Mod})\to D^b(\hh_c\dash\on{Mod})$
is a triangulated equivalence.
\end{theorem}

\begin{proof} We have
\begin{align}\label{ddproof}
\RG^i(\hilbb,\,\ah_c)\,\,\cong\,\,\RG^i(\hilbb,\,&\eend_\az(\CR_c))
\,\,\cong\,\, \Ext^i_{D^b(\az_c\dash\on{Mod})}(\CR_c,\CR_c)\\
&\xymatrix
{
\ar[rrr]_<>(0.5){\text{Theorem \ref{Ehe}}}^<>(0.5){\sim}&&&
\,\Ext^i_{D^b(\ehe\dash\on{Mod})}(\e\hh_c,\e\hh_c).
}\nonumber
\end{align}
The $\Ext$-group on the right vanishes for all $i>0$
since $\e\hh_c$ is a projective $\ehe$-module. This proves 
the vanishing statement in part (i).
The statement of part (i)  for
 $i=0$ follows from the  isomorphisms:
$$\xymatrix
{\hh_c
\ar[rrrr]_<>(0.5){\sim}^<>(0.5){\text{\cite[Theorem 1.5(iv)]{EG}}}
&&&&
\Hom_{\ehe}(\e\hh_c,\e\hh_c)
\ar[rr]_<>(0.5){\sim}^<>(0.5){\eqref{ddproof}}
&&\G(\hilbb,\,\ah_c).}$$

To prove part (ii) we use a commutative
diagram
$$
\xymatrix
{
D^b(\ehe\dash\on{Mod})\ar@{=}[rrrrr]^<>(0.5){\text{Morita equivalence}}
\ar[d]_<>(0.5){\loc}&&&&&
D^b(\hh_c\dash\on{Mod})
\ar[d]_<>(0.5){\loc}\\
D^b(\az_c\dash\on{Mod})\ar@{=}[rrrrr]^<>(0.5){\text{Morita equivalence}}
&&&&&D^b(\ah_c\dash\on{Mod}).
}
$$
Since the left vertical arrow is an equivalence
by Theorem \ref{Ehe},
it follows that the right vertical arrow is an equivalence
as well. The functor $\RG$ is a right adjoint of $\loc$,
hence, it must also be an equivalence,
which is a quasi-inverse of $\loc$.
\end{proof}

\subsection{Splitting on the fibers of  the Hilbert-Chow map.}
We have the Hilbert-Chow map $\Upsilon: \Hilb\to (\h\times\h)/W$.
Given $\xi\in (\h\times\h)/W$,
write $\op{Hilb}_\xi:=\Upsilon\inv(\xi)$ for the fiber of
$\Upsilon$ over $\xi$, and let ${\wh{\op{Hilb}}}_\xi\tw$ denote
the completion of $\hilbb$ along this fiber, a formal scheme.

The next result is essentially due to \cite{bk}.

\begin{theorem}\label{splitting} For any point
$\xi\in (\h\times\h)/W$, the restriction of
the Azumaya algebra $\ah_c,$ resp. $\az_c$, to
the formal neighborhood of the
fiber $\op{Hilb}_\xi\tw\sset \hilbb$
splits, i.e., there is a vector bundle
$\CV_{c,\xi}$, resp. $\CW_{c,\xi}$,  on ${\wh{\op{Hilb}}}_\xi\tw$ such that one has
$$\ah_c\big|_{{\wh{\op{Hilb}}}_\xi\tw}\,\cong\,
\Bigl(\eend_{_{\oo_{{\wh{\op{Hilb}}}_\xi\tw}}}\CV_{c,\xi}\Bigr)^{\op{opp}},
\quad\text{resp.},\quad
\az_c\big|_{{\wh{\op{Hilb}}}_\xi\tw}\,\cong\,
\Bigl(\eend_{_{\oo_{{\wh{\op{Hilb}}}_\xi\tw}}}\CW_{c,\xi}\Bigr)^{\op{opp}}.
$$
\end{theorem}

The above vector bundle $\CV_{c,\xi}$, resp. $\CW_{c,\xi}$, is called a
{\em splitting bundle} for $\ah_c$, resp. for $\az_c$.

\begin{corollary}\label{rigid}
For any $\xi\in (\h\times\h)/W$ and $i>0$, we have
$\Ext^i(\CV_{c,\xi},\CV_{c,\xi})=0,$ resp.,
$\Ext^i(\CW_{c,\xi},\CW_{c,\xi})=0,$
where the Ext-groups are considered in the
category $\Coh({\wh{\op{Hilb}}}_\xi\tw)$.
\end{corollary}

The rank of a splitting bundle is equal to
the degree of the corresponding Azumaya algebra.
In particular, we have $\rk \CW_{c,\xi}=p^n$ and
$\rk\CV_{c,\xi}=p^n\cd n!$. This suggests
the following
\begin{conj}\label{procesi} For any $\xi\in [(\h\times\h)/W]\tw$, 
there is a vector bundle isomorphism
$\Fr^*\CV_{c,\xi}\cong (\Fr^*\CW_{c,\xi})\otimes(\CP\big|_{\wh{\op{Hilb}}^n_{\Fr(\xi)}})$,
where $\CP$ denotes 
 the {\em Procesi 
bundle}
on $\Hilb$, see \cite{H}.
\end{conj}

In view of Theorem \ref{roma}(iii) it is sufficient to prove 
the Conjecture  for $c=0$, that is, for the case where $\hh_c=\D(\h)\# W$.

\begin{proof}[Proof of Theorem \ref{splitting}.] 
Clearly, it suffices to prove the Theorem
for $\az_c$. We repeat the argument
in the proof of \cite[Proposition 5.4]{bk}.

First, recall that Morita equivalence classes
of Azumaya algebras on a scheme $Y$ are classified
by $Br(Y)$,  the Brauer group of $Y$. Further,
the Brauer group of a local complete $\k$-algebra is known to
be trivial. Thus,  proving the Theorem
amounts to showing that, for any $c\in {\mathbb{F}_p}$, the
class $[\az_c]\in Br(\hilbb)$ belongs to the 
image of the pull-back morphism
$\dis \Upsilon^*: Br\bigl([(\h\times\h)/W]\tw\bigr)
\to Br(\hilbb)$.

To prove this, we consider the following diagram
\beq{TS}
\xymatrix{\TS\;\;\ar@{_{(}->}[d]
\ar@{->>}[rr]_<>(0.5){q_{_S}}^<>(0.5){\text{$W$-covering}}&&
S\;\ar@{_{(}->}[d]&&S^\circ\;\ar@{_{(}->}[d]
\ar[ll]^<>(0.5){\Upsilon_{_S}}_<>(0.5){\sim}\\
[\h\times\h]\tw\ar@{->>}[rr]_<>(0.5){q}^<>(0.5){\text{finite map}}&&
[(\h\times\h)/W]\tw&&\hilbb.\ar[ll]^<>(0.5){\Upsilon}
}
\eeq
In this diagram, $S$ is a  Zariski open dense subset in
$[(\h\times\h)/W]\tw$ such that:

\pb{The Hilbert-Chow map
restricts to an isomorphism $S^\circ:=\Upsilon\inv(S)\iso S$,
to be denoted $\Upsilon_{_S}$,
and}

\pb{The projection  $q: [\h\times\h]\tw \map[(\h\times\h)/W]\tw$
is unramified over $S$. Thus, we have a Galois
covering $ \TS:=q\inv(S) \to S$,
to be denoted $q_{_S}$.}

Our goal is to construct a class
$\beta\in Br\bigl([(\h\times\h)/W]\tw\bigr)$
such that $[\az_c]=\Upsilon^*(\beta)$.
We will follow  the strategy of \cite{bk}.

Let $Y$ be an arbitrary  affine scheme,
and $H_{et}^2(Y,\Gm)\tors$ be the torsion subgroup
of the second \'etale cohomology
of $Y$ with coefficients in the multiplicative group.
By a theorem of Gabber \cite{Ga}, one has an isomorphism
$Br(Y)\cong H_{et}^2(Y,\Gm)\tors$.
Further,
it is a simple matter to see that the norm-map
associated to the projection  $q: [\h\times\h]\tw
\map[(\h\times\h)/W]\tw$
gives rise to a morphism on  \'etale cohomology,
cf. \cite{bk}:
$$q_*:\ H_{et}^2\bigl([\h\times\h]\tw,\Gm\bigr)^W\tors\map
H_{et}^2\bigl((\h\times\h)/W),\Gm\bigr)\tors.$$

Now let
 $\D\tw:=\Fr\idot\D_\h$ be the standard Azumaya algebra
on $T^*\h\tw=[\h\times\h]\tw$, arising from the sheaf
of crystalline differential operators on $\h$, cf. \S2.
The sheaf $\D\tw$ has a natural $W$-equivariant structure,
hence the corresponding class
$[\D\tw]$ is a $W$-invariant class in the Brauer group,
that is an element 
of $\dis Br\bigl([\h\times\h]\tw\bigr)^W
\cong
H_{et}^2\bigl([\h\times\h]\tw,\Gm\bigr)^W\tors.
$
We set 
$$\beta:=q_*([\D\tw])\in
H_{et}^2\bigl([(\h\times\h)/W]\tw,\Gm\bigr)\tors\cong
Br\bigl([(\h\times\h)/W]\tw\bigr).$$

The Theorem would follow provided we show that $\Upsilon^*\beta=[\az_c]$.
We first prove a weaker claim
\beq{weak}
(\Upsilon^*\beta)\big|_{S^\circ}=[\az_c]\big|_{S^\circ}
\quad\text{holds in}\quad Br(S^\circ).
\eeq

To see this, restrict the  Azumaya algebra $\D\tw$
to the open subset $\TS\sset [\h\times\h]\tw$. The map $q_{_S}$,
see \eqref{TS}, is a Galois covering with
the Galois group $W$.
It follows that the sheaf ${\scr B}:=\left((q_{_S})\idot(\D\tw)\big|_{\TS}\right)^W$
is an  Azumaya algebra on $S$. Furthermore, it is immediate from
the construction that, in $Br(S)$, one has an equality
$\beta\big|_S=[{\scr B}]$. Pulling back via the
isomorphism $\Upsilon_S$,
see \eqref{TS}, we deduce that
$\dis (\Upsilon^*\beta)\big|_{S^\circ}=(\Upsilon_S)^*\bigl(\beta\big|_{S}\bigr)=
(\Upsilon_S)^*[{\scr B}]$.

To complete the proof of \eqref{weak}, we 
use  Theorem \ref{roma}(iv)
and deduce that, for all $c\in{\mathbb{F}_p}$, the corresponding
Azumaya algebras $\az_c$ are Morita equivalent,
hence represent the same
class in $Br(\hilbb)$. Thus, we may assume without loss of generality
that
$c=0$. In that case
the corresponding algebra $\ehe$ is
isomorphic to $\D(\h)^W$.
Furthermore, going through the Hamiltonian reduction
construction of the Azumaya algebra
$\az_0$, it is easy to verify that
we have an  Azumaya algebra isomorphism 
$(\Upsilon_{_S})^*{\scr B}
\cong\az_0\big|_{S^\circ}$.
This yields an equality of the corresponding
classes in $Br({S^\circ})$, and \eqref{weak} follows.

To complete the proof of the Theorem we recall the 
well-known result saying that restriction
to a Zariski open dense subset induces
an {\em injective} morphism of the corresponding Brauer groups.
Thus, we have an injection $Br(\hilbb)\into Br(S^\circ),
\, \alpha\mapsto \alpha\big|_{S^\circ}$,
and we have shown above that $(\Upsilon^*\beta)\big|_{S^\circ}=
[\az_0]\big|_{S^\circ}$.
Hence, $\Upsilon^*\beta=[\az_0]=[\az_c],\,\forall c\in{\mathbb{F}_p}$, 
and the Theorem is proved.
\end{proof}

\subsection{Bigrading on  $\wh{\hh}_{c,0}$.}\label{bigrading_sec}
Let $\xi=0$ be the origin of $(\h\times\h)/W$ and
set $\wh{\hh}_{c,0}=\G(\widehat{\op{Hilb}}_0\tw,
\ah_c)$, the completion of the Cherednik algebra $\hh_c$
at the zero central character.
The isomorphism of Theorem \ref{splitting} gives
a continuous (right) $\wh{\hh}_{c,0}$-action on
the  splitting vector bundle $\CV_{c,0}$ by
vector bundle endomorphisms. In particular,
there is an action of the Symmetric group $S_n$ on
 $\CV_{c,0}$.  Thus, for any simple $S_n$-representation $\tau$,
we have the corresponding $\tau$-isotypic component
$\Hom_{S_n}(\tau, \CV_{c,0})$. This isotypic component
is again a vector bundle on $\widehat{\op{Hilb}}_0\tw,$
moreover, the natural evaluation map
gives an $S_n$-equivariant vector  bundle
isomorphism 
\beq{isotyp}
\bigoplus\nolimits_{\tau\in{\mathsf{Irrep}(S_n)}}\Hom_{S_n}(\tau, \CV_{c,0})\otimes\tau
\iso \CV_{c,0},
\eeq
where the direct sum on the left runs over the set
of (isomorphism classes of) simple $S_n$-representations $\tau$.

The  cohomology vanishing in Theorem 
\ref{equiv_intro} implies that the splitting vector bundle $\CV_{c,0}$
is {\em rigid}, i.e., we have
$\Ext^1(\CV_{c,0},\CV_{c,0})=0$. 
It follows that, for each $\tau$, the
 isotypic component $\Hom_{S_n}(\tau, \CV_{c,0})$
is also a rigid vector bundle.

The  tautological $GL_2$-action  on $\BA^2$ gives rise
to a natural  $GL_2$-action  on 
the punctual Hilbert scheme  $\hilbb_0$. We restrict attention
to the subgroup $\Gm\times\Gm\sset GL_2$, of diagonal matrices, and view
 $\hilbb_0$ as a projective $\Gm\times\Gm$-variety.
Thus, according to  Proposition \ref{vologod}, see \S\ref{volgd},
for each irredicible $S_n$-module $\tau$, the corresponding
 isotypic component $\Hom_{S_n}(\tau, \CV_{c,0})$
may  be equipped with  a 
 $\Gm\times\Gm$-equivariant structure. 
This gives, via the isomorphism \eqref{isotyp},
 a  $\Gm\times\Gm$-equivariant  structure
on the vector bundle $\CV_{c,0}$.
Further, the $\Gm\times\Gm$-equivariant   structure on   $\CV_{c,0}$
induces one  on the
 Azumaya algebra $\wh{\ah}_{c,0}=
\Bigl(\eend_{_{\oo_{{\wh{\op{Hilb}}}_\xi\tw}}}\CV_{c,\xi}\Bigr)^{\op{opp}}.$
Therefore, there is a continuous  $\Gm\times\Gm$-action
on the vector space $\wh{\hh}_{c,0}=\G(\widehat{\op{Hilb}}_0\tw,
\ah_c)$, of  global sections. 

Thus, we have defined an action of
the group  $\Gm\times\Gm$ on the
topological algebra  $\wh{\hh}_{c,0}$
by continuous algebra automorphisms.
Furthermore, the  $\Gm\times\Gm$-action on   $\wh{\hh}_{c,0}$
 fixes each element of the group $S_n
\sset \wh{\hh}_{c,0},$
since  the $\Gm\times\Gm$-equivariant   structure 
is compatible with the $S_n$-action, by construction.
Taking the direct sum of the weight spaces of the
  $\Gm\times\Gm$-action, we get
a $\Z^2$-graded  dense subalgebra  $\hh^\circ=\bigoplus_{k,l\in
  \Z}\hh_c^{k,l}
\sset \wh{\hh}_{c,0},$
with $S_n$-stable homogeneous components~$\hh_c^{k,l}.$

\section{Induction functor and  comparison with \cite{EG}}\label{comn}
In this subsection, we let  $\k$ be an arbitrary 
algebraically closed field,
either of characteristic zero or of characteristic $p$.

\subsection{} 
Let $P$ be a linear algebraic group with Lie algebra $\p$,
and $\DD$
an associative algebra equipped with a $P$-action 
by algebra automorphisms
 and with a $P$-equivariant algebra map
$\rho: \Up\to\DD$, as in Sect. \ref{quant_red}.
Recall our convention to write 
$\DD\cdot J$ instead of $\DD\cdot \rho(J)$.

\begin{defn}\label{ind} Given a two-sided ideal $J\sset\Up$,
let $\ind(\DD{\uparrow}J)\sset \DD$ denote the annihilator of
the left $\DD$-module $\DD/\DD\cdot J$. This is a
 two-sided ideal in $\DD$, called the  ideal {\em induced}
from $J$.
\end{defn}
  It follows from the definition 
that $\ind(\DD{\uparrow}J)$ is the maximal two-sided ideal of $\DD$
contained in the left ideal $\DD\cdot J$.

Assume next that $P$ is an algebraic subgroup
in another connected linear algebraic group $G$.
Set $\g:=\Lie G$, and let $\Ug$ be the corresponding enveloping algebra.
Thus, $\p\sset\g$ and $\Up\sset\Ug$.
Given a two-sided ideal $J\sset\Up$, as above,
we may form an induced ideal $\ind(\Ug\upi)\sset\Ug$.

Now, let $G$ act on an associative algebra
$\DD$, and let $\Ug\to \DD$ be a $G$-equivariant algebra map.
We consider the composite map
$\Up\into\Ug\to \DD$. Let  $J\sset \Up$ be a two-sided ideal
and $\ind(\Ug{\uparrow}J)\sset \Ug$ the corresponding
induced ideal. Since $J\sset \Ug\cdot J$ we 
have $\DD\cdot \ind(\Ug{\uparrow}J)\sset \DD\cdot J$. 
Hence, the
projection $\DD/\DD\cd\ind(\Ug{\uparrow}J)\onto \DD/\DD\cd J$ induces
an algebra map
\beq{theta}
 (\DD/\DD\cd\ind(\Ug{\uparrow}J))^G\too (\DD/\DD\cd J)^P.
\eeq

\subsection{}
We recall the setup of section \ref{calo}, so
 $V$ is an $n$-dimensional vector space over $\k$,
and we put $G=\GL(V)$ and $\g=\gl(V).$
We also fix $c\in\k$ and let $\chi_c:=c\cdot\Tr: \g\to\k$ be
the corresponding Lie algebra character  

     From now on, we fix a non-zero vector $\bv\in V$.
Let $P$ be a parabolic subgroup of $G$
formed by the maps $V\to V$ that preserve
the line $\k\bv$. Thus $G/P\cong\BP(V)$,
the $(n-1)$-dimensional projective space associated to $V$.
We put $\p:=\Lie P$, and write $\chi^\p_c:=\chi_c\big|_\p$
for the character $\chi_c$ restricted to the
subalgebra $\p\sset\g$.
Thus $\chi^\p_c\in\p^*$ is a $P$-fixed point for the coadjoint action of
$P$
on $\p^*$. 

We extend $\chi^\p_c$ to an
 associative algebra homomorphism $\chi^\p_c: \Up\to\k$,
and 
let $J_c:=\Ker(\Up\to\k),$  denote
the corresponding two-sided ideal generated by the elements\linebreak
${\{x-c\cdot\Tr(x)\}_{x\in\p},}$ see Definition \ref{IJ}.
Also, write ${\op{\mathsf{Ind}}}_c:=\ind(\Ug{\uparrow}J_c)$ 
for the two-sided ideal in $\Ug$
induced from $J_c$.
It is known that  ${\op{\mathsf{Ind}}}_c$ is
a primitive ideal in $\Ug$,
moreover, it is exactly the
primitive ideal considered in \cite{EG}.

The adjoint action of $G$  on $\g$ gives rise to
an associative algebra homomorphism
$\ad: \Ug \to \D(\g)$.
Thus, we may consider the homomorphism 
\eqref{theta} in the special case
$\DD:=\D(\g),$
$J:=J_c$. 
As usual, we abuse the notation and write
$\D(\g)\cdot{{\op{\mathsf{Ind}}}_c}$
instead of $\D(\g)\cdot\ad{{\op{\mathsf{Ind}}}_c}$
for the corresponding left ideal in $\D(\g).$

We propose the following

\begin{conj}
  For any $c\in\k$, the following canonical map 
is an algebra isomorphism: 
\beq{PG}
 \bigl(\D(\g)\big/\D(\g)\cdot{{\op{\mathsf{Ind}}}_c}\bigr)^G
\stackrel{\eqref{theta}}\tooo\bigl(\D(\g)\big/\D(\g)\cdot{J_c}\bigr)^{P}.
\eeq
\end{conj}

\subsection{}
In chapter 7 of [EG], the authors have constructed, for any $c\in\k$, an algebra
homomorphism, called the {\em deformed Harish-Chandra homomorphism}:
\beq{rad}
\Phi_c:\ \bigl(\D(\grs)\big/\D(\grs)\cdot{{\op{\mathsf{Ind}}}_c}\bigr)^G
\too \D(\hreg)^W.
\eeq
such that one has $\Phi_c(\Delta_\g)=\cm_c,$
the Calogero-Moser operator.

The reader should be warned that the map referred to  as the
deformed Harish-Chandra homomorphism in [EG]
was actually a map
$
\D(\grs)^G\too \D(\hreg)^W.
$
However, it has been shown in \cite{EG} that the latter
map vanishes on the two-sided ideal 
$$\bigl(\D(\grs)\cdot{{\op{\mathsf{Ind}}}_c}\bigr)^G=
\bigl(\D(\grs)\cdot{{\op{\mathsf{Ind}}}_c}\bigr)\;
\cap\;\D(\grs)^G,
$$
hence, descends to a well-defined homomorphism
\beq{rad0}
\bigl(\D(\grs)\big/\D(\grs)\cdot{{\op{\mathsf{Ind}}}_c}\bigr)^G=
\D(\grs)^G\big/\bigl(\D(\grs)\cdot{{\op{\mathsf{Ind}}}_c}\bigr)^G
\too \D(\hreg)^W,
\eeq
where the equality on the left   exploits 
 semisimplicity of the $\ad G$-action on $\D(\grs)$.

The  construction of \cite{EG} can be also carried out
over a field $\k$ of characteristic $p$.
In that case, the adjoint $G$-action on $\D(\grs)$ is not semisimple,
so the equality on the left of \eqref{rad0} does not hold,
in general. A more careful analysis of the construction
of  \cite{EG}, similar to that of Section \ref{hc_sec} of the present paper,
shows
 that it actually produces, without any
semisimplicity assumption, 
a homomorphism of the form
\eqref{rad}.

\subsection{} Recall the open subset $U\sset \grs\times\VC$
formed by all pairs $(x,v)$ such that $v$ is a cyclic vector for $x$, see
Definition \ref{eer}.
Let $\gc\sset \grs$ be the set
of all elements $x\in\g=\gl(V)$ such that
our fixed vector
$\bv\in V$ is a cyclic vector for $x$. Clearly, $\gc$ is
an $\Ad P$-stable Zariski open dense subset of $\g$.
Recall also the two-sided ideal $I_c\sset\Ug$ introduced 
in Sect. \ref{hc_sec}, and observe that
${\op{\mathsf{Ind}}}_c\sset I_c$.

The following result provides a relation between the
Harish-Chandra homomorphism $\Psi_c$ of Proposition
\ref{prop}, and the  homomorphism
\eqref{rad}.

\begin{prop}\label{compare} There is a natural algebra
map $F_c: \bigl(\D(\gc)\big/\D(\gc)\cdot{J_c}\bigr)^{P}
\too{\G(U,\D_U/\D_U\cdot I_c)^G}$ making the
following diagram commute
$$
\xymatrix{
{\bigl(\D(\grs)\big/\D(\grs)\cdot{J_c}\bigr)^{P}}\ar[r]^<>(0.5){\text{res}}
&{\bigl(\D(\gc)\big/\D(\gc)\cdot{J_c}\bigr)^{P}}
\ar[r]^<>(0.5){F_c}
&\G(U,\,\D_U/\D_U\cd I_c)^G\ar[d]_<>(0.5){\Psi_c}^<>(0.5){\text{Prop. \ref{prop}}}\\
{\bigl(\D(\grs)/\D(\grs)\cdot{{\op{\mathsf{Ind}}}_c}\bigr)^G}
\ar[u]_<>(0.5){\eqref{theta}}
\ar@{->>}[rr]^<>(0.5){\Phi_c}_<>(0.5){\eqref{rad}}&&{\;\D(\hreg)^W\;} 
.}
$$
\end{prop}
\begin{proof}[Sketch of Proof.]
The assignment $(x,v)\mto \k v$ clearly
gives a  $G$-equivariant
fibration $f: U\map \BP(V)$.
The fiber of this map over the class of the line $\k\bv\in\BP(V)$
 is the set
$U_\bv=\{(x,t\cdot\bv)\mid x\in\gc, t\in\k^\times\}$.
Thus, we have a $G$-equivariant isomorphism
$U\cong G\times_P U_\bv$.
Further,
it is clear that the map $(x,t)\mapsto (x,t\cdot\bv)$ gives a
$P$-equivariant isomorphism
$\gc\times \Gm \iso U_\bv$. Thus,
 we obtain the following
 diagram
\beq{big_map}
\xymatrix{
G\times \gc\times \Gm\ar@{=}[r]&G\times  U_\bv
\ar@{->>}[rrrr]_<>(0.5){P\text{-bundle}}^<>(0.5){(g,x)\mto
g\times_{_P}x}\ar[d]^<>(0.5){\pr_{_G}}&&&&
{\;G\times_\Parab U_\bv\;}\ar@{=}[r]
\ar[d]^<>(0.5){\pr}
&
U\ar[d]^<>(0.5){f}
\\
&G\ar[rrrr]^<>(0.5){g\mapsto
gP}
&&&&G/\Parab\ar@{=}[r]&
\BP(V).
}
\eeq

Let $X:=G\times  U_\bv\cong G\times \gc\times \Gm.$
The assignment $u\mto 1\otimes u\otimes 1$
gives an
algebra map 
$$ \D(\gc)\too \D(G)\otimes\D(\gc)\otimes\D(\Gm)
\iso \D(G\times  U_\bv)=\D(X).
$$
We compose this map with
isomorphism \eqref{XA} applied to $X$,
$\,Y:=G\times_\Parab U_\bv$, and to the $P$-bundle map
 in  the top row of diagram \eqref{big_map}.
This way, we get a chain of algebra morphisms
\beq{tau}
F'_c:\ \D(\gc)^P \too \D(X)^P\too \D(X\,\D_X/\D_X\cd J_c)^P 
\too
\D(U, \chi_c).
\eeq

One can show,
using the equality
$\ind(\Ug{\uparrow}J_c) = \cap_{g\in G}\, \Ad g(\Ug\cdot J_c),
$ that a  suitably refined version of  the above
construction  produces a well-defined algebra map
$$
F_c'':\  \G(\gc,\, \D_\gc/\D_\gc\cd J_c\bigr)^P
\too \G(U, \D_U/\D_U\cd {\op{\mathsf{Ind}}}_c)^G.
$$
Furthermore, one verifies that composing the map $F''_c$
with the natural restriction map
$\G(\grs,\, \D_\grs/\D_\grs\cdot J_c\bigr)^P\too
 \G(\gc,\, \D_\gc/\D_\gc\cdot J_c\bigr)^P$ one obtains
a homomorphism $F_c$ that makes
the diagram of the Proposition commute.
\end{proof}
\bigskip
\pagebreak[3]

\section{Appendix: The $p$-center of symplectic reflection algebras}

\centerline{\bf {by Pavel Etingof}}
\vskip 6pt

\subsection{} Let $\k$ be an algebraically closed field of characteristic $p>0$,
$V$ a finite-dimensional symplectic vector space over $\k$,
and $\G\sset Sp(V)$ a finite subgroup.
We assume that $p$ is odd and prime to $|\Gamma|$.

Write $S\sset\G$ for the set of symplectic reflections in $\G$.
Let $\hh_{t,c}(V,\G)$ be the symplectic reflection algebra
with parameters $t$ and $ c\in\k[S]^\G$, as defined in \cite{EG}.
For $t=c=0$, we have $\hh_{0,0}\simeq SV\#\G$,
a cross-product of $SV$ with the group $\G$,
to be denoted $\ohh$. 

In general,
the algebra $\hh_{t,c}$ comes equipped with an increasing filtration
such that $\gr\hh_{t,c}\simeq \ohh=SV\#\G.$
Let ${\mathsf{Z}}_{t,c}$ be the center of $\hh_{t,c}$,
equipped with the induced filtration. 

 The main result of this appendix is the following theorem.

\begin{theorem}\label{pasha}
There is a graded algebra isomorphism $\gr ({\mathsf{Z}}_{1,c})=((SV)^p)^\Gamma.$
\end{theorem}

\begin{corollary}[Satake isomorphism]
The map ${\mathsf{Z}}_{1,c}\to\ehe,\,z\mto \e\cdot z$
gives an isomorphism between the centers
of the algebras $\hh_{1,c}$ and $\e\hh_{1,c}\e$, respectively.
\end{corollary}
\begin{proof}[Proof of Corollary.] The associated graded map is an isomorphism.
\end{proof}

\begin{rem} One may consider the field $\k$ as the residue
class field in $W(\k)$, the ring of Witt vectors of $\k$. 
The algebra  $\hh_{1,c}$ may thus be regarded as
a specialization of a  $W(\k)$-algebra.
Applying the well-known,
see \cite{Ha}, Hayashi construction
 to this situation one obtains
a natural Poisson bracket on the $p$-center
${\mathsf{Z}}_{1,c}\sset \hh_{1,c}$.
\erem

\subsection{}
To prove  Theorem \ref{pasha} we need the following result.

\begin{prop}\label{hom} The Hochschild 
cohomology of the algebra $\ohh$ is given
by the formula
$$
\HH^m(\ohh)=(\oplus_{g\in \Gamma: {\rm codim}V^g\le m}
\Omega^{m-{\rm codim}V^g}(V^g))^\Gamma.
$$
(here $\Omega^j$ denotes the space of differential forms 
of rank $j$) Moreover, elements $z\in {\mathsf{Z}}_{0,0}=\HH^0(\ohh)$ 
act on this cohomology by multiplying differential forms on $V^g$ 
by the restriction of $z$ to $V^g$.
\end{prop}

\begin{proof} The first statement is 
straightforward by using Koszul resolutions, similarly to
\cite{AFLS}. The second statement is easy.
\end{proof} 
\subsection{Proof of Theorem \ref{pasha}.}
Since $\gr \hh_{t,c}=\ohh$, we have 
the Brylinski spectral sequence\cite{Br}. This spectral sequence has
$$
E_1^{r,q}=\HH^{r+q}(\ohh),\quad \text{and}\quad
E_\infty^{r,q}=\gr (\HH^{r+q}(\hh_{t,c})).
$$
Since the algebra $\hh_{t,c}$ is $\Z/2\Z$-graded, 
the differentials $d_j$ with odd subscripts $j$ in this spectral sequence
vanish automatically. Therefore we will abuse notation 
by writing $d_i,E_i^{r,q}$ instead of $d_{2i},E_{2i}^{r,q}$. 

Let us now consider the differential $d_1$. 
Obviously, we have 
$d_1=td_1^{(1)}+\sum_{s\in S/\Gamma}
c_sd_1^{(s)}$, where $d_1^{(1)}$ is the differential corresponding to
$c=0$, $t=1$. It was checked by Brylinski that $d_1^{(1)}$ is
the De Rham differential. 

By a result of \cite{EG} (which generalizes in a straightforward manner to 
positive characteristic), we have
$\gr {\mathsf{Z}}_{0,c}={\mathsf{Z}}_{0,0}$. Thus 
the differentials $d_1^{(s)}$ are zero in cohomological degree
zero, so they are morphisms of modules 
over ${\mathsf{Z}}_{0,0}$. Thus in cohomological degree 1, these differentials must 
land in twisted (torsion) components (as generically 
on $M_c:=Spec(Z_{0,c})$ the algebra $\hh_{0,c}$ is Azumaya). This implies that 
$E_2^{-q,q}=((SV)^p)^\Gamma$, and $E_2^{-q+1,q}$ is a submodule 
(over $((SV)^p)^\Gamma$) of $\Omega_{(p)}^1(V)^\Gamma$, 
where $\Omega^1_{(p)}(V)$ denotes
the first cohomology of the De Rham complex of $V$. 

We will now show that all the higher
differentials $d_m$, $m\ge 2$, vanish in cohomological degree zero, 
and hence $E_\infty^{-q,q}=((SV)^p)^\Gamma$, as desired. 
Assume the contrary, and let $m\ge 2$ be the smallest number 
such that $d_m$ does not vanish in cohomological degree zero. 
Then we can view $d_m$ as a derivation 
$d_m: ((SV)^p)^\Gamma\to \Omega_{(p)}^1(V)^\Gamma$. 

Recall now (\cite{K}, Theorem 7.2)
that if $X$ is a smooth affine variety over $\k$ then 
the cohomology of the De Rham complex of $X$ twisted by 
the Frobenius map can be identified with the module of
differential forms on $X$ (as a graded ${\mathcal O}_X$-module) 
via a map $C^{-1}: \Omega^\bullet(X)\to H^\bullet(\Fr\idot
\Omega^\bullet(X))$,
called the Cartier operator. This map is defined
by the formulas $C^{-1}(a)=a^p$, $C^{-1}(da)=a^{p-1}da$ for functions 
$a$ on $X$. 

Let $d_m'=Cd_mC^{-1}$. Then $d_m': (SV)^\Gamma\to \Omega^1(V)^\Gamma$ is
a derivation. Thus, $d_m'$ gives rise to a $\Gamma$-equivariant
regular function $f$ on $V_{reg}$ with values in $V\otimes V^*$,
where $V_{reg}$ is the set of points of $V$ which have
the trivial stabilizer in $\Gamma$. Since the 
complement of $V_{reg}$ has codimension 2, the function $f$
extends to $V$ and defines an element of $SV\otimes V\otimes
V^*$. Hence $d_m'$ is obtained by restricting a map 
$d_m': SV\to \Omega^1(V)$ to $\Gamma$-invariants. 
Therefore, $d_m$ is obtained
by restricting a map $d_m: (SV)^p\to \Omega_{(p)}^1(V)$ 
to $\Gamma$-invariants.

Now observe that the map $d_m$ must have degree $-2m+2<0$. 
On the other hand, the generators of $(SV)^p$ sit in degree $p$,
while the lowest degree in $\Omega^1_{(p)}(V)$ is also equal to
$p$. This means that $d_m=0$, which is a contradiction. 
The theorem is proved. 
\qed

{\footnotesize{

}}

{\small{
\noindent
{\bf R.B.}:  Department of Mathematics,  MIT,
77 Mass. Ave, Cambridge, MA 02139, USA;\\
\hphantom{x}\quad\quad\enspace {\tt bezrukav@math.mit.edu}
\smallskip

\noindent
{\bf M.F.}: Independent University of Moscow, 11 Bolshoy Vlasyevskiy per.,
119002 Moscow, Russia;\\
\hphantom{x}\quad\quad\enspace {\tt fnklberg@mccme.ru}
\smallskip

\noindent
{\bf V.G.}: Department of Mathematics, University of Chicago,
Chicago, IL
60637, USA;\\
\hphantom{x}\quad\quad\enspace {\tt ginzburg@math.uchicago.edu}}
\smallskip

\noindent
{\bf P.E.}: Department of Mathematics,  MIT,
77 Mass. Ave, Cambridge, MA 02139, USA;\\
\hphantom{x}\quad\quad\enspace {\tt etingof@math.mit.edu}
\smallskip

\begin{thebibliography}{APK}

\bibitem[AFLS]{AFLS} J. Alev, M.A. Farinati, T. Lambre, and A.L. Solotar:
{\it Homologie des invariants d'une alg\`ebre de Weyl sous l'action d'un
groupe fini.} J. of Algebra {\bf 232} (2000), 564--577.



\bibitem[BB]{BB} A. Beilinson,  J. Bernstein,
{\it A proof of Jantzen conjectures.}
 I. M. Gelfand Seminar, 1--50, Adv. Soviet Math., {\bf 16}, Part 1,
Amer. Math. Soc., Providence, RI,
1993.

\bibitem[BGS]{BGS} A. Beilinson, V. Ginzburg, W. Soergel,
     {\it Koszul duality patterns in Representation Theory}, J. Amer. Math. Soc.
{\bf 9} (1996), 473--527.


                         

\bibitem[BEG]{BEG}
Yu. Berest, P. Etingof, V. Ginzburg,
\textit{Cherednik algebras and differential operators on quasi-invariants}.
 Duke Math. J. {\bf  118}  (2003),  279--337.
{\tt{arXiv:math.QA/0111005}}.

\bibitem[BEG2]{BEG2}
Yu. Berest, P. Etingof, V. Ginzburg,
\textit{Finite-dimensional representations of rational Cherednik
algebras.}
  Int. Math. Res. Not.  2003,  no. 19, 1053--1088. 
 

\bibitem[BK]{bk} R.~Bezrukavnikov, D.~Kaledin,
{\em MacKay equivalence for symplectic quotient singularities},
 Proc. of the Steklov Inst. of Math. {\bf 246} (2004), 13--33 (in Russian).
 {\tt{arXiv:math.AG/0401002}}.


\bibitem[BMR]{BMR} R. Bezrukavnikov, I.  Mirkovi\'c, D. Rumynin,
{\em Localization of modules for a semisimple Lie algebra in prime
characteristic},
Preprint {\tt{arXiv:math.RT/0205144}}.

\bibitem[BKR]{BKR} T. Bridgeland, A.  King, M.  Reid,
{\em  The McKay correspondence as an equivalence of derived categories.}
 J. Amer. Math. Soc.  {\bf 14}  (2001),  535--554.

\bibitem[Br]{Br} J.-L. Brylinski: {\it
A differential complex for Poisson manifolds.}
 J. Diff. Geom. 28 (1988), 93--114.

\bibitem[CG]{CG} N. Chriss, V. Ginzburg, {\it
Representation theory and complex geometry.} 
Birkh\"auser Boston, 1997.

\bibitem[Ch]{Ch} I. Cherednik, {\it Double affine Hecke
 algebras,
Knizhnik-Zamolodchikov equations, and Macdonald operators },
IMRN (Duke math.J.) {\bf 9} (1992), p.171-180.







\bibitem[DO]{DO} C. Dunkl, E. Opdam,
{\it Dunkl operators for complex reflection groups.}
Proc. London Math. Soc.  {\bf  86}  (2003), 70--108.
[{\tt{arXiv:math.RT/0108185}}].

\bibitem[EV]{EV} H. Esnault, E. Viehweg, 
\textsl{Lectures on vanishing theorems.}
DMV Seminar, {\bf 20}.
Birkha\"user Verlag, Basel, 1992.

\bibitem[EG]{EG} P. Etingof, V. Ginzburg,
{\it Symplectic reflection algebras, 
Calogero-Moser space, and deformed Harish-Chandra
   homomorphism.} 
Invent. Math. {\bf 147} (2002), 243-348, [{\tt{arXiv:math.AG/0011114}}].
\bibitem[Ga]{Ga} O. Gabber, {\em
 Some theorems on Azumaya algebras.  The Brauer group}
 (Sem., Les Plans-sur-Bex, 1980),  pp. 129--209, 
Lecture Notes in Math., {\bf 844}, Springer, Berlin-New York, 1981.
\bibitem[GIT]{GIT} D. Mumford, J. Fogarty, F. Kirwan, {\em
Geometric invariant theory.}
Third edition. Ergebnisse der Mathematik und ihrer Grenzgebiete
{\bf 34}
Springer-Verlag, Berlin, 1994. 

\bibitem[GG]{GG} W. L.
Gan, V. Ginzburg, {\em
Quantization of Slodowy slices.}
Int. Math. Res. Not. 2002, no. 5, 243--255.


\bibitem[GS]{g} I.~Gordon, T. Stafford,
{\em Rational Cherednik algebras and Hilbert schemes
of points.} Preprint 2004, {\tt{arXiv:math.RA/0407516}}
(to appear in Adv. Math.).

 \bibitem[Gr]{Gr} A. Grothendieck, SGA 1, Lecture Notes in Mathematics,
{\bf 224} (1971).

\bibitem[J]{J} N.  Jacobson, \textsl{Lie algebras.}
 Interscience Tracts in Pure and Applied Mathematics, {\bf  10}
Interscience Publishers (a division of John Wiley \& Sons), New
York-London 1962.

\bibitem[Ja]{Ja} J.C. Jantzen,  {\em
 Representations of Lie algebras in prime characteristic.}
 Notes by Iain Gordon. NATO Adv. Sci. Inst. Ser. C Math. Phys. Sci.,
{\bf 514},  Representation theories and algebraic geometry (Montreal, PQ,
1997),  185--235, Kluwer Acad. Publ., Dordrecht, 1998.


\bibitem[H]{H} M. Haiman, \textit{Vanishing theorems and character
formulas for the Hilbert scheme of points in the plane}, 
 Invent. Mathem. {\bf 149} (2002), 371-408,
[{\tt{arXiv:math.AG/0201148}}].

\bibitem[Ha]{Ha} T. Hayashi, {\it Sugawara operators and Kac-Kazhdan
conjecture.} Invent. Math. {\bf  94} (1988), 13--52.



\bibitem[Har]{Har} R. Hartshorne, \textsl{Algebraic geometry.}
 Graduate Texts in Mathematics, {\bf 52}. Springer-Verlag, New York-Heidelberg, 1977.



\bibitem[K]{K}
N. Katz, {\em Nilpotent connections and the monodromy theorem,} Publ. Math.
IHES {\bf 39}(1970).

\bibitem[KKS]{KKS} D. Kazhdan, B. Kostant, and S. Sternberg,
{\it Hamiltonian
group actions and dynamical systems of Calogero type,} Comm. Pure Apll. Math.,
{\bf 31}(1978),  481-507.



\bibitem[La]{La} F.  Latour, {\em 
Representations of rational Cherednik algebras of rank 1 in positive
characteristic},
Preprint {\tt{arXiv:math.RT/0307390}}.





\bibitem[Na1]{Na1}  H.  Nakajima, {\it 
Lectures on Hilbert schemes of points on surfaces.}
 University Lecture Series, \textbf{18}. American Mathematical Society, Providence, RI, 1999.


\bibitem[Na2]{Na2}  H.  Nakajima, {\it 
 Quiver varieties and Kac-Moody algebras.}
Duke Math. J. \textbf{91} (1998),  515--560.
 
\bibitem[Ob]{Ob} A. Oblomkov, {\em
Double affine Hecke algebras and Calogero-Moser spaces.}
 {\tt{arXiv:math.RT/0303190}}.

\bibitem[Pr]{Pr} A. Premet, {\em
 Special transverse slices and their enveloping algebras.}
 Adv. Math.  {\bf  170}  (2002),   1--55.
 

\bibitem[PS]{PS} A. Premet, S.  Skryabin,
{\em Representations of restricted Lie algebras and families of 
associative $\scr L$-algebras.} J. Reine Angew. Math.   {\bf  507}  (1999), 189--218.



\bibitem[Q]{Q} D. Quillen, {\it 
On the endormorphism ring of a simple module over an enveloping
algebra.}
  Proc. Amer. Math. Soc. {\bf  21} 1969 171--172. 


\bibitem[Wi]{Wi} G. Wilson, {\it Collisions of Calogero-Moser
particles and an adelic Grassmannian,} Inv. Math.
{\bf 133} (1998), 1-41.

\end{thebibliography}
\end{document}